\documentclass[a4paper,12pt,reqno]{amsart}

\usepackage{amsmath,amssymb,amsthm,bm}
\usepackage{caption,graphicx}
\usepackage[english]{babel}
\usepackage{amscd}
\usepackage{amsgen}
\usepackage[final]{epsfig}
\usepackage{latexsym}
\usepackage{amsfonts}
\usepackage{url}
\usepackage{slashed}
\usepackage[all]{xy}
\usepackage{here}
\usepackage{comment}

\newtheorem{thm}{Theorem}
\newtheorem{lemm}{Lemma}
\newtheorem{prop}{Proposition}
\newtheorem{cor}{Corollary}

\newtheorem{defn}{Definition}

\allowdisplaybreaks

\begin{document}

\title{Fundamental groups, 3-braids, and effective estimates
of invariants}

\author{Burglind J\"oricke \\ {\ ´´}}

\address{Max-Planck-Institut f\"ur Mathematik, Vivatsgasse 7,
\\53111 Bonn\\ Germany     }

\email{joericke@googlemail.com}




\begin{abstract}
We define invariants of braids rather than invariants of conjugacy
classes of braids.  For any pure $3$-braid we give effective upper and
lower bounds for
these invariants. This is done in terms of a natural syllable decomposition
of the word representing the image of the braid in the
braid group modulo its center. The bounds differ by a multiplicative
constant not depending on the word. Respective bounds are given for all $3$-braids.
We also obtain effective upper and lower bounds for the entropy of
pure $3$-braids in these terms.  The proof leads to the study of the
extremal length (in the sense of Ahlfors) of classes of curves representing elements of the
fundamental group of the twice punctured complex plane.

\end{abstract}

\maketitle

\centerline \today

\section{Statement of the results}\label{intro}

\bigskip

In this paper we define invariants of $n$-braids (rather than
invariants of conjugacy classes of braids.) Note, that a popular
invariant of braids, the entropy, is a conjugacy invariant, i.e. it
does not distinguish different elements of a conjugacy class of braids.
The present invariants may be compared from a conceptional point of
view with the conformal module of conjugacy classes of braids which is
inverse proportional to the entropy (see \cite{Jo1}, \cite{Jo2}).

For
pure  $3$-braids we give effective upper and lower bounds for the invariants. This is done in
terms of a natural syllable decomposition of
the word representing the image of the braid in $\mathcal{B}_3 \diagup
\mathcal{Z}_3$, the group of $3$-braids modulo its center. The bounds
differ by
a universal multiplicative constant.

The estimates for the $3$-braid invariants also give bounds for the
entropy of arbitrary pure $3$-braids in the mentioned terms, and also
provide estimates of the
invariants for arbitrary (maybe, not pure) $3$-braids.

Notice that there is an algorithm
which detects in principle whether a braid is a pseudo-Anosov braid and
in this
case it computes the entropy (\cite{BH}). This approach uses train
tracks. 
For $3$-braids which can be represented by short words the entropy is
known and can be calculated explicitly without using train
tracks. As far as we know
despite these results explicit bounds for the entropy of arbitrary
$3$-braids in terms of representing words are new.

The braid invariants are defined as follows. For a subset $A$ of the
complex plane
$\mathbb{C}$ we consider the
configuration space $C_n (A) = \{(z_1 , \ldots , z_n) \in A^n
: z_i
\ne z_j$ for $i \ne j\}$ of $n$ particles moving along $A$
without collision. Denote by ${\mathcal S}_n$ the symmetric group.
Each
permutation in ${\mathcal S}_n$ acts on $C_n (A)$ by permuting
the
coordinates. Consider the quotient $C_n (A) \diagup {\mathcal
S}_n$.
The quotient $C_n (A) \diagup {\mathcal S}_n$ is
called
the symmetrized configuration space related to $A$.
The natural projection $C_n (\mathbb{C}) \to  C_n (\mathbb{C}) \diagup {\mathcal S}_n$ is denoted by $\mathcal{P}_{\mathcal{S}_n}$.

Choose a base point $E_n \in C_n (\mathbb {C}) \diagup {\mathcal S}_n$.
Regard braids on $n$ strands  ($n$-braids for short)
with base point $E_n$ as homotopy classes of loops with
base point $E_n$ in the symmetrized configuration space,
equivalently, as
elements of the fundamental group $\pi_1(C_n (\mathbb {C})
\diagup
{\mathcal S}_n, E_n)$ of the symmetrized configuration space
with
base point $E_n$.

The totally real subspace $\mathcal{E}^n_{tr}\stackrel{def}{=}\,C_n (\mathbb {R})
\diagup{\mathcal S}_n\,$ of $\,C_n (\mathbb {C}) \diagup {\mathcal
S}_n\,$ is connected and simply connected. Indeed, the totally real
subspace
$C_n (\mathbb {R})$ of $C_n (\mathbb {C})$ is the union of the
connected components $\,\{(x_1, \ldots , x_n) \in
\mathbb{R}^n: \,
x_{\sigma(1)} < x_{\sigma(2)} < \ldots < x_{\sigma(n)}\}\,$
over all
permutations $\sigma \in {\mathcal S}_n$. Thus $C_n (\mathbb
{R})$
is invariant under the action of ${\mathcal S}_n$ and the
quotient
is homeomorphic to $\{(x_1, \ldots , x_n) \in \mathbb{R}^n:
x_1 <
x_2 < \ldots < x_n\}\,$. Hence the claim.

The fundamental group $\pi_1(\,C_n (\mathbb {C}) \diagup
{\mathcal
S}_n\,,\; E_n\,)$ is isomorphic to the relative fundamental
group
$\pi_1(\,C_n (\mathbb {C}) \diagup {\mathcal S}_n\,,\; C_n
(\mathbb
{R}) \diagup {\mathcal S}_n \,)$. The elements of the latter
group
are homotopy classes of arcs in $\,C_n (\mathbb {C}) \diagup
{\mathcal S}_n\,$ with endpoints in $\,C_n (\mathbb {R})
\diagup
{\mathcal S}_n\,$.

The isomorphism between the two groups is obtained as follows.
Since
a change of the base point leads to an isomorphism of the
fundamental
group, we may assume that $E_n$ is contained
in the
totally real subspace $\, C_n (\mathbb {R}) \diagup {\mathcal
S}_n\,$. Each element of $\pi_1(\,C_n (\mathbb {C}) \diagup
{\mathcal S}_n\,,\; E_n\,)$ is a subset of an element of
$\pi_1(\,C_n (\mathbb {C}) \diagup {\mathcal S}_n\,,\; C_n
(\mathbb
{R}) \diagup {\mathcal S}_n \,)$. Vice versa, since $\,C_n
(\mathbb
{R}) \diagup {\mathcal S}_n\,$ is connected and $E_n$ is
contained in  $\,C_n (\mathbb{R}) \diagup {\mathcal S}_n\,$  each class
in $\pi_1(\,C_n (\mathbb {C}) \diagup
{\mathcal S}_n\,,\; C_n (\mathbb {R}) \diagup {\mathcal S}_n\,
)$
contains a class in $\pi_1(\,C_n (\mathbb {C}) \diagup
{\mathcal
S}_n\,,\; E_n\,)$. Since $\,C_n (\mathbb {R}) \diagup
{\mathcal
S}_n\,$ is simply connected, each class in $\pi_1(\,C_n
(\mathbb
{C}) \diagup {\mathcal S}_n\,,\; C_n (\mathbb {R}) \diagup
{\mathcal
S}_n \,)$ contains no more than one class of $\pi_1(\,C_n
(\mathbb
{C}) \diagup {\mathcal S}_n\,,\; E_n\,)$. Indeed, if two loops
in
$\,C_n (\mathbb {C}) \diagup {\mathcal S}_n\,$ with base point
$E_n$
are homotopic as loops in $\,C_n (\mathbb {C}) \diagup
{\mathcal
S}_n\,$ with varying base point in $\,C_n (\mathbb {R})
\diagup
{\mathcal S}_n\,$ then they are homotopic as loops in $\,C_n
(\mathbb {C}) \diagup {\mathcal S}_n\,$ with fixed base point
$E_n$.

Let $R$ be an open rectangle in the complex plane $\mathbb{C}$. Unless said otherwise
the considered rectangles will always have sides
parallel to the coordinate axes. Denote the length of the horizontal
sides of $R$
by $\sf b$ and the length of the vertical sides by
$\sf a$.
(For instance, we may consider $R= \{z=x+iy:0<x<\sf b,\,
0<y<\sf a\,\}$.) The conformal module of the rectangle $R$ introduced
by
Ahlfors \cite{A1} equals $m(R)= \frac{\sf b}{\sf a}$.
The extremal length of $R$ equals $\lambda(R)=\frac{\sf a}{\sf
b}$,
which is the inverse of the conformal module.

Let $b \in \mathcal{B}_n$ be a braid.
Denote its image in the relative fundamental group\\
$\pi_1(\,C_n
(\mathbb {C}) \diagup {\mathcal S}_n\,,\; C_n (\mathbb {R})
\diagup
{\mathcal S}_n\, )$ by $b_{tr}$. For a rectangle $R$ as above let $f:R
\to \,C_n
(\mathbb {C}) \diagup {\mathcal S}_n\,$ be a mapping which
admits a continuous extension to the closure $\bar R$ (denoted again by $f$) which
takes the (open) horizontal sides into $\,C_n (\mathbb {R})
\diagup
{\mathcal S}_n\,$. We say that the mapping represents $b_{tr}$
if for each maximal vertical line segment contained in $R$
(i.e. $R$
intersected with a vertical line in $\mathbb{C}$) the
restriction of
$f$ to the closure of the line segment represents $b_{tr}$.

We are now ready to define for any braid its  extremal length  with
totally real boundary values (and the conformal module  with totally
real boundary values, respectively).
\begin{defn}\label{def1} Let $b \in \mathcal{B}_n$ be an
$n$-braid. The extremal length $\Lambda_{tr}(b)$ with totally
real horizontal
boundary values is defined as
\begin{align}
\Lambda_{tr}(b)=& \inf \{\lambda(R): R\, \mbox{ a rectangle
which
admits a holomorphic mapping to} \nonumber \\
&C_n (\mathbb {C}) \diagup {\mathcal S}_n \,\mbox{ that
represents}\; b_{tr}\}\,.\nonumber
\end{align}
The conformal module $\mathcal{M}_{tr}(b)$ of $b$ with totally
real horizontal
boundary values, respectively, is defined as
\begin{align}
\mathcal{M}_{tr}(b)= &\sup \{m(R): R\, \mbox{ a rectangle which
admits a holomorphic mapping to}\; \nonumber \\
& C_n (\mathbb {C}) \diagup {\mathcal S}_n \,\mbox{ that
represents}
\; b_{tr}\}\,.\nonumber
\end{align}
\end{defn}
Note that the two invariants are
inverse to each other. It is more convenient to work with the extremal
length and we will mostly speak about the extremal length rather than
about the conformal module.

The choice of the boundary values is motivated by real algebraic
geometry.
Following Arnold \cite{Ar} a point in the symmetrized
configuration
space  $C_n (\mathbb {C}) \diagup {\mathcal S}_n$ can be
considered
as unordered $n$-tuple of pairwise distinct complex numbers and
can
be identified with the monic polynomial with these zeros.
Parametrize the space $\mathfrak{P}_n$ of monic polynomials of
degree $n$ without multiple zeros by their coefficients. We
obtain a
biholomorphic mapping between  $C_n (\mathbb {C}) \diagup
{\mathcal S}_n$ and
$\mathbb{C}^n \setminus \{\textsf{D}_n=0\}$, where
$\textsf{D}_n$ is
the discriminant.

The zero set of a polynomial with real coefficients consists
of a
set of real points and a set of pairs of complex conjugate
points. A
polynomial in $\mathfrak{P}_n$ all whose zeros are real can be
identified with a point in $C_n (\mathbb {R}) \diagup
{\mathcal
S}_n$. Thus, totally real boundary values correspond to one of
the
cases arising in real algebraic geometry. The set of
polynomials in
$\mathfrak{P}_n$ with real coefficients and $\ell \ge 2$ pairs
of
complex conjugate zeros is not simply connected. We will not
develop
the case of such boundary values in this paper.

However, for 3-braids it is convenient to consider the case of
boundary values corresponding to one real root and a pair of
complex
conjugate roots. The set $\mathcal{E}^3_{pb}$ of polynomials in
$\mathfrak{P}_3$
with one real zero and a pair of complex conjugate zeros is
simply
connected. ("$pb$" stands for perpendicular bisector. In fact,
the
real root lies on the perpendicular bisector of the line
segment
joining the two complex conjugate
roots.)
The braid group $\mathcal{B}_3$ is isomorphic to $\pi_1(C_3
(\mathbb
{C}) \diagup {\mathcal S}_3, \mathcal{E}^3_{pb})$. For a braid $b \in
\mathcal{B}_3$ we denote by $b_{pb}$ its image in $\pi_1(C_3
(\mathbb {C}) \diagup {\mathcal S}_3, \mathcal{E}^3_{pb})$ under this
isomorphism. The convention that a continuous mapping from a rectangle
to symmetrized configuration space represents $b_{pb}$ is made in the same way as the respective convention for $b_{tr}$.

We have the following definition.
\begin{defn}\label{def2} Let $b \in \mathcal{B}^3$ be a
$3$-braid.
The extremal length $\Lambda_{pb}(b)$ of $b$ with
perpendicular
bisector boundary values is defined as
\begin{align}
\Lambda_{pb}(b)= & \inf \{\lambda(R):\; R \, \mbox{a rectangle
which
admits a holomorphic mapping to}  \nonumber \\
& C_3 (\mathbb {C}) \diagup {\mathcal S}_3 \;\mbox{ that
represents}\;
b_{pb}\}. \nonumber
\end{align}
\end{defn}
Recall the definition of the extremal length of closed braids. Closed braids can be identified with conjugacy classes of braids or with free
homotopy classes of loops in $C_n(\mathbb{C}) \diagup \mathcal{S}_n$. (See \cite{Jo1}, \cite{Jo2}.)
The extremal length (conformal module, respectively) of a
conjugacy class of braids is defined as follows. We say that a
continuous mapping $f$ of an annulus $A= \{z \in \mathbb{C}:
\, r<|z|<R\},\;$ $0\leq r < R \leq \infty,\;$ into
$C_n(\mathbb{C}) \diagup \mathcal{S}_n$ represents a conjugacy class
$\hat b$ of $n$-braids if for
some (and hence for any) circle $\,\{|z|=\rho \} \subset A\,$
the
loop $\,f:\{|z|=\rho \} \rightarrow  C_n(\mathbb{C}) \diagup
\mathcal{S}_n
\,$
represents the conjugacy class $\,\hat b$. According to Ahlfors' definition the
conformal module of an annulus $\,A= \{z \in
\mathbb{C}:\; r<|z|<R\}\,\,$ equals $\,m(A)=
\frac{1}{2\pi}\,\log(\frac{R}{r})\,, $ and the extremal length $\lambda(A)$  equals $\frac{1}{m(A)}\,.$
Associate to each conjugacy class of elements of the fundamental group of
$C_n(\mathbb{C}) \diagup \mathcal{S}_n$, or, equivalently, to each
conjugacy class of $n$-braids, its extremal length, defined as
follows.\\
\begin{defn}\label{def3} Let $\hat b$ be a conjugacy
class of $n$-braids, $n \geq 2$. The extremal length $\Lambda(\hat
b)$ of $\hat b$ is defined as $ \Lambda(\hat b)= inf_{A \in
\mathcal{A}}\,
\lambda(A),$ where
$\mathcal{A}$ denotes the set of all annuli which admit a
holomorphic mapping into
$C_n(\mathbb{C}) \diagup \mathcal{S}_n$ that represents $\hat b$.
\end{defn}
It is proved in \cite{Jo1} that the extremal length $\Lambda(\hat b)$
is proportional to the entropy $h(\hat b)$
of the conjugacy class, namely $h(\hat b) = \frac{\pi}{2} \Lambda(\hat b)$.

Denote by $\Delta_n$ the Garside element in the braid group
$\mathcal{B}_n$. ( $\Delta_n^2$ is a full twist). Note that the
subgroup
$\langle \Delta_n^2 \rangle$ of $\mathcal{B}_n$ generated by
$\Delta_n^2$ is the
center $\mathcal{Z}_n$ of $\mathcal{B}_n$.  We have the
following lemma.

\begin{lemm}\label{lemm1} For each braid $b\in \mathcal{B}_n$
we
have $\Lambda(\hat b)= \Lambda(\widehat{b \Delta_n^2})$ and
$\Lambda_{tr}(b)= \Lambda_{tr}(b \Delta_n)$. \emph{\emph{}}
For each braid $b\in
\mathcal{B}_3$ we have $\Lambda_{pb}(b)= \Lambda_{pb}(b
\Delta_3)$.
\end{lemm}

\noindent {\bf Proof.} Indeed, let $R=\{x+iy: x \in (0, 1),\, y \in (0, \textsf{a})\}$, and suppose a holomorphic mapping $f:R
\to \,C_n(\mathbb{C}^n)
\diagup \mathcal{S}_n\, $ represents $b_{tr}$ (or $b_{pb}$, respectively, provided $n=3$).
Let $\tilde f=(\tilde f_1, \dots,\tilde f_n): R
\to
\,C_n(\mathbb{C}^n) $ be a lift of $f$
to a
mapping into $C_n(\mathbb{C}^n)$.
The mapping
$\zeta \to e^{\frac{\pi}{\textsf{a}}
\zeta}\tilde f(\zeta)= e^{\frac{\pi}{\textsf{a}}
\zeta}(\tilde f_1(\zeta),\ldots,\tilde f_n(\zeta)) $ is holomorphic on
$R$.
For the canonical projection
$\mathcal{P}_{\mathcal{S}_n}:\,C_n(\mathbb{C}^n)\to\,
C_n(\mathbb{C}^n)\diagup \mathcal{S}_n\,$ the mapping $\zeta
\to\,
\mathcal{P}_{\mathcal{S}_n}(e^{\frac{\pi}{\textsf{a}}
\zeta}(\tilde
f(\zeta)))\,$
represents $(b\, \Delta_n)_{tr}$ (or $(b\, \Delta_n)_{pb}$,
respectively, provided $n=3$).
The stated relation for conjugacy classes of braids is proved in the same way.
\hfill $\Box$

\medskip

The versions of the extremal length of braids are morally
related.
In this paper we will give details only for $3$-braids.

We consider now pure $3$-braids identified with elements of the
fundamental group $\pi_1(C_3(\mathbb{C}) \diagup \mathcal{S}_3, E_3)$
with base point $E_3 \in C_3(\mathbb{C}) \diagup \mathcal{S}_3$. Denote
the group of pure
$3$-braids by $\mathcal{PB}_3$.

The quotient $\mathcal{PB}_3 \diagup \langle \Delta_3^2
\rangle$ is
a free group in two generators, the class of $\sigma_1^2$ and
the
class of $\sigma_2^2$. Denote the two generators
by $a_1$ and $a_2$.

We will now state theorems on upper and lower bounds
for the versions
of the extremal length of any pure braid in terms of the word representing
its image in $\mathcal{PB}_3 \diagup
\langle \Delta_3^2 \rangle$ .
Take a non-trivial element of $\mathcal{PB}_3 \diagup \langle \Delta_3^2
\rangle$. Represent it
as reduced word $w= w_1^{n_1} \cdot
w_2^{n_2} \cdot \ldots ,\,$ where the $n_j$ are non-zero integers and
the $w_j$
are alternately equal to either $a_1$ or $a_2$. We refer to the
$w_j^{n_j}$ as the terms of the word.
We will estimate the extremal length in terms of a decomposition of the
word
into syllables. We describe now the syllable decomposition of the word.
\begin{itemize}
\item [(1)] Any term $w_j^{n_j}$ of the reduced word with
    $|n_j| \ge
2$ is a syllable.
\item [(2)] Any maximal sequence of at least two consecutive terms of the
    reduced
word which have equal power equal to either $+1$ or $-1$ is
a
syllable.
\item [(3)] Each remaining term of the reduced word enters
    with
power $+1$ or $-1$ and is characterized by the following property. The
neighbouring
term on the right (if there is one) and also the
neighbouring term
on the left (if there is one) has power different from that
of the
given one. Each term of this type is a syllable called
singleton.
\end{itemize}

Define the degree of a syllable $\mbox{deg}(\mbox{syllable})$
to be the sum of the absolute values of the powers of terms entering
the
syllable.

For example, the syllables of the word $a_2^{-1}\,a_1^2\,
a_2^{-3}\,a_1^{-1}\,a_2^{-1}\,a_1^{-1}\,a_2\,a_1^{-1}$  from
left
to right are the singleton $a_2^{-1}$, the syllable $a_1^2$ of
degree $2$, the
syllable $a_2^{-3}$ of degree $3$, the syllable
$a_1^{-1}\,a_2^{-1}\,a_1^{-1}$
of degree $3$, the singleton $a_2$ and the singleton
$a_1^{-1}$.

Label the syllables of a non-trivial word from left to right by consecutive
integral numbers $j=1,2,\ldots \;$ . Let
$d_j$ be the degree of the $j$-th syllable $\mathfrak{s}_j$. Put
\begin{equation}\label{eq1}
\mathcal{L}(w)\stackrel{def}{=} \sum_j \log (4 d_j -1).
\end{equation}
If $w$ is the identity we put $\mathcal{L}(w)=0$. Notice that for the word consisting of the single syllable $\mathfrak{s}_j$ we have $\mathcal{L}(\mathfrak{s}_j)= \log (4 d_j -1)$. Thus, $\mathcal{L}(w) = \sum \mathcal{L}(\mathfrak{s}_j)$ where the sum runs over the syllables of $w$.


The following theorem holds.

\begin{thm}\label{thm1} Let $b \in \mathcal{PB}_3$ be a pure
$3$-braid and let $w$ be the word representing
its image in $\mathcal{PB}_3 \diagup
\langle \Delta_3^2 \rangle$.
Then
\begin{equation}\label{eq1a}
\frac{1}{2\pi} \, \mathcal{L}(w) \leq \Lambda_{tr}(b)
=\frac{1}{\mathcal{M}_{tr}(b)} \leq 300 \cdot
\mathcal{L}(w),
\end{equation}
except  in the following cases: $w=a_1^n$ or
$w=a_2^n$
for an integer $n$. In these cases $\Lambda_{tr}(b)=0$, i.e
$\mathcal{M}_{tr}(b)=\infty$.

Moreover,
\begin{equation}\label{eq1b}
\frac{1}{2\pi} \, \mathcal{L}(w) \leq \Lambda_{pb}(b)=
\frac{1}{\mathcal{M}_{pb}(b)} \leq 300
\cdot \mathcal{L}(w),
\end{equation}
 except in the following case:
each term in the reduced word $w$ has the same power, which
equals either $+1$ or $-1$. In these cases $\Lambda_{pb}(b)=0$,
i.e. $\mathcal{M}_{pb}(b)=\infty$.
\end{thm}
The following propositions are immediate consequences of Theorem 1.

\begin{prop}\label{prop1} For a pure braid $b \in
\mathcal{PB}_3$
which is not one of the exceptional cases of Theorem
\ref{thm1} the
two versions of the extremal length are comparable:
$$
C_1 \,\Lambda_{tr}(b) \,\leq \,\Lambda_{pb}(b)\,\leq
\,C_2
\,\Lambda_{tr}(b)
$$
for positive constants $C_1$ and $C_2$ which do not depend on
$b$.
\end{prop}

\begin{prop}\label{prop2a}  For each element $b \in \mathcal{PB}_3$
whose image $w$ in the pure braid group modulo its center is not a singleton the estimate
$$
\frac{1}{2 \pi} \cdot \mathcal{L}(w) \leq \Lambda_{tr}(b) +
\Lambda_{pb}(b) \leq
600 \cdot \mathcal{L}(w)
$$
holds.
\end{prop}

Let $ b \in \mathcal{PB}_3 $ be a pure $3$-braid and let $w$ be
the word representing its image in  $\mathcal{PB}_3 \diagup
\langle \Delta_3^2 \rangle$. The conjugacy class $\hat w$  of elements in  $\mathcal{PB}_3 \diagup
\langle \Delta_3^2 \rangle$ corresponds to the conjugacy class $\hat b$.
Any word $\tilde w$ obtained from $w$ as follows will be
called a cyclically syllable reduced conjugate of the word $w$.
Write in reduced form the periodic word $\ldots
\,w \,w \, \ldots$ which is infinite in both directions and obtained by
repeating the
entry $w$ and no other entry. Cut off from the infinite
word a word $\tilde w$ consisting of consecutive terms of
the infinite
word, so that $\tilde w$ is conjugate to $w$ and the cuts are
between two syllables, not inside a syllable. The following
estimate of the entropy of pure braids holds.

\begin{thm}\label{thm2} Let $\hat b$ be a conjugacy class of pure $3$-braids, 
let $\hat w$ be the conjugacy class of elements of  $\mathcal{PB}_3 \diagup
\langle \Delta_3^2 \rangle$ corresponding to $\hat b$
and let $w$ be a cyclically
syllable reduced word representing the conjugacy class $\hat w$. Then
$$
\frac{1}{2 \pi} \cdot \mathcal{L}(w) \,\leq \, \Lambda (\hat b) =
\frac{2}{\pi}
h(\hat b)= (\mathcal{M}(\hat b))^{-1}
\, \leq \, 300 \cdot \mathcal{L}(w),
$$
with the following exceptions: $w=a_1^n$, $w=a_2^n$ and
$w=(a_1a_2)^n$. In these cases $\Lambda(\hat b)= h(\hat b)=0$
and
$\mathcal{M}(\hat b) =\infty$.
\end{thm}

Notice, that the exceptional braids for Theorem \ref{thm2} are
exactly the reducible pure $3$-braids. Notice also, that the
word
$w=(a_1a_2)^n a_1$ is not
cyclically syllable reduced. A
cyclically syllable reduced conjugate is $a_1^2 (a_2a_1)^{n-1}a_2$.
Hence, for a braid $b$ whose image in $\mathcal{PB}_3 \diagup
\langle \Delta_3^2 \rangle$ equals $w$ we obtain that
$\Lambda(\hat b)$ is positive and   $\Lambda_{pb}( b)= 0$. However, $\Lambda_{tr}(b)$ is
positive.

We want to point out here the following fact. Take any pure braid $b$ whose image in $\mathcal{PB}_3 \diagup
\langle \Delta_3^2 \rangle$ is a cyclically syllable reduced word
which is not one of the exceptional cases of Theorem 1. Then
the equality $\Lambda(\hat b) \geq \Lambda_{pb}(b)$ holds. However,
the extremal length of $\hat b$ may be strictly larger than the extremal length of $b_{pb}$, see the remark in the end of Section 4.
There is a respective remark for $\Lambda_{pb}$  replaced by $\Lambda_{tr}$.

We consider now arbitrary $3$-braids (not necessarily pure braids).
The following lemma holds.
\begin{lemm}\label{lemm1'}
Any braid $b\in \mathcal{B}_3$ which is not a power of $\Delta_3$ can be written in a unique way in the form
\begin{equation}\label{eq2'}
\sigma_j^k \, b_1 \, \Delta_3^{\ell}\,
\end{equation}
where $j=1$ or $j=2$, $k\neq 0$ is an integer, $\ell$ is a (not necessarily even) integer, and $b_1$ is a word in $\sigma_1^2$ and $\sigma_2^2$ in reduced form. If
$b_1$ is not the identity, then the first term of $b_1$ is a non-zero even
power of $\sigma_2$ if $j=1$, and  $b_1$ is a non-zero even  power of
$\sigma_1$ if $j=2$.
\end{lemm}

For an integer $j\neq 0$ we denote by $q(j)$ that integer neighbour
of $j$, which is closest to zero. In other words, $q(j)=j$
for each even integer $j\neq 0$. For each odd integer $j\,,$  $q(j)= j
-\mbox{sgn}(j)$, where  $\mbox{sgn}(j)$
for a non-zero integer number $j$ equals $1$ if $j$ is positive,
and $-1$ if $j$ is negative. For a braid in form \eqref{eq2'} we put
$\vartheta(b) \stackrel{def}{=}\sigma_j^{q(k)} \, b_1$.
The following theorem holds.

\begin{thm}\label{thm3}
Let $b \in \mathcal{B}_3$ be a (not necessarily pure) braid which is not a power of $\Delta_3$, and let $w$ be the word representing the image of $\vartheta(b)$ in $\mathcal{B}_3 \diagup \langle \Delta_3^2\rangle$.
Then
$$
\frac{1}{2\pi}\mathcal{L}(w) \leq \Lambda_{tr}(b) \leq 300 \cdot
\mathcal{L}(w)  \,,
$$
except in the case when $b=\sigma_j^{k}\,\Delta_3^{\ell}$ where $j=1$ or  $j=2$, $k\neq 0$ is an integer number, and $\ell$ is an arbitrary integer. In this case $\Lambda_{tr}(b)=0$.
\end{thm}

\medskip
\noindent {\bf Acknowledgment.} The author is grateful to the SFB "Raum-Zeit-Materie" at
Humboldt-University Berlin
for support in the beginning of the work on the paper, and to the
Max-Planck-Institute in Bonn where the main part of the work was done.
O.Viro stimulated the work asking about an invariant of braids rather than invariants of conjugacy classes of braids (like entropy and conformal module of conjugacy classes). The author would like to thank Alexander Wei{\ss}e for teaching how to draw
figures and for producing  the wonderful essential parts of the
difficult figures.

\medskip

\section{Coverings of $\mathbb{C} \setminus \{-1,1\}$ and
slalom curves}

It will be convenient to work with homotopy classes of curves in a more general setting.
Let $X$ be a topological space, and
let $\mathcal{E}_1$ and $\mathcal{E}_2$ be relatively closed subsets of $X$. 
Let $\textsf{h}=_{\mathcal{E}_1}{\textsf{h}}_{\,\mathcal{E}_2}$ be a homotopy class of curves
in $X$
with initial point in $\mathcal{E}_1$ and terminating point in $\mathcal{E}_2$. If
$\mathcal{E}_1=\mathcal{E}_2$ we write $\textsf{h}_{\,\mathcal{E}_1}$ instead of
$\textsf{h}=_{\mathcal{E}_1}\textsf{h}_{\,\mathcal{E}_2}$. A continuous mapping
$f$ from an open rectangle into $X$ which admits a continuous extension
to the
closure of the rectangle (denoted again by $f$) is said to represent $\textsf{h}$ if the lower
open
horizontal side is mapped to $\mathcal{E}_1$, the upper horizontal side is mapped
to $\mathcal{E}_2$
and the restriction of $f$ to the closure of each maximal vertical segment in the rectangle
represents $\textsf{h}$. In case $X=C_n(\mathbb{C})\diagup \mathcal{S}_n$ and $\mathcal{E}_1=\mathcal{E}_2$ is the totally real subspace $\mathcal{E}^n  _{tr}\stackrel{def}=C_n(\mathbb{R})\diagup \mathcal{S}_n $ we write for short
$\textsf{h}_{tr}$ instead of  $\textsf{h}_{\,\mathcal{E}^n_{tr}}$. In case $X=C_3(\mathbb{C})\diagup \mathcal{S}_3$ and  $\mathcal{E}_1=\mathcal{E}_2=\mathcal{E}^3_{pb}$ we write for short $\textsf{h}_{pb}$.

We give the following general definition of
extremal length for homotopy classes of curves.
\begin{defn}\label{def3a}For an open subset $X$ of a complex manifold
$\mathcal{X}$, two relatively closed subsets $\mathcal{E}_1$ and $\mathcal{E}_2$ of $X$ and
a
homotopy class $\textsf{h}=_{\mathcal{E}_1}\textsf{h}_{\mathcal{E}_2}$ of curves in $X$
with initial
point in $\mathcal{E}_1$ and terminating point in $\mathcal{E}_2$ the extremal length
$\Lambda(\textsf{h})$ is defined as
\begin{align}
\Lambda(\textsf{h})= & \inf \{\lambda(R): R\, \mbox{ a rectangle
which admits a holomorphic mapping to} \nonumber \\
& X \,\mbox{ that
represents}
\; \textsf{h}\}\,. \nonumber
\end{align}
\end{defn}

We come now back to the full braid group $\mathcal{B}_3$ identified  with the
relative fundamental group $\pi_1(C_3(\mathbb{C}) \diagup
\mathcal{S}_3, C_3(\mathbb{R}) \diagup \mathcal{S}_3)$. The elements of the relative fundamental group are represented by arcs with endpoints in the totally real space $C_3(\mathbb{R}) \diagup \mathcal{S}_3$. Recall that the preimage
$C_3(\mathbb{R})$ of the totally real set $C_3(\mathbb{R}) \diagup
\mathcal{S}_3$ under $\mathcal{P}_{\mathcal{S}_3}$
has several connected components.
Let $z=(z_1,z_2,z_3)\in C_3(\mathbb{C})$. Denote by $M_z$ the
M\"{o}bius
transformation that maps $z_1$ to $0$, $z_3$ to $1$ and fixes
$\infty$. Then $M_z(z_2)$ omits $0$, $1$ and $\infty$. Notice
that $M_z(z_2)$
 is equal to the cross ratio $(z_2,z_3;z_1,\infty)
=\frac{z_2-z_1}{z_3-z_1}\cdot
\frac{z_3-\infty}{z_2-\infty}=\frac{z_2-z_1}{z_3-z_1}$. For a point $z=(z_1,z_2,z_3) \in C_3(\mathbb{C})$ we put $\mathfrak{C}(z)=2 M_z(z_2)-1$. The mapping $\mathfrak{C}$ takes
$C_3(\mathbb{C})$ to $\mathbb{C}\setminus \{-1,1\}$ and
$C_3(\mathbb{R})^0 \stackrel{def}{=}
\{(x_1,x_2,x_3) \in \mathbb{R}^3: x_1<x_2<x_3\}    $ to $(-1,1)$.

Associate to a curve $\tilde{\gamma}(t)\,=\, (\tilde \gamma_1(t),\tilde \gamma_2(t), \tilde \gamma_3(t)),\,
t \in
[0,1],$ in $C_3(\mathbb{C})$ the curve
$\mathfrak{C}(\tilde \gamma)(t)\overset{def}= 2
\,\frac{\tilde \gamma_2(t)-\tilde \gamma_1(t)}{\tilde \gamma_3(t)-\tilde \gamma_1(t)}
\,-1\,,\,
t \in
[0,1],$ in $\mathbb{C}$  
which omits the points $-1$ and $1$.

Let $\gamma$ be a curve in $C_3(\mathbb{C})\diagup \mathcal{S}_3$ with endpoints on the totally real set $C_3(\mathbb{R}) \diagup \mathcal{S}_3$. Consider its lift $\tilde \gamma$
to  $C_3(\mathbb{C})$ with initial point in a chosen connected component of $C_3(\mathbb{R})$. The curve $\mathfrak{C}(\tilde \gamma)$ is a curve in $\mathbb{C} \setminus \{-1,1\}$ with initial point in a connected component of $\mathbb{R} \setminus \{-1,1\}$ .
In this way we obtain mappings from homotopy classes $b_{tr}$ of curves in $C_3(\mathbb{C}) \diagup \mathcal{S}_3$ with endpoints in $C_3(\mathbb{R}) \diagup \mathcal{S}_3$ to homotopy classes in $\mathbb{C} \setminus \{-1,1\}$ with endpoints in
$\mathbb{R} \setminus \{-1,1\}$. Since each of the mappings is defined by a lift and a composition with a holomorphic mapping the extremal length of the obtained class of curves in the complex plane coincides with the extremal length of the class of curves in the symmetrized configuration space. This will allow us to study the extremal length of $3$-braids
with totally real boundary values by studying the extremal length of
homotopy classes of curves with end points on  $\mathbb{R}\setminus
\{-1,1\}$.
For each homotopy class $b_{tr}$  we may choose the most convenient among these mappings.

For a pure braid $b_{tr}$ we usually choose the connected component  $C_3^0(\mathbb
{R})$ and denote by $\gamma^0$ the lift to $C_3(\mathbb{C})$ with initial point in $C_3^0(\mathbb
{R})$ of a curve $\gamma$ representing $b_{tr}$.
The curve $\mathfrak{C}(\gamma^0)$ has endpoints on
$(-1,1)$. 

In the case of pure braids with $pb$ boundary values we usually choose the connected component $\overset{\circ}{\mathcal{E}^3}_{pb}=\{(z_1,z_2,z_3) \in C_3(\mathbb{C}): z_2 \in \mathbb{R}, z_1 \, \mbox{and}\, z_3 \,\mbox{complex}\, \mbox{ conjugate}\}$ of the preimage $\mathcal{P}^{-1}({\mathcal{E}}^3_{pb})$. The mapping $\mathfrak{C}$ takes
$\overset{\circ}{\mathcal{E}^3}_{pb}$ onto $i\mathbb{R}$.

The following lemma states that the
mapping $\gamma \to \mathfrak{C}(\gamma^0)$ induces a
homomorphism $\mathfrak{C}_*$ from $\mathcal{PB}_3$ 
to the fundamental group of $\mathbb {C}\setminus \{-1,1\}$
whose kernel is the subgroup $\langle \Delta_3^2
\rangle$ of
$\mathcal{PB}_3$ generated by $\Delta_3^2$.

\begin{lemm}\label{lemm1a} The group $\mathcal{PB}_3\diagup \langle
\Delta_3^2 \rangle$ is isomorphic to the fundamental group of
$\mathbb{C} \setminus \{-1,1\}$ with base point $0$.
\end{lemm}
\noindent {\bf Proof.} Identify the pure braid group $\mathcal{PB}_3$ with a subgroup of the fundamental group of $C_3(\mathbb{C}) \diagup \mathcal{S}_3$ with base point in the totally real set. For curves $\gamma$ representing
elements of this subgroup we consider the lift $\gamma^0$. In this way we
may identify the pure braid group with base point in the totally real set with the fundamental group of $C_3(\mathbb{C})$ with base point in $C_3(\mathbb{R})^0$. Choose the base point $(-1,0,1)  \in
C_3(\mathbb{C})$ and identify $\mathcal{PB}_3$ with
$\pi_1(C_3(\mathbb{C}),
(-1,0,1))$.

Let $\overset{\circ}{\gamma}(t)\,=\, (\overset{\circ}{\gamma}_1(t),\overset{\circ}{\gamma}_2(t),\overset{\circ}{\gamma} _3(t)),\,
t \in
[0,1],$ be a curve in $C_3(\mathbb{C)}$
and let
$\mathfrak{C}(\overset{\circ}{\gamma})$ be the associated curve in $\mathbb{C} \setminus \{-1,1\}$.
If $\overset{\circ}{\gamma}$ is a loop
with base
point $\overset{\circ}{\gamma}(0)=(-1,0,1) \in C_3(\mathbb{C})^0$,  then $\mathfrak{C}(\overset{\circ}{\gamma})$ is a
loop with base
point $\mathfrak{C}(\overset{\circ}{\gamma})(0)=0$. The homotopy class of
$\mathfrak{C}(\overset{\circ}{\gamma})$ in $\mathbb{C}
\setminus \{-1,1\}$ with base point $0$ depends only on the
homotopy class of $\overset{\circ}{\gamma}$ with base point $(-1,0,1)$
in the configuration space $C_3(\mathbb{C})$. We obtain a
homomorphism $\mathfrak{C}_*$ from the fundamental group
$\pi_1(C_3(\mathbb{C}) , (-1,0,1))$ of
$C_3(\mathbb{C})$ with base point $(-1,0,1)$ to the fundamental group
$\pi_1 \stackrel{def}{=}\pi_1(\mathbb{C}\setminus
\{-1,1\},0)$
of the twice punctured complex plane with base point $0$.

The braids represented by the two loops $\overset{\circ}{\gamma}$ and
$\tilde{\gamma}$,
$\tilde{\gamma}(t)\stackrel{def}{=}(-1,\mathfrak{C}(\overset{\circ}{\gamma})(t),1),
\,t \in [0,1],\,$ differ by a power $\Delta_3^{2N}$ of the
full
twist $\Delta_3^2$. The number $N$ can be interpreted
as
the linking number of the first and the third strands of the
geometric braid $\overset{\circ}{\gamma}(t), \, t \in [0,1].\,$ This linking
number is obtained as follows. Discard the second strand. The resulting
braid
equals $\sigma ^{2N}$ where $N$ is the mentioned linking
number. For
the
geometric braid 
$\tilde \gamma$ the linking number of
the
first and third strand is zero. It follows that $\mathfrak{C}_*$ is
surjective
and its kernel equals $\langle \Delta_3^2 \rangle$, the
subgroup
of $\mathcal{B}_3$ generated by the full twist.  We obtain an
isomorphism
from 
$\mathcal{PB}_3 \diagup \langle \Delta_3^2\rangle
$ to the fundamental group
$\pi_1 =\pi_1(\mathbb{C}\setminus
\{-1,1\},0)$
of the twice punctured complex plane with base point $0$. \hfill $\Box$

\medskip

Notice that the mapping $\mathfrak{C}$ takes free homotopy classes of loops in
$C_3(\mathbb{C})\diagup \mathcal{S}_3$ corresponding to conjugacy classes of pure braids to free homotopy classes of loops in $\mathbb{C}\setminus \{-1,1\}$ corresponding to conjugacy classes of elements of the fundamental group of  $\mathbb{C}\setminus \{-1,1\}$.

Further, for an arc
$\overset{\circ}{\gamma}$ in $C_3(\mathbb{C})$ with both endpoints in
$C_3(\mathbb{R})^0$ the arc $\mathfrak{C}(\overset{\circ}{\gamma})$ in $\mathbb{C}\setminus
\{-1,1\}$ has endpoints in $(-1,1)$. This gives an isomorphism  from $\pi_1(C_3(\mathbb{C}),C_3(\mathbb{R})^0)$
to $\pi_1^{tr} \stackrel{def}{=}   \pi_1(\mathbb{C}\setminus\{-1,1\}, (-1,1))$
which we denote again by $\mathfrak{C}_*$ if no confusion arises.
In some cases it may be more convenient to consider the homomorphism obtained by using another connected component of $C_3(\mathbb{R})$ instead of $C_3(\mathbb{R})^0$ .

If $\overset{\circ}{\gamma}$ has endpoints on the perpendicular bisector space
$\overset{\circ}{\mathcal{E}^3}_{pb}$ then $\mathfrak{C}(\overset{\circ}{\gamma})$ has endpoints on
the imaginary axis $i\,\mathbb{R}$. We obtain a homomorphism,
denoted again by $\mathfrak{C}_*$ if no confusion arises, from
the
relative fundamental group  $\pi_1(\,C_3 (\mathbb {C})\,,\;
\overset{\circ}{\mathcal{E}^3}_{pb}\,)$
to the relative fundamental group $\pi_1^{pb}
\stackrel{def}{=}\pi_1(\,\mathbb
{C}\setminus
\{-1,1\}\,, \;i\, \mathbb{R}\,).$
In each case the
kernel of the homomorphism is the subgroup $\langle \Delta_3^2 \rangle$ of
$\mathcal{PB}_3$ generated by $\Delta_3^2$. Notice that $\pi_1$, $\pi_1^{tr}$ and $\pi_1^{pb}$ are isomorphic.

We will also consider the relative fundamental
groups $^{pb}\pi_1^{tr}$ and $^{tr}\pi_1^{pb}$ with mixed horizontal
boundary values (and the respective relative fundamental groups of the configuration space). Here the elements of $^{pb}\pi_1^{tr}$ are the homotopy
classes $_{i\mathbb{R}}\textsf{h}_{(-1,1)}$ and the elements of
$^{tr}\pi_1^{pb}$ are the homotopy classes
$_{(-1,1)}\textsf{h}_{i\mathbb{R}}$ in the space $X=\mathbb{C} \setminus \{-1,1\}$ .

The following theorem holds for mixed boundary values. Note that for mixed boundary values the estimate of the extremal length holds always while for $pb$ or $tr$ boundary values there are exceptional cases.

\medskip

\noindent {\bf Theorem $1'$.} {\it For all $w \in \pi_1$ the following inequalities hold}
\begin{align}
\frac{1}{2\pi} \, \mathcal{L}(w) \leq \Lambda(_{tr}(w)_{pb})
\leq 300 \cdot
\mathcal{L}(w), \label{eq13a}\\
\frac{1}{2\pi} \, \mathcal{L}(w) \leq \Lambda(_{pb}(w)_{tr})
\leq 300
\cdot \mathcal{L}(w).\label{eq13b}
\end{align}

\medskip

Lemma \ref{lemm1} and Lemma \ref{lemm1a}  together with the arguments above
imply the following lemma.
\begin{lemm}\label{lemm2} The invariants $\Lambda_{tr}$
$\Lambda_{pb}$  descend to invariants of
the quotient $\mathcal{PB}_3 \diagup \langle \Delta_3^2
\rangle$. For any pure braid $b \in \mathcal{PB}_3$ the following
equalities hold
\begin{align}
\Lambda_{tr}(b)& = \Lambda(\mathfrak{C}_*(b_{tr}))\,,\nonumber\\
\Lambda_{pb}(b)& = \Lambda(\mathfrak{C}_*(b_{pb}))\,,\nonumber\\
\Lambda(\hat b)& = \Lambda(\mathfrak{C}_*(\hat b))\,.\nonumber
\end{align}
\end{lemm}
For the proof of the theorems 
we will work with the fundamental group of the twice punctured complex plane rather than with the fundamental group of the symmetrized configuration space.
We will lift mappings with image in  $\mathbb{C} \setminus \{-1,1\}$ to
mappings with image in a suitable covering
of the twice punctured plane. It will be convenient to use instead
of the universal covering of $\mathbb{C} \setminus \{-1,1\}$
two different coverings of $\mathbb{C} \setminus \{-1,1\}$ by
$\mathbb{C}
\setminus i \mathbb{Z}$. Though the universal covering is well
studied,
factorizing the universal covering through two different coverings by
$\mathbb{C} \setminus i \mathbb{Z}$ has the advantage to make the
contribution
of the syllables to the extremal length transparent. This is especially
crucial for the proof of the upper bound.

To obtain the first covering $\mathbb{C} \setminus i \mathbb{Z} \to
\mathbb{C}
\setminus \{-1,1\}$ we take the
universal covering of the twice punctured Riemann sphere
$\mathbb{P}^1 \setminus \{-1,1\}$ (the logarithmic covering) and
remove all preimages of $\infty$ under the covering map.
Geometrically the logarithmic covering of $\mathbb{P}^1
\setminus
\{-1,1\}$ can be described as follows. Take copies of
$\mathbb{P}^1
\setminus[-1,1]$ labeled by the set $\mathbb{Z}$ of integer
numbers. Attach to each copy two copies of
$(-1,1)$, the
$+$-edge (the accumulation set of points of the upper
half-plane)
and the $-$-edge (the accumulation set of points of the lower
half-plane). For each $k \in \mathbb{Z}$ we glue the $+$-edge
of the
$k$-th copy to the $-$-edge of the $k+1$-st copy (using the
identity
mapping on $(-1,1)$ to identify points on different edges).
Denote
by $U_{\log}$ the set obtained from the described covering by
removing all preimages of $\infty$.

Choose curves $\alpha_j(t), t \in [0,1],\, j=1,2,$  in
$\mathbb{C}
\setminus \{-1,1\}$  which represent the
generators $a_j,j=1,2,$ of the fundamental group $\pi_1(\mathbb{C}
\setminus \{-1,1\},\,0\,\}$ and have the following properties. The
initial point and the terminating point of each of the curves are the
only points of the curve on the
interval $(-1,1)$ and are also the only points of the curve on the
imaginary axis.
The following proposition holds.

\begin{prop}\label{prop2} The set $U_{\log}$ is conformally
equivalent to $\mathbb{C}\setminus  i\,\mathbb{Z}$. The
mapping
$f_1 \circ f_2$, $f_2(z)= \frac{e^{\pi z} -1}{e^{\pi z} +1},
\, z
\in
\mathbb{C}\setminus  i\,\mathbb{Z}$, $f_1(w)=
\frac{1}{2}(w+\frac{1}{w}),\, w \in \mathbb{C}\setminus\{-1,0,1\}$,
is
a covering map from $\mathbb{C}\setminus  i \,\mathbb{Z}$ to
$\mathbb{C} \setminus \{-1,1\}$.

For each $k \in \mathbb{Z}$ the lift of $\alpha_1$ with initial point $\frac{- i}{2} +  i
k$ is a curve which joins $\frac{- i}{2} +  i k$ with
$\frac{- i}{2} +  i (k+1)$ and is contained in the closed
left half-plane. The only points on the imaginary axis are
the
endpoints.

The lift of $\alpha_2$ with initial point $\frac{- i}{2} +  i
k$ is a curve which joins $\frac{- i}{2} +  i k$ with
$\frac{- i}{2} +  i (k-1)$ and is contained in the closed
right half-plane. The only points on the imaginary axis are
the
endpoints.
\end{prop}

\begin{figure}[H]
\begin{center}
\includegraphics[width=15cm]{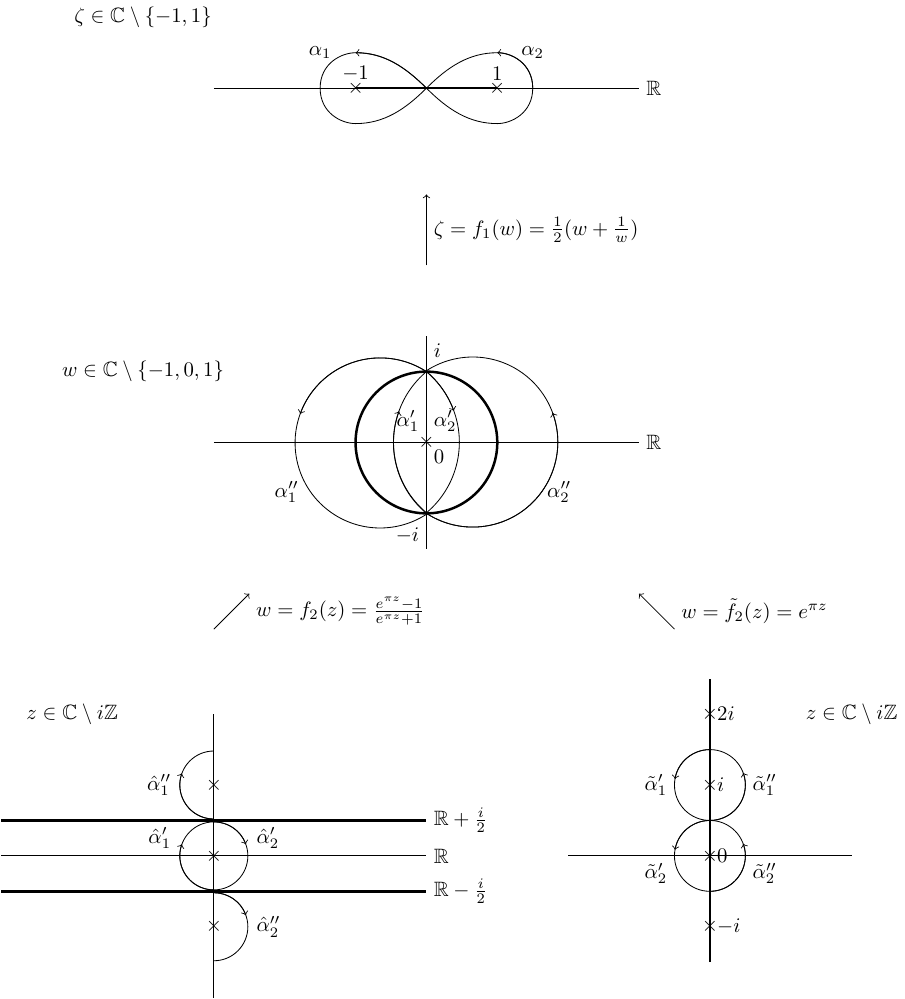}
\end{center}
\end{figure}

\centerline {Figure 1}

\medskip

Figure 1 shows the curves $\alpha_1$ and $\alpha_2$ which represent the
generators of the fundamental group $\pi_1(\mathbb{C} \setminus
\{-1,1\} ,0)$
and their lifts under the covering maps $f_1$ and $f_2 \circ f_1$ . For
$j=1,2$
the curves $\alpha_j '$ and $\alpha_j ''$ are the two lifts of
$\alpha_j$
under the double covering $f_1: \mathbb{C} \setminus \{-1,0,1\} \to
\mathbb{C}  \setminus \{-1,1\}$. The curve $\hat {\alpha}
_1'$ is the lift of $\alpha _1'$ under the mapping $f_2$ with initial
point
$\frac{-i}{2}$, the curve $\hat {\alpha} _1''$ is
the lift of $\alpha _1''$ under the mapping $f_2$ with initial point
$\frac{i}{2}$,
the curve $\hat {\alpha} _2'$ lifts $\alpha _2'$ and has
initial point $\frac{i}{2}$, and the curve $\hat {\alpha} _2''$ lifts
$\alpha
_2''$ and has
initial point $-\frac{i}{2}$.

\medskip

\noindent \textbf{Proof}.
The mapping $f_1$ is the restriction of the Zhukovsky function to
$\mathbb{C}\setminus\{-1,0,1\}$. The
Zhukovski function defines a double
branched covering of the Riemann sphere $\mathbb{P}^1$ with branch
locus $\{-1,1\}$. In particular, the Zhukovsky function provides a
conformal mapping
from the unit disc
$\mathbb{D}$ onto $\mathbb{P}^1 \setminus [-1,1]$. It maps
$-1$ to
$-1$, $1$ to $1$ and $0$ to $\infty$. The upper half-circle is
mapped onto the $-$-edge, the upper half-disc is mapped onto the
lower
half-plane, the lower half-circle is mapped onto the $+$-edge
and the
lower half-disc is mapped onto the upper half-plane.
Similarly, it provides a conformal mapping of the exterior of the
closed unit
disc onto $\mathbb{P}^1 \setminus [-1,1]$ which preserves the upper
half-plane
and also preserves the lower half-plane.

The mapping $f_2$ extends to an infinite covering of $\mathbb{P}^1
\setminus \{-1,1\}$ by $\mathbb{C}$. By an abuse of notation we denote this extension also by $f_2$. (The extension of) $f_2$ provides a
conformal mapping $f_2^0$ from the strip $\{z \in \mathbb{C}:
-\frac{1}{2}
< \mbox{Im} \, z < \frac{1}{2}\}$ onto the unit disc, which takes the real
axis onto the segment $(-1,1)$ , such that
$\lim _{x \in \mathbb{R}, x \to -\infty}= -1$, $\lim _{x \in
\mathbb{R}, x
\to +\infty}=1$. Further, $f_2$ maps
$\frac{i}{2}$ to $i$, $-\frac{i}{2}$ to $-i$, and $0$ to $0$.
The line $\{z \in \mathbb{C}: \mbox{Im} \, z
=\frac{1}{2}\}$ is mapped onto the upper half-circle, the upper
half-strip $\{z \in \mathbb{C}: 0 < \mbox{Im}\, z <
\frac{1}{2}\}$ is mapped onto the upper half-disc, the line $\{z
\in
\mathbb{C}: \mbox{Im} \, z =-\frac{1}{2}\}$ is mapped onto the lower
half-circle and the lower half-strip is mapped onto the lower
half-disc.

We obtain the following properties for the composition $f=f_1 \circ
f_2: \mathbb{C} \setminus i \mathbb{Z} \to \mathbb{C} \setminus \{-1,1\}$. The mapping $f$ takes the punctured strip $\{z \in \mathbb{C} \setminus \{0\}: -\frac{1}{2}< \mbox{Im} \, z < \frac{1}{2}\}$ onto $\mathbb{C} \setminus [-1,1]$, so that
the upper half-strip $\{z \in \mathbb{C}: 0
< \mbox{Im} \, z < \frac{1}{2}\}$ is mapped onto the lower half-plane and the
lower half-strip is mapped onto the upper half-plane. Since both mappings
$f_1$ and $f_2$ take points in the left half-plane to points in the
left half-plane, the mapping $f$ takes the left half-strip $\{z \in
\mathbb{C}:- \frac{1}{2} < \mbox{Im} \, z < \frac{1}{2},\,
\mbox{Re}\, z<0\}$ onto the subset $\{z \in
\mathbb{C}:\mbox{Re} \, z<0\}
\setminus [-1,0]$ of the left half-plane, respectively, it
takes the
right half-strip $\{z \in \mathbb{C}:- \frac{1}{2} <
\mbox{Im} \, z
< \frac{1}{2},\, \mbox{Re} \, z>0\}$ onto the subset $\{z \in
\mathbb{C}:\mbox{Re} \, z>0\} \setminus (0,1]$ of the right
half-plane.
This implies that the lift
of the curve $\alpha_1$ under $f_1 \circ f_2$ with initial point $\frac{- i}{2}$ is a
curve which joins $\frac{- i}{2}$ with $\frac{+ i}{2}$ and is
contained in the intersection of  the closed left half-plane
with
the strip  $\{z \in \mathbb{C}:- \frac{1}{2} \leq \mbox{Im} \,z
\leq
\frac{1}{2},\,\}$. The only points on the imaginary axis are
the
endpoints. Respectively, the lift of the curve $\alpha_2$ with
initial point $\frac{+ i}{2}$ is a curve which joins
$\frac{i}{2}$
with $\frac{-i}{2}$ and is contained in the right
half-strip $\{z \in \mathbb{C}:- \frac{1}{2} \leq \mbox{Im} \,z
\leq
\frac{1}{2},\, \mbox{Re} \, z \ge 0 \}$. Again, the only points on
the
imaginary axis are the endpoints.

The mapping $f=f_1 \circ f_2$ has period $i$. Indeed, $f_2$ has period $2i$,
$f_2(z+i)=\frac{1}{f_2(z)}$, and $f_1(\frac{1}{w})=f_1(w)$.

We proved the statement concerning the lift of curves under $f_1 \circ f_2$.

To see that $U_{\log}$ is conformally equivalent to
$\mathbb{C}
\setminus  i\,\mathbb{Z} $ we identify the set $\mathbb{C} \setminus
[-1,1] $ with the sheet of $U_{\log}$ labeled by $0$. The
conformal
mapping $f \mid \{z \in \mathbb{C} \setminus \{0\}:- \frac{1}{2} <
\mbox{Im} \, z
< \frac{1}{2}\}$ (whose image is $\mathbb{C} \setminus
[-1,1] $ ) extends by Schwarz's reflection principle through
the line $\{z \in \mathbb{C}: \mbox{Im} \, z =\frac{1}{2}\}$ which is
mapped onto the $-$-edge of $\mathbb{C} \setminus [-1,1] $. The
reflected mapping takes
$\{z \in \mathbb{C}\setminus\{i\}: \frac{1}{2} < \mbox{Im} \, z
< \frac{3}{2}\}$ conformally onto  $\mathbb{C} \setminus
[-1,1]) $. 
We identify the image of the punctured strip $\{z \in \mathbb{C}\setminus \{i\}:
\frac{1}{2} < \mbox{Im} \, z
< \frac{3}{2}\}$ with the sheet
of $U_{\log}$ labeled by $-1$.

Induction
on reflection through the lines  $\{z \in \mathbb{C}:
\mbox{Im} \,z
=\frac{1}{2} +j\},\,j \in \mathbb{Z},$ gives the conformal
mapping
from $C \setminus  i \mathbb{Z} $ onto $U_{\log}$ . \hfill $\Box$

\medskip

The second covering is given by the mapping $f_1\circ \tilde f_2:
\mathbb{C}
\setminus i \mathbb{Z} \to \mathbb{C} \setminus\{-1,1\}$ , where $f_1$
is as
before and $\tilde f_2$ is the exponential map, $\tilde f_2(z)=e^{\pi
z}$. Recall that each curve $\alpha_j,\, j=1,2,$ has two lifts
$\alpha_j'$
and
$\alpha_j''$ under $f_1$. For an illustration of the following
proposition see the right part of Figure 1.

\begin{prop}\label{prop5a} The mapping
$\tilde f_2$ takes $\mathbb{C} \setminus i\mathbb{Z}$ to
$\mathbb{C}\setminus
\{-1,0,1\}$, and $f_1$ takes the latter set to $\mathbb{C}\setminus
\{-1,1\}$.
The lifts  $\tilde
\alpha_j'$ and
$\tilde \alpha_j''$, $j=1,2,$ of $\alpha_j'$ and
$\alpha_j''$ under $\tilde f_2$ have the following properties. The
lifts
$\tilde
\alpha_j',\, j=1,2,\,$   are contained in the closed left half-plane
and are
directed
downwards (i.e in the direction of decreasing $y$), the lifts  $\tilde
\alpha_j'',\,j=1,2,\,$ are contained in the closed right half-plane and
directed
upwards.
The initial point of $\tilde \alpha_1'$ is  $i + \frac{1}{2} i$, the
initial
point of $\tilde \alpha_1''$ and $\tilde \alpha_2'$ is
 $\frac{1}{2} i$, the initial point of $\tilde \alpha _2''$ is
 $-\frac{1}{2}
 i$. All other lifts are obtained by translation by an integral
 multiple of
 $2i$.
\end{prop}

The straightforward proof is left to the reader.

Notice that the mapping $f_1\circ \tilde{ f_2}$ takes the intervals $( i k, i (k+1)), \, k\in \mathbb{Z},$ to the interval $(-1,1)$.


Consider the curve $\alpha_1^n,\;n\in \mathbb{Z}\setminus
\{0\}$. It
runs $n$ times along the curve $\alpha_1$ if $n>0$, and $|n|$
times
along the curve $\alpha_1$ with inverted orientation if $n<0$. It is
homotopic in $\mathbb{C} \setminus \{-1,1\}$ with base point
$0$ to a curve whose interior (i.e. the complement of its endpoints) is
contained in the open left
half-plane. We call a representative of $a_1^n$ whose interior is in the open left half-plane a standard representative of
$a_1^n$.
In the same way we define standard
representatives of $a_2^n$.
For each $k \in
\mathbb{Z}$ the curve $\alpha_1^n$ lifts under $f_1\circ f_2$ to a
curve with
initial point
 $\frac{- i}{2} +  i k$ and terminating point  $\frac{- i}{2}
 +  i k +
 i n$
which is contained in the closed left half-plane and omits the
points in $ i \mathbb{Z}$. Respectively, $\alpha_2^n,\;n\in
\mathbb{Z}\setminus \{0\},$ lifts under $f_1\circ f_2$ to a curve with
initial
point
 $\frac{+ i}{2} +  i k$ and terminating point  $\frac{+ i}{2}
 +  i k -
 i n$
which is contained in the closed right half-plane and omits the
points in $ i \mathbb{Z}$. The mentioned lifts are homotopic
through curves in $\mathbb{C} \setminus  i \mathbb{Z}$ with
endpoints on $i\,\mathbb{R} \setminus  i \mathbb{Z}$ to curves
with interior contained in the open (right, respectively,
left)
half-plane. Standard representatives of $a_1^n$ and $a_2^n$ lift to
such curves.

\begin{defn}\label{def4} A simple curve in  $\mathbb{C}
\setminus  i
\mathbb{Z}$ with endpoints on different connected components
of
$i\mathbb{R}
\setminus  i \mathbb{Z}$ is called an elementary slalom curve
if
its interior is
contained in
one of the open half-planes $\{z \in \mathbb{C}: \,\mbox{Re} \, z
>0\}$ or
$\{z \in \mathbb{C}: \,\mbox{Re} \, z <0\}$.

A curve in  $\mathbb{C} \setminus  i \mathbb{Z}$ is called an
elementary half slalom curve if one of the endpoints is
contained in a
horizontal line $\{z \in \mathbb{C}:\mbox{Im} \, z =k+\frac{1}{2}\}$ for
an
integer $k$
and the union of the curve with its (suitably oriented)
reflection in the line
$\{z \in \mathbb{C}:\mbox{Im} \, z = k+\frac{1}{2}\}$
is an elementary slalom curve.

 A slalom curve in  $\mathbb{C}
\setminus  i \mathbb{Z}$ is a curve which can be divided into
a
finite number of elementary slalom curves so that consecutive
elementary slalom curves are contained in different
half-planes.

A curve which is homotopic to a slalom curve (elementary slalom curve,
respectively) in $\mathbb{C} \setminus  i
\mathbb{Z}$ through curves with endpoints in $i\mathbb{R}
\setminus
 i \mathbb{Z}$ is called a homotopy slalom curve (elementary homotopy
 slalom
 curve, respectively).

A curve which is homotopic  to an elementary half-slalom curve in
$\mathbb{C}
\setminus  i
\mathbb{Z}$ through curves with one endpoint in $i\mathbb{R}
\setminus
 i \mathbb{Z}$ and the other endpoint on the line
 $\{z \in \mathbb{C}:\mbox{Im} \,z =k+\frac{1}{2}\}$ for an
integer $k$ is called an elementary homotopy half-slalom curve.
\end{defn}

We call an elementary slalom curve non-trivial if its endpoints are
contained in
intervals $(ik,i(k+1))$ and $(i\ell, i(\ell+1))$ with $|k-\ell|\geq 2$.
Note
that the union of an elementary half-slalom curve with its
reflection in
the horizontal line that contains one endpoint is always a non-trivial
elementary slalom curve (i.e. an elementary half-slalom curve is "half" of a
non-trivial elementary slalom curve).
We saw that the lifts under $f_1 \circ f_2$  of representatives of
terms
$_{pb}(a_j^n)_{pb} \in \pi_1^{pb}$ with $|n|\geq 2$ are non-trivial
elementary
homotopy slalom curves.

The lifts of representatives of elements of $\pi_1^{pb}$ under
$f_1\circ \tilde
f_2$ look different.
The representatives of $a_j^n$, $|n|\geq 1,$ lift to curves which
make $|n|$ half-turns around a point in $i\mathbb{Z}$ (positive
half-turns if
$n>0$, and negative half-turns if $n<0$). Hence, each such
representative lifts
under $f_1\circ \tilde f_2$ to the composition of $|n|$ trivial
elementary homotopy
slalom curves.

Take any syllable of from (2), i.e. any maximal sequence of at least two
consecutive terms of
the word which enter with equal power being either $1$ or $-1$. Recall
that $d$ denotes the sum of the absolute values of the powers of the terms of
the syllable. There is a representing curve that makes $d$ half-turns
around the interval $[-1,1]$ (positive half-turns, if the exponents of terms
in the
syllable are $1$, and negative half-turns otherwise). The lift under $f_1 \circ \tilde f_2$
of this representative is a non-trivial elementary slalom curve.

We will call homotopy classes of homotopy slalom curves for short slalom classes.

\begin{figure}[H]
\begin{center}
\includegraphics[width=10cm]{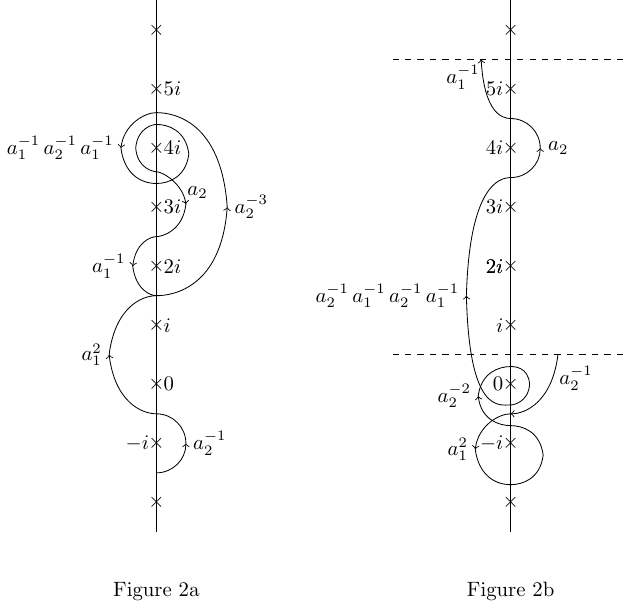}
\end{center}
\end{figure}

\centerline {Figure 2}
\bigskip

Figure 2 shows two slalom curves.  The curve in Figure 2a  is a lift
under $f_1
\circ f_2$ of a  curve in
$\mathbb{C}\setminus  \{-1,1\}$ with initial  point and terminating
point equal to $0 \in i\mathbb{R}$ representing the word
$
a_2^{-1}\,a_1^2\,
a_2^{-3}\,a_1^{-1}\,a_2^{-1}\,a_1^{-1}\,a_2\,a_1^{-1}\,
$
in the relative fundamental group $\pi_1(\mathbb{C}\setminus
\{-1,1\},i\mathbb{R})$.
The curve in Figure 2b is a lift under $f_1 \circ \tilde f_2$ of a
curve with initial and terminating point in $i\mathbb{R} \setminus \{0\}$
representing the same word. For each elementary piece of the slalom curves in Figure 2 we indicate
the element of $\pi_1$ which lifts under the considered mapping to the respective elementary slalom class.

Non-trivial elementary slalom classes and elementary half-slalom classes have positive
extremal length (in the sense of Definition 4)
which can be effectively estimated from above and from below.
In this sense homotopy classes of curves in $\mathbb{C}\setminus \{-1,1\}$ whose lifts
under $f_1 \circ f_2$ or $f_1 \circ \tilde f_2$  are
non-trivial elementary slalom classes or elementary half-slalom classes serve as building blocks. We obtained the following fact. For an element of a relative fundamental group of the
twice
punctured complex plane the pieces representing syllables of form (1)
with $pb$
boundary values lift to non-trivial elementary slalom classes
under
$f_1\circ f_2$, while the pieces representing syllables of form (2)
with $tr$
boundary values  lift to non-trivial elementary slalom classes
under
$f_1\circ \tilde f_2$.
For singletons we may select pieces of representing
curves which
lift to non-trivial elementary homotopy half-slalom curves. These facts
will be
used to obtain a lower bound for the extremal length of elements of the
relative fundamental groups.
Recall that the extremal length of a homotopy class of curves is equal
to the
extremal length of the class of their lifts under a holomorphic
covering.

The method to obtain the upper bound is roughly to patch together in a
quasiconformal way the holomorphic mappings of rectangles representing
syllables and to perturb the obtained quasiconformal mapping to a holomorphic mapping.

We conclude this section with  relating the two explicitly given
coverings of
$\mathbb{C}\setminus \{-1,1\}$ by $\mathbb{C} \setminus i\mathbb{Z}$ to
the
universal covering of the twice punctured plane.
 Denote by $U$ the
universal covering of $\mathbb{C}
\setminus  i\mathbb{Z}$. Geometrically it is obtained in the
following
way. Consider the left half-plane $\mathbb{C}_{\ell}$ and call it the Riemann surface of generation $0$. The first step is the following. For each
integer $k$  we take a copy of the right
half-plane $\mathbb{C}_r$  and glue it to $\mathbb{C}_{\ell}$ along the
interval $(ki,(k +1)i)$ (using the identity mapping for gluing). We obtain a Riemann
surface with a natural projection to
$\mathbb{C} \setminus i \mathbb{Z}$ called the Riemann surface of first generation. At the second step we consider the
Riemann surface of first generation and proceed similarly by gluing
the left half-plane along each copy
of
intervals
$(ki, (k + 1)i)$ which is an end of the Riemann
surface of first generation. By
induction
we obtain the universal covering of $\mathbb{C} \setminus i
\mathbb{Z}$.

Denote by $\widetilde{\mathbb{C}}_{\ell}$ the lift of the left
half-plane to the
first sheet of $U$ over $\mathbb{C}_{\ell}$.
Let $\mathbb{C}_{\ell}^{Cl}$ be the closure of
$\mathbb{C}_{\ell}$
in $\mathbb{C} \setminus i \mathbb{Z}$,  let
$\widetilde{\mathbb{C}}_{\ell}^{Cl}$ be	the closure of
$\widetilde {\mathbb{C}}_{\ell}$  in the
universal
covering of $\mathbb{C} \setminus i \mathbb{Z}$ and let
${(ki, (k + 1)i)}{\;\widetilde{}} \subset \widetilde{\mathbb{C}}_{\ell}^{Cl}$
be
the
lift of the intervals $(ki, (k + 1)i)$.

For each $k$ we denote by $\mathfrak{D}^{\ell}_k$  the half-disc
$\{z \in \mathbb{C}_{\ell}: |z-i(k+\frac{1}{2})|<\frac{1}{2}\}$ with
diameter
$(ik,i(k+1))$ which is contained in the left half-plane and by
$\rho_k$ the
respective open half-circle $\{z \in \mathbb{C}_{ \ell}:
|z-i(k+\frac{1}{2})|=\frac{1}{2}\}$. We call the $\rho_k$ half-circles of generation $0$.

\begin{lemm}\label{lemm3} There is a conformal mapping
$\varphi : U
\to \mathbb{C}_{\ell}$ which takes
$\widetilde{\mathbb{C}}_{\ell}$ onto $\mathbb{C}_{\ell} \setminus
\bigcup
_{k=-\infty}^{\infty} \overline{\mathfrak{D}^{\ell}_k}$ so that
${(ki, (k + 1)i)}{\;\widetilde{}}$ is mapped onto $\rho_k$ for each $k$.
\end{lemm}
\noindent {\bf Proof.} Consider the half-strip
$\mathfrak{H}_0\stackrel{def}{=}\{z \in {\mathbb{C}}_{\ell}:
0<\mbox{Im} \, z<1\}$.
By a theorem of Caratheodory (\cite{Go}, Chapter II.3, Theorem 4 and Theorem $4'$,  and also \cite{Ma} , Theorem 2.24 and Theorem 2.25) each conformal mapping which takes its lift
$\widetilde{\mathfrak{H}}_0$ onto the set $\mathfrak{H}_0\setminus
\overline{\mathfrak{D}^{\ell}_0}$, extends continuously to a
homeomorphism
between closures. Let
$\varphi_0$  be the conformal mapping whose extension to the boundary
takes the
point $\tilde 0$ over $0$ to $0$, the point $\tilde 1$ to $1$ and the
point
$\tilde{\infty}$ (considered as prime end of
$\widetilde{\mathfrak{H}}_0$)   to
$\infty$. The extension of the conformal mapping to the boundary takes
the lift ${\{z \in {\mathbb{C}}_{\ell}: \mbox{Im} \, z=1\}}{\;\widetilde{}}$
of ${\{z
\in
{\mathbb{C}}_{\ell}: \mbox{Im} \, z=1\}}$
onto ${\{z
\in
{\mathbb{C}}_{\ell}: \mbox{Im} \, z=1\}}$, it takes the lift ${\{z \in
{\mathbb{C}}_{\ell}: \mbox{Im} \, z=0\}}{\;\widetilde{}}$ of  ${\{z \in
{\mathbb{C}}_{\ell}: \mbox{Im} \, z=0\}}$ onto ${\{z \in {\mathbb{C}}_{\ell}:
\mbox{Im} \, z=0\}}$, and maps the lift ${(0, i)}{\;\widetilde{}}$ of $(0,i)$ onto $\rho_0$. By
induction we
extend $\varphi_0$ by Schwarz's reflection principle across the
half-lines
${\{z \in {\mathbb{C}}_{\ell}: \mbox{Im} \, z=k\}}{\;\widetilde{}}$ which are the lifts of the respective half-lines in ${\mathbb{C}}_{\ell}$. We obtain a conformal
mapping of $\widetilde{\mathbb{C}}_{\ell}$ onto $\mathbb{C}_{\ell}
\setminus
\bigcup _{k=-\infty}^{\infty} \overline{\mathfrak{D}^{\ell}_k}$, denoted again by $\varphi_0$, whose extension
to $\widetilde{\mathbb{C}}_{\ell}^{Cl}$  takes the segment
${(ik,
i(k+1))}{\;\widetilde{}}$ onto the half-circle $\rho_k$ for each integer number $k$ .

Schwarz's reflection principle across each segment  ${(ik,
i(k+1))}{\;\widetilde{}}$ provides extension of the conformal mapping  $\varphi_0$ to the Riemann surface of first generation. Note that for each $k_0$ the image of the half-circles $\rho_{k}$, $k\neq k_0,$
under reflection of
$\mathbb{C}_{\ell} \setminus \bigcup _{k=-\infty}^{\infty}
\overline{\mathfrak{D}^{\ell}_k}$ across $\rho_{k_0}$
are half-circles with diameter on the imaginary axis. We call them half-circles of the first generation. We obtained a conformal mapping of the Riemann surface of first generation to the unbounded connected component of the left half-plane with the half-circles of first generation removed.

Apply the reflection principle by induction to the Riemann surface of generation $n$ and all copies of intervals $(ik,i(k+1))$ in its "boundary". We obtain the Riemann surface of generation $n+1$, half-circles of generation  $n+1$ and a conformal mapping of this Riemann surface onto the unbounded connected component of the left half-plane with the half-circles of generation $n+1$ removed.

\begin{figure}[H]
\begin{center}
\includegraphics[width=4cm]{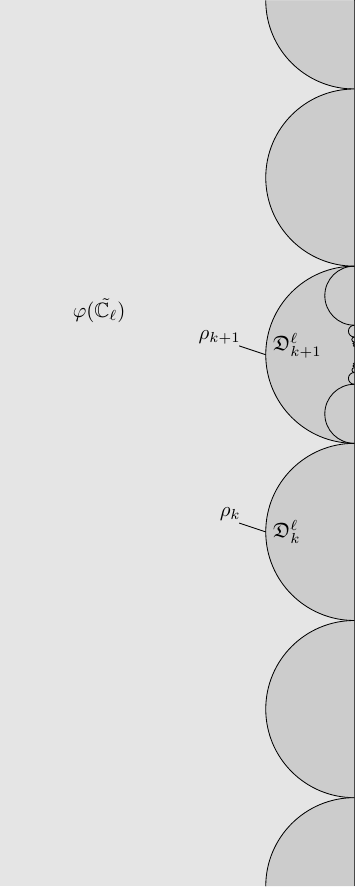}
\end{center}
\end{figure}
\centerline {Figure 3}

\medskip

The supremum of the diameters of the half-circles of generation $n$ does not exceed $2^{-n}$. Indeed, the half-circles of generation $n$ are obtained as follows. Take a half-circle $\rho$ of generation $n-1$ and reflect all other half-circles of generation $n-1$ across it. Each half-circle of generation $n-1$ different from $\rho$ has diameter on one side of the diameter of  $\rho$.
Hence, the image of each half-circle of generation
$n-1$ different from $\rho$ under reflection across $\rho$ is contained in a quarter disc, hence has diameter not exceeding half of the diameter of $\rho$. The statement is obtained by induction. By the statement we obtain in the limit a conformal mapping $\varphi$ from
$U$ onto ${\mathbb{C}}_{\ell}$.  \hfill $\Box$

\section{Building blocks and lower bound of their extremal length}

Consider an elementary slalom curve in the left half-plane with initial
point in the interval $(ik, i(k + 1))$ and with terminating point in the interval $(i j,i(j + 1))$.
Suppose that
$|k-j| \geq 2$, i.e.  the slalom curve is non-trivial.
After a
translation by a (half-integer or integer) multiple of $i$ we
obtain  a
curve 
that has endpoints in the intervals $( -i(M +1),-
iM)$ and $(iM, i(M +1))$ with $M = \frac{|k-j|-1}{2}\geq
\frac{1}{2}$. Assume that $j-k\geq2$, i.e. the curve is oriented
"upwards". The general case of elementary slalom curves in either the left or the right half-plane with initial point in the interval $(ik, i(k + 1))$ and with terminating point in the interval $(i j,i(j + 1))$ and either $k>j$ or $k<j$ is treated in the same way and leads to the same estimates related to the lower bound for the extremal length.

Let $M$ be any positive number.
Denote a curve with interior in the open left half-plane, with initial point in $( -i(M +1),-
iM)$ and terminating point in $(iM, i(M +1))$, by $\gamma_{\ell,M}$.
Notice that with the afore mentioned value $M = \frac{|k-j|-1}{2}$ the curve $\gamma_{\ell,M}$ is not an elementary slalom curve for
even $|k-j|$, but the curve $\gamma_{\ell,M} + iM$ is always an elementary
slalom curve. It will
be convenient to work  with the normalized curve $\gamma_{\ell,M}$
for estimating the extremal length.

For any integer or half-integer number $M\geq \frac{1}{2}$ we assign to the curve $\gamma_{\ell,M}$ the
class
$\gamma^*_{\ell,M}$ of curves which are homotopic to
$\gamma_{\ell,M}$
through curves in $\mathbb{C} \setminus i(\mathbb{Z}-M)$ with initial
point
in
$(-i(M + 1), -iM)$ and terminating point in $(iM, i(M+ 1))$.  The class
$\gamma^*_{\ell,M}$ is represented by the conformal mapping of an open
rectangle
onto the open left half-plane with the following property. The
conformal mapping
admits a continuous extension to the closed rectangle which takes
the open lower  side onto $(-i(M + 1), -iM)$, and the open upper side onto
$(iM, i(M + 1))$, respectively. Hence, the extremal length of the class
$\gamma^*_{\ell,M}$ is
bounded from above by the extremal length of this rectangle.

Respectively, let  $\gamma_{r,M}$ be a curve with
interior  contained
in the right half-plane ${\mathbb{C}_{r}}$ with initial point  in the
interval
$(-i(M+1),-iM)$ and terminating point in $(iM,i(M+1)$. Consider the
class
$\gamma^*_{r,M}$ of curves which are homotopic to
$\gamma_{r,M}$
through curves in $\mathbb{C} \setminus i(\mathbb{Z}-M)$ with initial
point
in
$(-i(M + 1), -iM)$ and terminating point in $(iM, i(M+ 1))$.

The conformal mappings representing $\gamma^*_{\ell,M}$ and $\gamma^*_{r,M}$ are related to elliptic integrals. With a
suitable normalization the inverse of the conformal mapping representing $\gamma^*_{\ell,M}$ is
equal to the elliptic integral
\begin{align}\label{eq3}
\mathcal{F}_M(z)& = \int_0^z
\frac{d\zeta}{\sqrt{(\zeta^2-(iM)^2)(\zeta^2-(i(M+1))^2)}}
\nonumber\\
& = \frac{i}{M+1} \int_0^{\frac{z}{iM}}
\frac{dw}{\sqrt{(1-w^2)((1 - (\frac{M}{M+1})^2 w^2)}}, \; z
\in
\mathbb{C}_{\ell}.
\end{align}

We use the branch of the square root which is positive on the
positive
real axis. The function $\mathcal{F}_M$ extends
holomorphically along any path in
$\mathbb{C}\setminus\{\pm i(M+1),\pm iM\}$, and extends
continuously to the
imaginary axis (the integral converges). The extended map takes
the union of the
closed left half-plane homeomorphically with $\infty$ to the closure of a
rectangle contained in the left half-plane. The points $-i(M+1)$,
$-iM$,
$iM$ and $i(M+1)$ are mapped by the extension of $\mathcal{F}_M$ to the vertices of the rectangle.

Note that
\begin{equation}\label{eq5a}
\frac{M+1}{i}\mathcal{F}_M(i M)= \int_0^1
\frac{dx}{\sqrt{(1-x^2)(1-k^2 x^2)}}=K(k)
\end{equation}
with $k=\frac{M}{M+1}$ is a complete elliptic integral of
first kind. Its values can be found in tables.

It is convenient to use elliptic integrals for obtaining an upper bound
for
the extremal length of arbitrary slalom curves by suitably normalizing and patching together pieces in a quasiconformal way. (See Section 5.)

As for the lower bound, the representation of $\gamma^*_{\ell,M}$ by
the inverse
of \eqref{eq3} does not give a lower bound for the extremal length of
the
elementary slalom curve by two reasons. First, the representing
mappings of
rectangles are not supposed to be conformal, they are merely
holomorphic, and
secondly, their images are not required to be contained in the
half-plane, but
merely, they are contained in $\mathbb{C}\setminus i(\mathbb{Z}-M)$ and
therefore
the mappings lift to the universal covering of $\mathbb{C} \setminus i
(\mathbb{Z}-M)$.

The following simple lemma deals with the difficulties.

\begin{lemm}\label{lemm5} Let $R$ and $R'$ be
rectangles with sides parallel to the axes. Suppose $S'$ is the
vertical strip bounded by the two vertical lines which are
prolongations of the vertical sides of the rectangle $R'$.
Let $f: R \to S'$ be a holomorphic mapping with continuous extension to the closure that takes the two horizontal sides of $R$ into different horizontal sides of $R'$. Then
$$
\lambda (R) \geq \lambda (R') \,.
$$
Equality holds if and only if the mapping is a surjective conformal
mapping from $R$ to $R'$.
\end{lemm}

The following lemma will be applied to classes of curves corresponding to conjugacy
classes of elements of the fundamental group.
\begin{lemm}\label{lemm6} Let  $A$ and $A'$ be two annuli and let $f$
be a
holomorphic mapping from $A$ into $A'$ which induces an isomorphism on
fundamental
groups.
Then
$$
\lambda(A) \geq \lambda(A')\,.
$$
Equality holds if and only if the mapping is a conformal mapping from
$A$ onto
$A'$.
\end{lemm}

\bigskip
\noindent \textbf{Proof of Lemma \ref{lemm5}}.
Normalize the rectangles and the mapping
so that ${R}= \{x+iy: x \in (0,1), y \in (0,\textsf{a}) \}\, $ and
${R}'=
\{x+iy: x \in (0,1), y \in (0,\textsf{a}')\} \, $. Denote the
continuous extension of $f$ to the closure of $R$ again by $f$.

We may assume that  $f$ maps the upper side of
${R}$ to the upper side of ${R}'$ and the lower side of ${R}$ to the
lower side of ${R}'$. Put $u= \mbox{Re}\, f$ and $v=\mbox{Im}\, f$. Then
\begin{align}\label{eq5b}
\textsf{a}'  =  \int_0^1  \textsf{a}'\,dx &= \int_0^1 (v(x, \textsf{a})
-v(x,0))dx =
\int_0^1 dx \int _0 ^{\textsf{a}} dy \frac{\partial}{\partial y}v(x,y) \\
\nonumber
& = \int_0^{\textsf{a}} dy\int_0^1 dx\, \frac{\partial}{\partial x}u(x,y) =
\int_0^{\textsf{a}} dy
\,(u(1,y)-u(0,y)) \leq \int_0^{\textsf{a}} 1\, dy =\textsf{a}.
\end{align}
 We used the Cauchy-Riemann
equations. To justify, for instance, the third equality we take the limit of the equality $v(x, \textsf{a}-\varepsilon)
-v(x,\varepsilon)=\int _{\varepsilon} ^{\textsf{a}-\varepsilon} dy\frac{\partial}{\partial y}v(x,y)$ for $\varepsilon \to +0$ and use that $v$ is continuous on the closure of $R$.

The relation \eqref{eq5a} implies the inequality $\textsf{a}' \leq \textsf{a}$.

If $\textsf{a}' = \textsf{a}$ then $u(1,y)-u(0,y)=1$ for each $y \in
(0,a)$.
Hence, the left side of ${R}$ is mapped to the left side of ${R}'$ and
the right
side  of ${R}$ is mapped to the right side of ${R}'$. Since also the
lower side
of   ${R}$ is mapped to the lower side of ${R}'$ and the upper side of
${R}$ is
mapped to the upper side of ${R}'$ the image of the positively oriented
boundary
curve of  ${R}$ has index $1$ with respect to any point of ${R}'$ and
index $0$ with respect to each point in $\mathbb{C}\setminus \bar R$.
By
the
argument principle $f(R)=R'$ and $f$ takes each value in ${R}'$ exactly
once. Hence,
$f$ is a conformal map of $R$
onto $R'$.       \hfill $\Box$

\medskip
\noindent {\bf Proof of Lemma \ref{lemm6}.}
Assume the annuli $A$ and $A'$ have center $0$,
smaller radius $1$ and larger radius $r$ and $r'$, respectively. The
set $A \setminus (0,\infty)$ is conformally equivalent to the rectangle
${R}=\{\xi+i\eta: \xi
\in (0,\log r), \eta \in (0,2 \pi)\}$. The exponential function covers
the annulus $A'$ by the strip $S'=\{\xi+i\eta: \xi \in (0,\log r'),
\eta \in \mathbb{R}\}$. We
obtain a holomorphic mapping  $g=U+iV$ from ${R}$ to $S'$ for which either
$V(\xi, 2\pi)=V(\xi,0)+2\pi$ or $V(\xi, 2\pi)=V(\xi,0)-2\pi$ for $\xi \in (0,\log r)$. Assume without loss of generality that the first option holds. Then
\begin{align}\label{eq5c}
2 \pi \log r   &=   \int _0^{\log r} (V(\xi,2\pi)-V(\xi,0))d\xi =
\int_0^{\log r} d\xi \int _0^{2\pi} d\eta \frac{\partial}{\partial
\eta}V(\xi,\eta) \\ \nonumber
& = \int_0^{\log r} d\xi \int _0^{2\pi}d\eta \frac{\partial}{\partial
\xi}U(\xi,\eta)= \int_0^{2\pi} (U(\log r,\eta)- U(0,\eta))d\eta \\
\nonumber
&\leq \int_0^{2\pi} \log r' d\eta = 2\pi \log r'.
\end{align}

Equality $r=r'$ holds iff $f$ maps the bigger circle of $A$  to the
bigger circle of $A'$ and maps the smaller circle of $A$  to the
smaller circle of $A'$. Since the map $f$ induces an isomorphism of
fundamental groups an application of the argument principle shows that
$f$ is a conformal mapping of $A$ onto $A'$.
\hfill $\Box$

Lemma \ref{lemm5} has the following two consequences.
\begin{cor}\label{cor1} For any positive integer or half-integer $M$ the extremal length of the class $\gamma^*_{\ell,M}$ is
equal to the
extremal length of a rectangle which admits a conformal mapping $f$ onto
$
\mathbb{C}_{\ell}^M \stackrel{def}{ =} \mathbb{C}_{\ell} \setminus
\big(\{|z-i(M+\frac{1}{2})|\leq \frac{1}{2}\}\cup
\{|z+i(M+\frac{1}{2})|\leq
\frac{1}{2}\}\big)$, which takes the open horizontal sides onto the two
half-circles in the
boundary of $\mathbb{C}_{\ell}^M$.
\end{cor}
\begin{cor}\label{cor2}  The extremal length of an element of each of
the relative fundamental groups
$\pi_1^{pb}$, $\pi_1^{tr}$, $^{tr}\pi_1^{pb}$, and  $^{pb}\pi_1^{tr}$
is
realized on a locally conformal mapping of a rectangle representing the
element. The extremal mapping extends locally conformally across the open
horizontal sides of the rectangle.

The extremal length of a conjugacy class of elements of the
fundamental group of the twice punctured plane is realized on a locally
conformal mapping of an annulus into the twice punctured plane.
\end{cor}

Before proving the two corollaries we recall Ahlfors'
definition of  extremal
length of a family of curves in the complex plane.

Let $\Gamma$ be a family each
member of which consists of the union of no more than countably
many (connected) locally rectifiable curves in the complex
plane. (We do not require that this union reparametrizes to a
single (connected) curve.)
In this context we will call also the elements of $\Gamma$
"curves".
Ahlfors defined the extremal length of the family $\Gamma$
as follows. For a non-negative measurable function $\varrho$ in the complex plane
he defines $A(\varrho) = \int\!\!\!\int_{\mathbb{C}} \varrho^2$. For an element
in $\gamma \in \Gamma$ and such a function $\varrho$  he puts
$L_{\gamma}(\varrho)= \int_{\gamma} \varrho |dz|$, if
$\varrho$ is measurable on $\gamma$ with respect to arc length
and $L(\varrho)=\infty$ otherwise. Put  $L(\varrho)=
\inf_{\gamma \in \Gamma}L_{\gamma}(\varrho)$.
The extremal length of the family $\Gamma$ is the following
value
$$
\lambda(\Gamma)=
\sup_{\varrho}\frac{L(\varrho)^2}{A(\varrho)},
$$
where the supremum is taken over all non-negative measurable
functions $\varrho$ for which $A(\varrho)$ is finite and does
not vanish.

It is not hard to see from this definition that the extremal length
is invariant under conformal mappings (Theorem 3 in
\cite{A1}). Further, a small  computation shows that the
previously mentioned definition of the extremal length of a
rectangle is equal to the extremal length of the family of
curves contained in the rectangle which join the opposite
horizontal sides (see Example 1 in \cite{A1}). Similarly, the extremal
length of an annulus is equal to the extremal length of the family of
curves which are contained in the annulus and represent
the conjugacy class of the positively oriented generator of the
fundamental group of the annulus.

Further, by Corollary \ref{cor1} the extremal length of the family
$\gamma^*_{\ell,M}$ in the
sense of Definition \ref{def3a} is equal to the extremal
length of
the family of curves contained in the half-plane with two half-discs
removed and
joining the two half-circles that are contained in the boundary of the half-discs. We have chosen
the general form of
Definition \ref{def3a} because it applies to classes of curves in
complex
manifolds of arbitrary dimension. Recall that the family
$\gamma^*_{\ell,M}$ was obtained from a family in a
$3$-dimensional manifold by holomorphic mappings and lifting. Notice that the
curves in the family $\gamma^*_{\ell,M}$ are not required to be
contained in $\mathbb{C}_{\ell}^M$, the half-plane
with two half-discs removed.

For later use we formulate two theorems of Ahlfors.

For two families  $\Gamma_1$ and $\Gamma_2$ as above the
following relation is introduced by Ahlfors:
$\Gamma_1 < \Gamma_2$ if each "curve"  $\gamma_2 \in
\Gamma_2$ contains a "curve" $\gamma_1 \in \Gamma_1$.

Ahlfors defines the sum $\Gamma_1+\Gamma_2$ of two such
families as follows. Each element of $\Gamma_1+\Gamma_2$ is
the union of a "curve" $\gamma_1 \in \Gamma_1$ and a "curve"
$\gamma_2 \in \Gamma_2$. The set  $\Gamma_1+\Gamma_2$ contains
all possible such unions.

The following theorems were proved by Ahlfors.

\medskip

\noindent {\bf Theorem A.}(\cite{A1}, Ch.1 Theorem 2) {\it If
$\Gamma' < \Gamma$ then $\lambda(\Gamma') < \lambda(\Gamma)$.}

\medskip

\noindent {\bf Theorem B.}(\cite{A1}, Ch.1 Theorem 4) {\it If
the families $\Gamma_j$ are contained in disjoint measurable
sets then
$\sum \lambda(\Gamma_j) \leq \lambda(\sum \Gamma_j)$.}

\medskip

\noindent {\bf Proof of Corollary \ref{cor1}.}
Lift the class $\gamma^*_{\ell,M}+iM$ (obtained from
$\gamma^*_{\ell,M}$ by translation by $iM$) to
the universal covering $U$ of $\mathbb{C} \setminus i \mathbb{Z}$ and denote the lifted class by $ ({\gamma^*_{\ell,M}+iM}) {\;\widetilde{}\;}$.
 Consider the composition  $\varphi \circ (({\gamma^*_{\ell,M}+iM}){\;\widetilde{}}\;)$ where
$\varphi$ is the conformal mapping from  $U$ onto $\mathbb{C}_{\ell}$. The obtained class is the class of curves in $\mathbb{C}_{\ell}$
which joins the half-circles $\{z \in \mathbb{C}_{\ell}:
|z+\frac{i}{2}|=\frac{1}{2}\}$ and $\{z \in \mathbb{C}_{\ell}:
|z-i(2M+\frac{1}{2})|=\frac{1}{2}\}$.
Hence, any holomorphic mapping of a rectangle into
$\mathbb{C}\setminus i\mathbb{Z}$ which
represents the class of curves $\gamma^*_{\ell,M}+iM$ results after lifting to $U$ and composing with a conformal mapping in a
holomorphic mapping of the rectangle to $\mathbb{C}_{\ell}$ which represents the class of curves in $\mathbb{C}_{\ell}$
which join the half-circles
$\rho^-\stackrel{def}{=}\{z \in \mathbb{C}_{\ell}:
|z+i(M+\frac{1}{2})|=\frac{1}{2}\}$ and $\rho^+\stackrel{def}{=}\{z \in
\mathbb{C}_{\ell}: |z-i(M+\frac{1}{2})|=\frac{1}{2}\}$.

The proof follows now from
Lemma \ref{lemm5}. To see this, identify $\mathbb{C}_{\ell}^M$
with a rectangle $R(M)$ by the conformal mapping $\mathfrak{c}$ from $R(M)$ onto
$\mathbb{C}_{\ell}^M$ whose extension to the boundary takes the horizontal
sides to the two half-circles. Denote by $\lambda_M$ the extremal
length of $R(M)$.
Notice that $\lambda_M$ is positive.
Extend the mapping $\mathfrak{c}$ by Schwarz's reflection principle across the
horizontal
sides of $R(M)$ to a rectangle $3R(M)$ of vertical side length equal to
$3$ times
the
vertical side length of $R(M)$. The rectangle added on top of $R(M)$
is mapped to the domain obtained by reflection of $\mathbb{C}_{\ell}^M$
in the half-circle $\rho^+_1=\{z \in \mathbb{C}_{\ell}: |z -
i(M+\frac{1}{2})|= \frac{1}{2}\}$. The half-circle $\rho^-_1=\{z \in
\mathbb{C}_{\ell}: |z + i(M+\frac{1}{2})|= \frac{1}{2}\}$
is mapped under this reflection to a half-circle $\rho^-_2$.
The respective fact holds for the rectangle added on the bottom of
$R(M)$.
It follows that the rectangle $3R(M)$ is mapped conformally onto the
component of $\mathbb{C}_{\ell} \setminus (\rho^+_2 \cup \rho^-_2)$
which contains $\mathbb{C}_{\ell}^M$ . This component is the left
half-plane with two half-discs removed. The half-circles $\rho^+_2$ and
$\rho^-_2$ are symmetric with respect to the real axis.

After repeating $n$ times the application of the reflection principle we obtain a
conformal mapping from a rectangle $ 3^{n} \,R(M)$ onto a domain $\Omega_n$
which is equal to the left half-plane with two half-discs removed. The
half-discs are symmetric with respect to the real axis. The domains
$\Omega_n$ are increasing.  The diameter of the removed half-discs at step $n$ tends to zero for $n \to
\infty$ since the extremal length of $ 3^{n} \,R(M)$ tends to $\infty$. We
obtain a conformal mapping of an
infinite strip $S'$ onto the left half-plane which we denote again by
$\mathfrak{c}$.

Let $f$ be any holomorphic mapping from a rectangle to the left
half-plane whose extension to the closure takes the upper side to
$\rho^+_1$ and the lower side to $\rho^-_1$. The corollary follows by
applying Lemma \ref{lemm5} to the mapping $\mathfrak{c}^{-1} \circ f$.
The form of the extremal mapping in Lemma  \ref{lemm5} shows that the extremal mapping of the corollary extends locally holomorphically across the horizontal sides of the rectangle. \hfill $\Box$

\medskip
\noindent {\bf Proof of Corollary \ref{cor2}.} For this proof it is more convenient to consider the fundamental group of the complex plane punctured at $0$ and $1$ rather than at $-1$ and $1$, and to consider the upper half plane $\mathbb{C}_+$ as universal covering of $\mathbb{C} \setminus \{0,1\}$.

Each holomorphic mapping of a
rectangle into $\mathbb{C} \setminus \{0,1\}$ that represents an element of the fundamental group lifts to the universal covering. The lift
takes the
open horizontal sides of the rectangle to certain half-circles with
diameter on
the imaginary axis (maybe, after a conformal self-map of the half-plane)
and
represents the class of curves that are contained in the half-plane and join the two
half-circles.
As in the proof of Corollary \ref{cor1} the extremal length is realized on a conformal
mapping of a rectangle
onto the half-plane with two deleted half-discs such that the horizontal sides are mapped onto the half-circles. Composing with the
covering map
we obtain a locally conformal mapping that extends locally conformally across the open horizontal sides of the rectangle.

We will now prove the statement concerning conjugacy classes of elements of the fundamental group of $\mathbb{C} \setminus
\{0,1\}$ with base point.

Note first that each element of the fundamental group corresponds to a covering
transformation.  Each covering transformation is a holomorphic
self-homeomorphism of the universal covering $\mathbb{C}_+$  that
extends to a broken linear transformation $T(z)=\frac{az+b}{cz+d}$ of
the Riemann sphere with integral coefficients $a,\,b,\,c,\,d\, ,$ such that
$ad-bc=1$. Moreover, $T$ is either parabolic (i.e. $T$ has one fixed
point and is conjugate to the mapping $z \to z+b'$ for a constant
$b'$), or $T$ is hyperbolic (i.e., T has two fixed points and is
conjugate to $z \to\kappa z$ for a positive real number $\kappa$), or
$T$ is elliptic (i.e., T has two complex fixed points symmetric with
respect to the imaginary axis and is conjugate to $z \to e^{i\theta}z$
for a real number $\theta$). See \cite{Le}, Chapter II, 9D and 9E.

Let $\hat a$ be a conjugacy class of
elements of the fundamental group of $\mathbb{C} \setminus
\{0,1\}$ and let $a$ be an element of the fundamental group
that represents $\hat a$. Denote by $T_a$
the covering transformation corresponding to $a$, and by $\langle
T_a\rangle$ the subgroup of the group of covering transformations
generated by $T_a$. Then the quotient $\mathbb{C}_+ \diagup
\langle T_a\rangle$ is an annulus. It has extremal length $0$ if $T_a$
is parabolic or elliptic and has positive extremal length if $T_a$ is
hyperbolic. If $f:A \to \mathbb{C} \setminus\{0,1\}$ is a holomorphic
mapping of an annulus $A$ to  $\mathbb{C} \setminus\{0,1\}$ that
represents $\hat a$, then $f$ lifts to a holomorphic map of $A$ into
$\mathbb{C}_+ \diagup \langle T_a\rangle$. The lift represents the
class of  a generator of the fundamental group of $\mathbb{C}_+
\diagup \langle T_a\rangle$. The corollary follows from Lemma 6.
\hfill $\Box$

\medskip

\medskip
The following proposition gives effective upper and lower bounds for the
extremal
length of the family $\gamma^*_{\ell,M}$ in dependence on the number $M\geq 2$.

\begin{prop}\label{prop4} Let $M$ be a positive integer or half-integer. The extremal
length
$\Lambda(\gamma^*_{\ell,M})$ of the class of curves $\gamma^*_{\ell,M}$
satisfies the following
inequalities:
\begin{equation}\label{eq6}
\frac{2}{\pi} \log(4M+1) \leq \Lambda(\gamma^*_{\ell,M})\leq
\frac{2}{\pi}\log(4M+3)\,.
\end{equation}
\end{prop}
\medskip

\noindent {\bf Proof.} By Corollary \ref{cor1} the extremal length
$\Lambda(\gamma^*_{\ell,M})$ equals the extremal length of a rectangle
which admits a
conformal mapping $f$ onto $ \mathbb{C}_{\ell}^M \stackrel{def}{ =}
\mathbb{C}_{\ell} \setminus \big(\{|z-i(M+\frac{1}{2})|\leq
\frac{1}{2}\}\cup
\{|z+i(M+\frac{1}{2})|\leq \frac{1}{2}\}\big)$, which takes horizontal
sides to
the two half-circles in the boundary of $\mathbb{C}_{\ell}^M$.

Put $y^{-}_M\stackrel{def}{=} \frac{1}{(4M+1)(4M+2) -1}  $ and
$y^{+}_M\stackrel{def}{=} \frac{1}{(4M+3)(4M+2) -1}  $. Notice that
$y^{+}_M
<y^{-}_M$. Denote by $D^{\ell}_M$ the open half-disc contained in the
left
half-plane with diameter $(-y^{-}_M \,i,  y^{+}_M\,i)$ and by
$D^{\ell}_1$ the
open unit half-disc in the left half-plane.
We will prove that there is a conformal self-map of the left half-plane
which takes
$\mathbb{C}_{\ell}^M$ onto $\Omega_M \stackrel{def}{=} D^{\ell}_1
\setminus
{\overline {D^{\ell}_M}}$.

Indeed, consider the conformal self-mapping $\phi_{1,M}(z) =
\frac{1}{2z+i(2M+1)}$ of the left half-plane. Its extension to the boundary
takes $iM$ to $-\frac{i}{4M+1}$, it takes $i(M+1)$ to
$-\frac{i}{4M+3}$, it takes $-iM$ to $-i$ and $-i(M+1)$ to $i$. Hence,
$\phi_{1,M}$ maps
the half-circle $\{|z+i(M+\frac{1}{2})|=\frac{1}{2}, \mbox{Re} \, z <
0\}$ to the
unit half-circle $\{|z|=1, \mbox{Re} \, z < 0\}$.   Since $\infty$ is
mapped to
$0$, $\phi_{1,M}$ maps $\mathbb{C}_{\ell} \setminus
\{|z+i(M+\frac{1}{2})|\leq \frac{1}{2}\}$ onto
$D^{\ell}_1$. Further $\phi_{1,M}$ maps the half-circle
$\{|z-i(M+\frac{1}{2})|=\frac{1}{2}, \mbox{Re}\, z < 0\}$ onto the
half-circle
with diameter $(-\frac{i}{4M+1},- \frac{i}{4M+3})$.

The mapping $\phi_{2,M}(z)=  \frac{z+\frac{i}{4M+2}}{1 -z \frac{i}{4M+2}}$
maps the
left half-plane conformally onto itself and preserves the unit
half-disc
$D^{\ell}_1$. It takes the point
$-\frac{i}{4M+1}$ to the point $- i y^{-}_M  $, and it takes the point
$-\frac{i}{4M+3}$ to the point $ i y^{+}_M $. Hence, $\phi_{2,M}$ takes
the half-disc with diameter $(-\frac{i}{4M+1},- \frac{i}{4M+3})$ onto
the half-disc with diameter $(-iy^{-}_M,i y^{+}_M)$, i.e. onto
$D^{\ell}_M$.
Hence, the composition $\phi_{2,M} \circ \phi_{1,M}$ maps  $\mathbb{C}_{\ell}^M$ conformally
onto
$\Omega_M = D^{\ell}_1 \setminus \overline{D^{\ell}_M}$.

Consider the family $\Gamma_M$ of curves in $\overline{\Omega_M}$ which join
the two
half-circles that are contained in the  boundary of $\Omega_M$.
By Corollary \ref{cor1} the extremal length of $\Gamma_M$ in the sense of
Ahlfors
equals $\Lambda(\gamma^*_{\ell,M})$. 

The following inclusions hold for the domain $\Omega_M$:
$$
D^{\ell}_1 \setminus \{|z|\leq y^{-}_M\} \subset \Omega_M \subset
D^{\ell}_1
\setminus  \{|z|\leq y^{+}_M\}\,.
$$
The first set is conformally equivalent to a rectangle with horizontal
side
length $\pi$ and vertical side length $\log ((4M+1)(4M+2)-1)$, the
second is
conformally equivalent to a rectangle with horizontal side length $\pi$
and
vertical side length $\log ((4M+2)(4M+3)-1)$.
Ahlfors's Theorem A implies the proposition. \hfill $\Box$

\medskip

Let $\tilde \gamma^*_{\ell,M}$ be the homotopy class of curves in $ \mathbb{C} \setminus i(\mathbb{Z}-M)$ that join the interval
$(-i(M+1),-iM)$ with the real axis such that  $\tilde \gamma^*_{\ell,M}$ contains a  curve $\tilde \gamma$ whose
union with its reflection in the real axis belongs to  $\gamma^*_{\ell,M}$. The following
lemma holds.

\begin{lemm}\label{lemm7}
\begin{equation}\label{eq7}
\Lambda(\tilde  \gamma^*_{\ell,M}) = \frac{1}{2} \Lambda(
\gamma^*_{\ell,M}) \,.
\end{equation}
\end{lemm}
\noindent {\bf Proof.} The extremal length $\Lambda(\tilde  \gamma^*_{\ell,M})$ equals the extremal length of a
rectangle which
admits a conformal mapping onto $\{z \in \mathbb{C}: \mbox{Re} \, z <0,
\mbox{Im} \, z
<0\} \setminus \{|z+i(M+\frac{1}{2})|\leq\frac{1}{2}\}$ such that the open
lower
side is mapped onto the half-circle in the boundary of the domain and
the open upper
side is mapped onto the negative real half-axis. This follows from Lemma \ref{lemm5} in the same way as Corollary \ref{cor1} does.
Apply Schwarz's reflection principle. By reflection in the upper side
of
the rectangle
the mapping extends to a
conformal mapping from the doubled rectangle (i.e. from the
rectangle with the same horizontal side length and double
vertical side length) onto the set $\mathbb{C}_M$, which equals the
half-plane with two removed half-discs. By Lemma \ref{lemm5} the
conformal mapping of the
doubled
rectangle realises the extremal length of $ \gamma^*_{\ell,M}$. Hence,
$$
\Lambda( \gamma^*_{\ell,M})=2 \Lambda(\tilde  \gamma^*_{\ell,M})\,.
$$
The lemma is proved. \hfill $\Box$

\medskip

The following proposition considers elements of the fundamental group
$\pi_1$
which have the form of syllables and gives upper and lower bounds for
the
extremal length of the corresponding classes of curves in the relative
fundamental groups with $pb$, $tr$ or mixed boundary values.

\begin{prop}\label{prop4b} The following statements hold for
elements of the relative fundamental groups.

\begin{itemize}
\item [1.] {\bf Syllables of form (1) with $pb$-boundary values.} The
    class
    $_{pb}(a_1^n)_{pb} \in \pi_1^{pb}$, $d=|n|\geq 2$, lifts under
    the
    covering $f_1 \circ f_2$ to a class of elementary homotopy slalom
    curves with the following property. The class
    is mapped by translation by an integer or half-integer
    multiple of
    $i$ to the class $\gamma^*_{\ell,M}$ with $M=\frac{d-1}{2}$ if
    $n>0$ and to this class with inverted orientation if $n<0$. The
    respective statement holds for $_{pb}(a_2^n)_{pb} \in
    \pi_1^{pb}$,
    $d=|n|\geq 2$ with $\gamma^*_{\ell,M}$ replaced by
    $\gamma^*_{r,M}$ with inverted orientation. (See Figure 4a for $a_2^{-3}$.) For
    the extremal length the following inequalities hold:
$$
\frac{2}{\pi} \log(2d-1)\leq \Lambda(_{pb}(a_j^n)_{pb})\leq
\frac{2}{\pi}\log(2d+1)\,.
$$

\item [2.] {\bf Syllables of form (1) or (3) with mixed boundary
    values. } For any $d=|n|\geq 1$ the class $_{pb}(a_1^n)_{tr} \in
    \,^{pb}\pi_1^{tr}$
    lifts under the covering $f_1 \circ f_2$ to a class of elementary
    homotopy
    half-slalom curves. For $n>0$ it is mapped by translation by an
    integer or
    half-integer multiple of $i$ to the class $\tilde
    \gamma^*_{\ell,M}$ with $M=d - \frac{1}{2}$. For $n<0$ we obtain the "upper half" of the class $\gamma^*_{\ell,M}$ with the orientation reversed.

    The respective statements hold for
    $_{pb}(a_2^n)_{tr}
    \,\in\, ^{pb}\pi_1^{tr}$, $d=|n|\geq 1$, with the respective halfs of
    $\tilde \gamma^*_{r,M}$ with suitable orientation.
    For the classes $_{tr}(a_j^n)_{pb}$  respective
    statements hold with the classes of translated homotopy slalom
    curves
    reflected in the real axis. (See Figure 4b for $_{pb}(a_2^{-1})_{tr}$.)

For the extremal length the following inequalities hold:
$$
\frac{1}{\pi} \log(4d-1)\leq \Lambda(_{pb}(a_j^n)_{tr})\leq
\frac{1}{\pi}\log(4d+1)\,.
$$
The same inequalities hold for the remaining cases of mixed
horizontal
boundary
values.

\item [3.] {\bf Syllables of form (2) with mixed boundary values.}
Let $_{pb}\mathfrak{s}_{tr}$ be a syllable of type (2) of degree
$d\geq 2$
with mixed boundary values. It lifts under $f_1 \circ \tilde f_2$ to
a class
of elementary homotopy half-slalom curves a suitable translation of
which is
either $\tilde \gamma^*_{\ell,M}$ or  a reflection of this class in
the real
axis or in the imaginary axis or in both, or one of these classes
with inverse
orientation. (See Figure 4c for $_{tr}(a_2a_1a_2a_1a_2)_{pb}$.) The number $M$ equals $M=d-\frac{1}{2}$.
The extremal length satisfies the inequalities
$$
\frac{1}{\pi} \log(4d-1)\leq\lambda(_{pb}(\mathfrak{s})_{tr})\leq
\frac{1}{\pi}\log(4d+1)\,.
$$

\item [4.] {\bf Syllables of form (2) with $tr$ boundary values.}
Let $_{tr}(\mathfrak{s})_{tr}$ be a syllable of type (2) with totally
real
boundary values. Then it lifts under $f_1 \circ \tilde f_2$ to a homotopy
class of
elementary slalom curves with $M=\frac{d-1}{2}$. The extremal length
of the
class satisfies the inequalities
$$
\frac{2}{\pi} \log(2d-1)\leq\lambda(_{tr}(\mathfrak{s})_{tr})\leq
\frac{2}{\pi}\log(2d+1)\,.
$$
\end{itemize}
\end{prop}

\noindent{\bf Proof.} Consider a curve $\alpha_1$ with interior in the
left
half-plane
which represents $a_1$ in the fundamental group with base point $0$.
The curve
$\alpha_1^n$ represents $a_1^n$ with $pb$ boundary values and also with
$tr$ and
mixed boundary values. Make the respective choice for $\alpha_2$.

\noindent {\bf Proof of Statement 1.} Assume first that  $n=d>0$. The curve  $\alpha_1^n$
lifts
under $f_1 \circ f_2$ to a curve in the closed left half-plane with
initial
point
$i(k-\frac{1}{2})$ and terminating point $i(k+n-\frac{1}{2})$ for an
integer number $k$. A
translation of the lifted curve  represents $\gamma^*_{\ell,M}$ with
$M=\frac{d-1}{2}$. If $n<0$ we obtain the class $\gamma^*_{\ell,M}$
with
inverted orientation. For $a_1$ replaced by $a_2$ we obtain the
respective
classes in the right half-plane. See Figure 4a for an elementary slalom
curve with $d=3$ that is the lift under $f_1 \circ f_2$ of a
representative of $_{pb}(a_2^{-3})_{pb}$ with initial point not equal
to $0$.

The estimate for the extremal length follows
from formula \eqref{eq6}.

\noindent{\bf Proof of Statement 2.} For $n=d>0$ the curve  $\alpha_1^n$ lifts under $f_1
\circ
f_2$ to an elementary homotopy half-slalom curve in $\mathbb{C}_{\ell}$ with initial point
$i(k+\frac{1}{2}) \in (ik,i(k+1))$ and terminating point
$i(k+d+\frac{1}{2}) \in i(k+d+\frac{1}{2})+\mathbb{R}$,
i.e the lift represents
$\tilde \gamma^*_{\ell,M}+i(d+k+\frac{1}{2})$ with $M=d-\frac{1}{2}$.  See Figure 4b for $_{pb}(a_2^{-1})_{tr}$ with
$d=1$. The remaining cases ( $n<0$, or $a_1$ replaced
by $a_2$, or both) can be obtained  from $\tilde
\gamma^*_{\ell,M}$
by suitable choices of inverting orientation, reflection in the real line
or reflection in the imaginary line.  The
estimate for the extremal length follows by lemma \ref{lemm7} and
\eqref{eq6}.

\begin{figure}[H]
\begin{center}
\includegraphics[width=13cm]{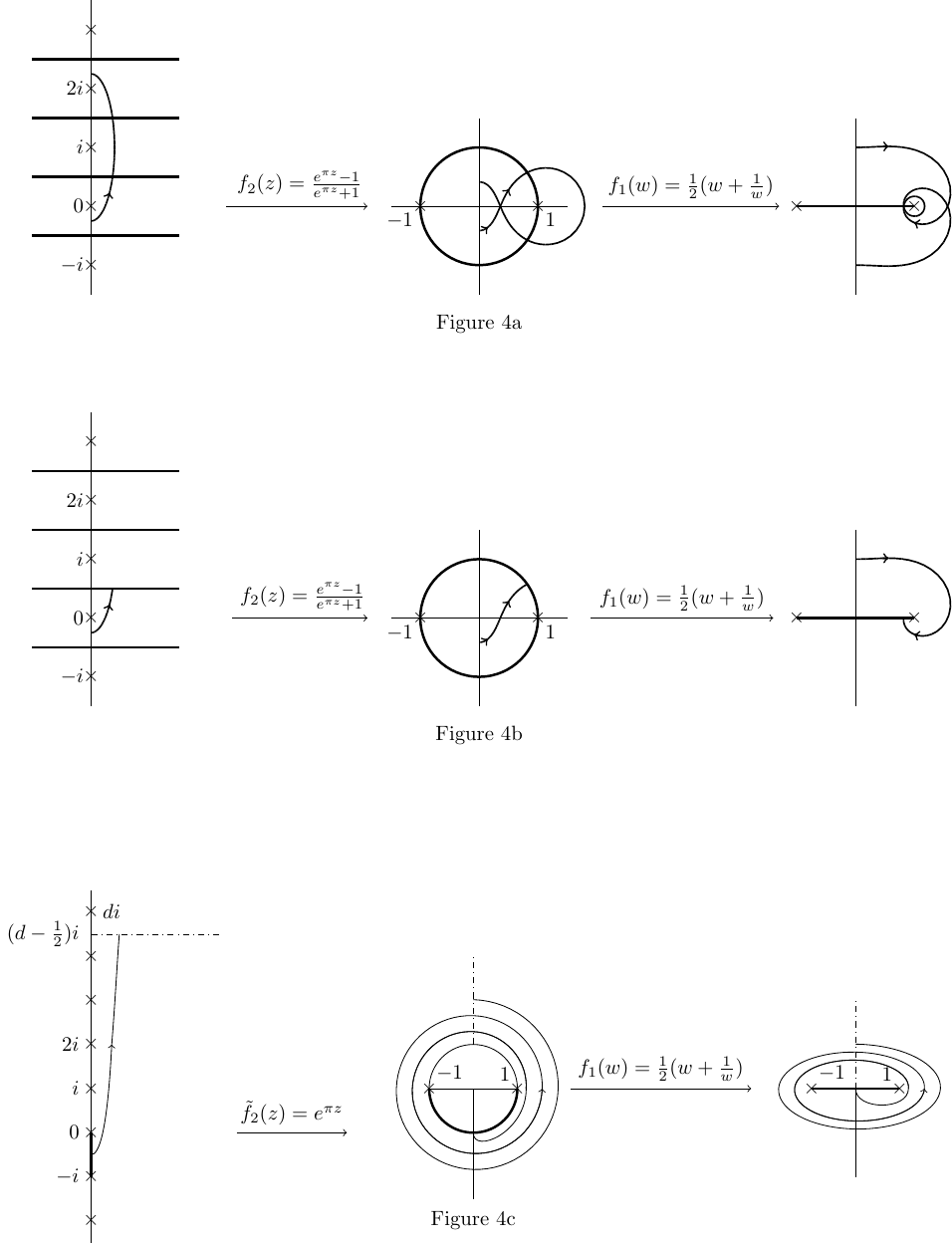}
\end{center}
\end{figure}
\centerline {Figure 4}

\bigskip

\noindent {\bf Proof of Statement 3.} Let  $\alpha_1$ and  $\alpha_2$ be as before. The
curve
$\alpha_1\,\alpha_2 \ldots$ represents
$_{tr}(\mathfrak{s})_{pb}=\,_{tr}(a_1 \,
a_2 \ldots)_{pb}$.  Its lift under $f_1 \circ \tilde f_2$  with initial
point
$i(2k -\frac{1}{2})\in ((2k-1)i,2ki)$ for an integer number $k$ can be seen as an
elementary
homotopy half-slalom curve in the closed left half-plane with
terminating
point
$i(2k-d-\frac{1}{2}) \in i(2k-d-\frac{1}{2}) + \mathbb{R}$ where $d$ is the number of letters in
$\mathfrak{s}$. Here $M= d-\frac{1}{2}$.
Notice that the lifts of
$(-1,1)$
under $f_1 \circ \tilde f_2$ are the intervals $(ik,i(k+1))$ for
integers $k$.
The lift of $\alpha_1\,\alpha_2 \ldots$ under $f_1 \circ \tilde f_2$
with
initial point $i(2k +\frac{1}{2}) \in (2k i,(2k+1)i)$ for an integer number $k$ is a curve in the closed right half-plane with terminating point $i(2k+d+\frac{1}{2}) \in i(2k+d+\frac{1}{2})+ \mathbb{R}$.

Let the first letter in the syllable be $a_2$ instead of $a_1$. The lift of the curve
$\alpha_2\,\alpha_1 \ldots$    under $f_1 \circ \tilde f_2$  with
initial point
$i(2k -\frac{1}{2}) \in ((2k-1)i,2ki)$ for an integer number $k$ is a curve in the right half-plane
with terminating point $i(2k+d-\frac{1}{2}) \in i(2k+d-\frac{1}{2})+\mathbb{R}$. See Figure 4c with number of half-turns equal to $d=5$, and also
proposition 4, as well as Figure 1.
The lift with initial
point
$i(2k +\frac{1}{2})$ is contained in the left half-plane and "moves
down".

If all generators enter with power $-1$ the orientation is reversed.
In all cases the estimate of the extremal length follows from
\eqref{eq6}.

\noindent{\bf The proof of Statement 4} is related to the proof of Statement 3 in the
same way
as the proof of Statement 1 is related to the proof of Statement 2. We
leave it
to the reader. \hfill $\Box$

\medskip

\section{The extremal length of arbitrary words in $\pi_1$. The lower bound}
Take an element of a relative fundamental group
of the twice punctured complex plane.
We will break representing curves into elementary pieces. The pieces will be
chosen so that we have a good lower bound of the extremal length of the
homotopy class of each piece. Ahlfors' theorem B will give a lower
bound for the extremal length of the element of the relative
fundamental group by the sum of the extremal lengths of the classes of
the elementary pieces.

We will use the following terminology.
Let $R$ be a rectangle. By a curvilinear rectangle contained in $R$ we mean a simply
connected domain in $R$ whose boundary looks as follows. It consist of two
vertical segments, one in each vertical side of $R$, and either two
simple arcs with interior in $R$ and endpoints on
opposite open vertical sides of $R$, or one such arc and a horizontal side
of the rectangle. The rectangle $R$ itself may also be considered as a
curvilinear rectangle contained in $R$. The arcs in $R$ with endpoints on
opposite open vertical sides are called curvilinear horizontal sides of the
curvilinear rectangle. Each curvilinear rectangle admits a conformal mapping onto a true rectangle that takes curvilinear horizontal sides to horizontal sides and vertical sides to vertical sides.

Lemma \ref{lemm14} below treats curves with $pb$ boundary
values
and is given in terms of their lifts under $f_1 \circ f_2$
to slalom curves. Lemma \ref{lemm14a}
is given in terms of relative fundamental groups of
$\mathbb{C}\setminus
\{-1,1\}$ and is more general.

We will use the symbol $\#$ for the boundary values if we are
free to choose either $pb$ or $tr$ boundary values.

\begin{lemm}\label{lemm14} Each homotopy slalom curve
is the union of elementary homotopy slalom curves.
\end{lemm}

Let $v_1$ and $v_2$ be words in $\pi_1$. For the word $_{\#}(v)_{\#} =
\,_{\#}(v_1 \,v_2)_{\#} \in \,^{\#}\pi_1\,^{\#}$ we say
that there is a sign change of exponents for the pair $(v_1,v_2)$ if
the sign of
the exponent of the last term of $v_1$ is different from the sign of
the
exponent of the first term of $v_2$.

\begin{lemm}\label{lemm14a} We have the following statements about
breaking curves into elementary pieces.
\begin{itemize}
\item [1.]\textbf{(${pb}$-boundary values between
terms, any powers.) }
Write an element  $w\in \pi_1$ as reduced word $w=w_1^{n_1} \,
w_2^{n_2} \ldots
w_k^{n_k}$, where each $w_j$
is one of the generators and the generators in consecutive
terms are different, and the $n_j$ are integral numbers
different from zero. Let $_{\#}(w)_{\#}$ be the element of the
relative
fundamental group with
$tr$, $pb$, or mixed boundary values corresponding to $w$.
Then any curve that represents $w$
in the chosen relative fundamental group is the union of
curves $\beta_j$ with the following property. For all $j,\, 1<j<k,$
$\beta_j$
represents $_{pb}(w_j^{n_j})_{pb}$. The left boundary values
of the first curve $\beta_1$ and the right boundary values of
the last curve $\beta_k$ are prescribed by the choice of the
relative fundamental group, the remaining boundary values of
$\beta_1$ and $\beta_k$ are also $pb$.

\item [2.]\textbf{($tr$-boundary values between terms,
sign change.)} Let $w_1$ be one of the standard generators $a_1$ or
$a_2$ of $\pi_1$ and let
$w_2$ be the
other standard generator. If for non-zero integers $n_1$ and $n_2$
there is a
sign change
of exponents for the pair $(w_1^{n_1},w_2^{n_2})$,  then any
representative of
$_{\#}
({w_1^{n_1}\,w_2^{n_2}})_{\#}$ is the union of two
curves
representing $_{\#} ({w_1 ^{n_1} })_{tr}$ and $ _{tr} ( {w_2 ^{n_2}})
_{\#}$,
respectively.

\item [3.]\textbf{($tr$-boundary values between terms,
one power of
absolute value bigger than 1.)} Let $w_1$ be one of the standard
generators $a_1$ or $a_2$ of $\pi_1$ and let
$w_2$ be the
other standard generator. If $|n_1|>1$ then
any representative of $ _{\#}(w_1^{n_1} \,w_2^{n_2})_{\#}$ is
the union of two
curves representing  $_{\#}(w_1^{(|n_1|-1) \begin{tiny}{\mbox{sgn}}\end{tiny}\,(n_1)})_{tr}$  and   $ _{tr}
({w_1^{\begin{tiny}{\mbox{sgn}}\end{tiny}\,(n_1)} \, w_2^{n_2}})_{\#}$, respectively.
If $|n_2|>1$ then
any representative of $ _{\#}(w_1^{n_1} \,w_2^{n_2})_{\#}$ is
the union of two
curves representing
$_{\#}(w_1^{n_1}\,w_2^{\begin{tiny}{\mbox{sgn}}\end{tiny}\,(n_2)})_{tr}$
and
$ _{tr}(w_2^{(|n_2|-1)\begin{tiny}{\mbox{sgn}}\end{tiny}(n_2)})_{\#}$, respectively.

\end{itemize}
\end{lemm}

\medskip

\noindent{\bf{Proof of Statement 1 of Lemma \ref{lemm14a}}.}
If Statement 1 is true for any smooth curve that represents $_{\#}w_{\#}$
and is transversal to $(-1,1) \cup i\mathbb{R}$,  then it is true for any representative of $_{\#}w_{\#}$. (By a smooth curve we mean a curve that is parametrized by a smooth mapping
with nowhere vanishing derivative.)

The argument is the following. Let $\beta^1:(0, \sf{a} ) \to \mathbb{C}\setminus \{-1,1\}$ be any curve representing $_{\#}w_{\#}$.
Fix a small neighbourhood $V$ of $(-1,1) \cup i\mathbb{R}$.
We may replace each component $\beta'$ of  $\beta^1 \cap V$
by a smooth curve $\beta''$  which is homotopic to $\beta'$  with the same endpoints and intersects $(-1,1) \cup i\mathbb{R}$ transversely and has no more intersection points with this set then $\beta'$.
Approximate the obtained curve by a smooth curve keeping it fixed near $(-1,1) \cup i\mathbb{R}$. We obtain a new curve that is smooth and transverse to $(-1,1) \cup i\mathbb{R}$. If the lemma holds for the new curve it holds also for the original one.

We will assume now that the curve $\beta^1$ is smooth and transverse to $(-1,1) \cup i\mathbb{R}$.
We will treat the case when $\beta^1$ has $pb$ boundary values. The other
cases are similar and we restrict ourselves to some small hints concerning the changes in the proof for those cases. Take a
smooth curve $\beta^0$ which represents  $_{pb}(w)_{pb}$, is transverse to $(-1,1) \cup i\mathbb{R}$ and lifts
under $f_1 \circ f_2$ to a smooth slalom curve (not
merely to a homotopy slalom curve).
This means that the lift of $\beta^0$ is the
union of elementary slalom curves. The collection of the elementary slalom curves projects to a collection
$\beta_j^0$ of curves representing the terms $w_j^{n_j}$ of
$w$ in $\pi_1^{pb}$. (In case some boundary values of $w$  are $tr$ we require instead of the lift of $\beta^0$ being a slalom curve that the lift of $\beta^0$ is the union of some elementary slalom curves and some  elementary half-slalom curve(s).)

Consider a smooth homotopy $h$ which joins $\beta^0$ with
$\beta^1$. In other words, there is
a smooth mapping $h$ from the closed rectangle $\bar R \stackrel{def}= \{x+iy \in \mathbb{C}: x \in [0,1],\, y \in[0,\sf{a}]\}$ to
$\mathbb{C} \setminus \{-1,1\}$ whose restriction to the left
vertical side equals $\beta^0$,  whose restriction to the
right vertical side equals $\beta^1$, and whose values on the
horizontal sides are contained in the imaginary axis
$i\mathbb{R}$ . (If some boundary values of $w$ are $tr$ then the values of the homotopy $h$ on the respective horizontal side are contained in $(-1,1)$ instead.) Since $\beta^0$ and $\beta^1$ are smooth
curves that intersect $(-1,1) \cup i\mathbb{R}$ transversely we may assume that zero is a regular value of both, $\mbox{Re}\,h$ and $\mbox{Im}\,h$. Indeed, there is a
neighbourhood $U$ of the two closed vertical sides of the rectangle
such that $0$ is a regular value of both, $\mbox{Re}\,h$ and $\mbox{Im}\,h$,  on $U$, and
we may assume by choosing the homotopy suitably that
in a neighbourhood $U'$ of the closed horizontal sides this is so.
Let $ \Psi$ be a non-negative smooth function on the closed
rectangle which equals $1$ outside $U \cup U'$ and equals zero
near the boundary of the rectangle. If $\varepsilon$ is a sufficiently small complex value such that $\mbox{Re}\,\varepsilon $ is a regular value  $\mbox{Re}\,h$ on the closed rectangle, and $\mbox{Im}\,\varepsilon $ is a regular value  $\mbox{Im}\,h$ on the closed rectangle,
then $h-\varepsilon \Psi$ is another
smooth homotopy joining $\beta^0$ with $\beta^1$ and
zero is
a regular value of both, the real part and the imaginary part of this function
on the closed rectangle.

Each connected component of the level set $\{L_0 \stackrel{def}= z \in \bar R: \mbox{Re} \,h(z)= 0\}$ on the closed rectangle is
either a closed arc that
joins two boundary points, or it is a circle. Consider the closed arcs. The intersection points of $L_0$ with
the left side of the rectangle divide the left side into connected
components. The restrictions of $h$ to these components
are the $\beta^0_j$.
Each of the components of $L_0$ that are closed arcs has at most one endpoint on the open left side of the rectangle.
Indeed, assume the contrary. Then the interval $I$ on the left side of
the
rectangle between the two endpoints of the arc contains either an
interval
with $\mbox{Re}\,h>0$  or an interval with $\mbox{Re}\,h<0$. That
means that the restriction of $h$ to the interval $I$ lifts under $f_1
\circ f_2$ to a slalom
curve and, hence, this restriction
is not homotopic to a constant through curves
in $\mathbb{C}\setminus \{-1,1\}$ with end points in the imaginary
axis. On the other hand, the existence of an arc
in the level set $\mbox{Re}\,h= 0$ joining the two endpoints of $I$
would
provide a relative homotopy in $\mathbb{C}\setminus \{-1,1\}$
(with endpoints in the imaginary axis) joining the restriction
of $h$ to the interval $I$ with a curve contained in the imaginary
axis. Such a curve would be homotopic to a constant in $\mathbb{C}\setminus \{-1,1\}$
(with endpoints in the imaginary axis) which is a contradiction.

The same reasoning shows that an arc in the level set $\{\mbox{Re}\,h=0\}$ in the rectangle $R$  with
one endpoint on the open left side of the rectangle cannot have its other endpoint on a closed horizontal side of $R$.
Indeed, as before the restriction of $h$ to
the interval on the vertical side between the endpoint of the arc and the vertex of the rectangle belonging to the closed horizontal side is not contractible within curves with endpoints on the imaginary axis, but on the other hand  it
would be homotopic through curves with endpoints on the imaginary axis to a curve contained in the imaginary axis. (In case the boundary values of $w$ corresponding to the respective horizontal side are $tr$, we get a  homotopy  through curves with endpoints in $(-1,1) \cup i\mathbb{R}$  to a curve contained in $(-1,1) \cup i\mathbb{R}$.)

We may ignore the arcs which have no endpoint on the open left side of $R$
and will also ignore the circles in the level
set of  $\mbox{Re}\,h=0$ on the closure $\bar R$.
We consider the arcs in $R$ contained in the level set $\mbox{Re}\,h=0$ that have one endpoint
on the open left side of the rectangle and the other endpoint on the open right
side. These arcs divide the rectangle into
curvilinear rectangles which are in bijective correspondence
to the intervals of division on the left side. Take the
curvilinear rectangle whose left side corresponds to
$\beta^0_j$. The restriction of $h$ to this curvilinear rectangle
provides a homotopy with boundary values in
the imaginary axes joining $\beta^0_j$ to the restriction $\beta^1_j$ of $h$ to
the right side
of the curvilinear rectangle. (In case some boundary values of $w$ are $tr$ the homotopy on the respective curvilinear rectangle may have mixed boundary values, respectively.) Since the division of the rectangle into
curvilinear rectangles induces a division of the right side of
the rectangle into intervals, the curve $\beta^1$ is the
composition of the curves $\beta_j \stackrel{def}{=}\beta^1_j$. Statement 1 is proved.


\noindent {\bf Proof of Statement 2.} Consider the case when
$n_1$ is positive and
$w_1=a_1$. The remaining cases are similar. Any representative of  $_{\#}(a_1^{n_1}\,a_2^{n_2})_{\#}$ has initial and terminating point in $(-1,1) \cup i\mathbb{R}$. For each $k \in \mathbb{Z}$ there is a lift (under $f_1\circ f_2$) of
the representative with initial point in $(ik, i(k+1)) \cup \{i(k+\frac{1}{2}) +\mathbb{R}\}$. This lift has its terminating point in $(i(k+n_1-n_2), i(k+n_1-n_2+1)) \cup (\{i(k+n_1-n_2+\frac{1}{2})\} + i\mathbb{R})$. Under the condition of Statement 2 the second exponent $n_2$ is negative. It is now easy to see that the lift of each curve
representing $_{\#}(a_1^{n_1}\,a_2^{n_2})_{\#}$
intersects the line $\{\mbox{Im}\, z = k+n_1
+\frac{1}{2}\}$. Cut the lift at an intersection point with this line and project both pieces to $\mathbb{C}\setminus \{-1,1\}$. We obtain a
decomposition of  $_{\#}(a_1^{n_1}\,a_2^{n_2})_{\#}$ into two
pieces representing $_{\#}{a_1^{n_1} }_{tr}$, and
$_{tr}{a_2^{n_2}}_{\#}$ respectively.

\noindent {\bf Proof of Statement 3.} Assume without loss of
generality that $w_1=a_1$
and $w_2=a_2$. By Statement 1 any curve representing
$_{\#}(w_1^{n_1}\,w_2^{n_2})_{\#}$ is the union
of a curve representing $_{\#}(w_1^{n_1})_{pb}$ and a curve
representing $_{pb}(w_2^{n_2})_{\#}$.

If $|n_1|>1$ and $k \in \mathbb{Z}$ the line $\{\mbox{Im}\, z =
k+\frac{1}{2}+ (|n_1|-1)\mbox{sgn}\,(n_1)\}$ does not meet the interval
$(ik,i(k+1))$.
The lift with initial point in $(ik,i(k+1))\cup \{i(k +\frac{1}{2}) + \mathbb{R}\} $ of any curve representing
$_{\#}(w_1^{n_1})_{pb}$ has its terminating point in the interval $(i(k+n_1),i(k+n_1+1))$ and, hence, intersects the line
$\{\mbox{Im}\, z =
k+\frac{1}{2}+ (|n_1|-1)\mbox{sgn}\,(n_1)\}$.
Hence, any curve representing
$_{\#}(w_1^{n_1})_{pb}$
is the union of a curve representing
$_{\#}(w_1^{(|n_1|-1)\begin{tiny}{\mbox{sgn}}\end{tiny}(n_1)})_{tr}$ and a curve
representing
$_{tr}(w_1^{\begin{tiny}{\mbox{sgn}}\end{tiny}(n_1)})_{pb}$. Hence,
any curve representing
 $_{\#}(w_1^{n_1}\,w_2^{n_2})_{\#}$ is the union
 of a curve representing $_{\#}(w_1^{(|n_1|-1) \begin{tiny}{\mbox{sgn}}\end{tiny}(n_1)})_{tr}$ and a curve
 representing $_{tr}(w_1^{ \begin{tiny}{\mbox{sgn}}\end{tiny}(n_1)}\,w_2^{n_2})_{\#}$.
 The argument for $|n_2|>1$ is similar. \hfill $\Box$

\medskip

The following lemma is a key part for the proof of the lower bound of the
extremal length of slalom classes, respectively, of elements of the relative fundamental group of $\mathbb{C}\setminus \{-1,1\}$.

\begin{lemm}\label{lemm15}
Consider a word  $w \in \pi_1$ with at least two syllables.
Take a locally conformal mapping $g:R \to \mathbb{C} \setminus
\{-1,1\}$ of a rectangle $R$ into the twice punctured plane
which  represents the word $w$  with $pb$, $tr$ or mixed
boundary values. Assume that $g$ extends locally conformally to a neighbourhood of the open horizontal sides.

The following statements hold.
\begin{itemize}
\item[1.]  $R$ can be divided into a collection of pairwise disjoint
curvilinear rectangles $R_j$ which are in bijection to the
terms of $w$ in such a way
that the restriction of g to $R_j$ represents the $j$-th
term with $pb$ boundary values if the term is not at
the (right or left) end of the word and, possibly, with mixed
boundary values if the term is at the end of the word.

\item [2.]
If $R$ contains at least one syllable of form (1) then $R$
contains a collection of mutually disjoint curvilinear rectangles
$R_{j}', \, j \in J,\,$
which are in bijective correspondence to the collection of
syllables $\mathfrak{s}_{j}$ of $w$ that are not of the form
(1) such that either the
restriction of $g$ to $R_{j}'$ represents the respective
syllable $\mathfrak{s}_{j}$ with mixed boundary values,
or the restriction represents a
syllable of form (2) with one more letter than
$\mathfrak{s}_{j}$ and mixed
boundary values.

\item[3.] If all syllables $\mathfrak{s}_j,\, j=1,\ldots,N,$ of the word $w$  are of form (2) or (3) then there is
    a
division of $R$ into curvilinear rectangles $R_j', \, j=1,\ldots,N,\,$
such that
for $j=1,\ldots,N-1,\,$ the restriction $g\mid R_j'$  represents the
$j$-th
syllable with mixed boundary values. In the same way there is a
decomposition
of $R$ into curvilinear rectangles $R_j'', \, j=1,\ldots,N,$ such
that for
$j=2,\ldots,N,$ the restriction $g\mid R_j''$  represents the $j$-th
syllable
with mixed boundary values.
\end{itemize}
\end{lemm}

\medskip

\noindent {\bf Proof}. We start with the {\bf proof of Statement 1.} 
Shrinking the
rectangle slightly in the horizontal direction we may assume
that the mapping $g$ extends locally conformally to a
neighbourhood of the closed rectangle. Hence, $0$ is a regular value of both, $\mbox{Re}\,g(z)$ and $\mbox{Im}\,g(z)$  in a neighbourhood of the closed rectangle $\bar R$.
Note that the level set
$\{z \in  R:\mbox{Re}g(z) =0\}$ in $\bar R$
does not contain circles by the maximum principle applied to the
functions $z \to exp(\pm g(z))$.

Normalize the rectangle so that we have $R=\{z \in \mathbb{C}: x \in (1,2),\, y \in (0,\sf{a})\}$. As in the proof of Lemma \ref{lemm14a} we associate to $w$ a curve $\beta^0$  which represents $w$ with the respective boundary values and lifts to a slalom curve in case the boundary values are $pb$ (in case some boundary values of $w$  are $tr$ we require instead that the lift of $\beta^0$ is the union of some elementary slalom curves and some  elementary half-slalom curve(s)).
Associate to the restriction $\beta^1$ of $g$ to the left vertical side of $R$ the homotopy $h:\overline{R_1}\stackrel{def}= \{z \in \mathbb{C}: x \in (0,1),\, y \in (0,\sf{a})\}$ of the proof of Statement 1 of Lemma \ref{lemm14a}. We may
choose the homotopy $h$ so that we obtain a smooth function $\tilde g$ on the closed rectangle
$\overline {\tilde  R}\stackrel{def}=\bar R \cup \bar{R_1}$ which equals $g$ on $\bar R$ and $h$ on $\bar{R_1}$ such that $0$ is a regular value of both, $\mbox{Re}\,{\tilde g}$ and $\mbox{Im}\,{\tilde g}$. We may consider $\tilde g$ as a homotopy of $\beta^0$ to the restriction $\beta^2$ of $\tilde g$ to the right side of the rectangle  $\tilde  R$. Consider the connected components of $\{z \in \overline{\tilde R}: \mbox{Re}\,  {\tilde g}=0\}$ that are closed arcs with one endpoint on the open left side of $\tilde R$. By the proof of Statement 1 of Lemma \ref{lemm14a} the other endpoint of each such arc lies on the open right side. We get a collection $E$ of points on the right side of the rectangle $\tilde R$ (which is also the right side of the rectangle $R$). The  mentioned arcs divide $\tilde R$ into curvilinear rectangles. We call  these arcs dividing arcs for $\tilde R$.

The restriction of $\tilde g$ to the $j$-th curvilinear rectangle in $\tilde R$ can be considered as a homotopy with endpoints in the imaginary axis
joining $\beta^0_j$ with the restriction $\beta_j^2$ of $\tilde g$ to the right side $I_j$ of the curvilinear rectangle. (If some boundary values of $w$ are $tr$ than the respective homotopy has boundary values in $(-1,1)$ on a horizontal side.)
Note that the union of the $I_j$
is the complement of $E$ in the right side of the rectangle $R$. $\beta^2_j$ represents the $j$-th term $w_j^{n_j}$ of $w$ since $\beta^0_j$  does.

Consider all connected components of the level set $\{z \in \bar R: \mbox{Re} \, {g}=0\}$ in the original rectangle $R$ that have an endpoint on the right side of $R$ contained in $E$. Each such component is a part of a dividing arc for $\tilde R$. Hence, each such component has its other endpoint on the left side of the original rectangle $R$. We call these components the dividing arcs for $R$. The dividing arcs for $R$ provide a division of the rectangle $R$ into curvilinear rectangles $R_j$. The right side of each $R_j$ is the connected components $I_j$ of the complement of $E$ in the right side of $R$ and the restriction of $g$ to $I_j$  equals $\beta^2_j$. Since $\beta^2_j$ represents the term $w_j^{n_j}$  we obtained the required collection of curvilinear rectangles $R_j$.  Statement 1 is proved.

\noindent {\bf Proof of Statement 2.} Choose a syllable $\mathfrak{s}_k$ of form (1).
Suppose there are syllables on the left of $\mathfrak{s}_k$ which are not of form (1). Consider them in the order from left to right.

If the left boundary values of $w$ are $pb$ then
we consider the most left syllable $\mathfrak{s}_{j_1}$
which is not of form (1). Suppose $j_1+1 <k$. If the next syllable $\mathfrak{s}_{j_1+1}$ to the right of
$\mathfrak{s}_{j_1}$ is also not of form (1) there is a sign change of exponents for the pair $(\mathfrak{s}_{j_1},\mathfrak{s}_{j_1+1})$ of consecutive syllables.
For the curvilinear rectangles $R_j$ of Statement 1 the restriction of
$g$ to $R_{j_1,j_1+1} \stackrel{def}{=} \mbox{Int}(\bar R_{j_1} \cup \bar
R_{j_1+1})$ represents $_{pb}(\mathfrak{s}_{j_1} \mathfrak{s}_{j_1+1})_{pb}$.
(Here $ \mbox{Int}X$ denotes the interior of a subset $X$ of a topological space.)
By Lemma \ref{lemm14a} Statement 2 and a similar argument as in the
proof of
Statement 1 of Lemma \ref{lemm15} the curvilinear rectangle $R_{j_1,j_1+1}$ can be divided into two curvilinear rectangles such that the restriction of $g$ to
them represents the syllables
$_{pb}(\mathfrak{s}_{j_1})_{tr}$ and
$_{tr}(\mathfrak{s}_{j_1+1})_{pb}$. We represented the two syllables
$\mathfrak{s}_{j_1}$ and $\mathfrak{s}_{j_1+1}$ which are both not of
form (1) as required.

Suppose the next syllable $\mathfrak{s}_{j_1+1}$ to the right of
$\mathfrak{s}_{j_1}$, is of
form (1), i.e
it equals $w_{j_1+1}^n$ for an integer $n\geq 2$ with $w_{j_1+1}$ being
a standard generator of $ \pi_1$. If there
is a sign change of exponents for the pair $(\mathfrak{s}_j,w_{j+1}^n)$
the preceding argument applies.

Suppose there is no sign change. By Lemma
\ref{lemm14a} Statement 3 and a similar argument as in the proof of
Statement 1
of Lemma \ref{lemm15} the curvilinear rectangle $R_{j_1,j_1+1}$ can be
divided into
two curvilinear rectangles such that the restriction of $g$ to them
represents the syllables
$_{pb}(\mathfrak{s}_{j_1}\, w_{j_1+1}^{\begin{tiny}{\mbox{sgn}}\end{tiny}(n)})_{tr}$
and $_{tr}(w_{j_1+1}^{n-{{\begin{tiny}{\mbox{sgn}}\end{tiny} (n)}}})_{pb}$ with mixed boundary values.

If the left boundary values of $w$ are $tr$ then $g\mid R_1$ represents $\mathfrak{s}_1$ with mixed boundary values.
Consider the most left syllable $\mathfrak{s}_{j_1}$ with $j_1>1$ which is not of form (1) and proceed as in the previous case.

The Statement 2 holds now for all syllables not of form (1) with label
$j\leq j_1+1$. The disjoint curvilinear rectangles contained in $R$
used for representing these syllables (or the syllables with one letter
more) are contained in the union of the closure of the rectangles
$R_j,\,j\leq j_1+1,$ of Statement 1 of the lemma.


We proceed by induction as follows. Suppose for some
$j'<k$ we achieved the following. We found disjoint curvilinear
rectangles contained in $\bigcup_{j\leq j'} \overline{R}_j $  which are
in bijective correspondence to all syllables
$\mathfrak{s}_j$ with $j\leq j'$ that are not of form (1) such that
the restrictions of $g$ to these rectangles represent the syllable, or
the syllable with one more letter, with mixed boundary values.

Consider the first syllable $\mathfrak{s}_{j_2}$ with $j_2<k$ on the right
of $\mathfrak{s}_{j'}$ with $\mathfrak{s}_{j_2}$
not of form (1) (if there is any). If $j_2+1<k$ we
proceed with $\mathfrak{s}_{j_2}$ in the same way as it was done
for $\mathfrak{s}_{j_1}$ in the case $j_1+1< k$ and continue the process.

Suppose there is an $\mathfrak{s}_{j_2}$ with $j_2<k$ on the right
of $\mathfrak{s}_{j'}$ with $\mathfrak{s}_{j_2}$
not of form (1) but the right neighbour of $\mathfrak{s}_{j_2}$ is the syllable
$\mathfrak{s}_k$ of form (1), say it equals $w_k^{n_k}$ with $|n_k| \geq 2$. Then we
represent $_{pb}(\mathfrak{s}_{j_2}\,w_k^{\begin{tiny}{\mbox{sgn}}\end{tiny}( n_k)})_{tr}$ and stop the process.

Make the same procedure from right to left until each
syllable not of form (1) on the right of $\mathfrak{s}_k$ is
represented in the
desired way. The curvilinear rectangles obtained by the construction
which starts from the left do not intersect the curvilinear rectangles
obtained by the construction which starts from the right because
$\mathfrak{s}_k=w_k^{n_k}$ with $|n_k|\geq 2$.
Statement 2 is proved.



\noindent {\bf Proof of Statement 3.} Under the conditions of Statement 3 there is a sign change of exponents for any pair of
consecutive syllables.

If the left boundary values of $w$ are $pb$ we
may consider curvilinear rectangles $R_{2j-1,2j}, \, j=1,2,\ldots ,$
such that
the restriction of $g$ to $R_{2j-1,2j}, \, j=1,2,\ldots , $ represents
$_{pb}(\mathfrak{s}_{2j-1} \,\mathfrak{s}_{2j})_{pb}$. Divide
$R_{2j-1,2j}$
by an arc in $\mbox{Im}\,g=0$ into two curvilinear  rectangles such
that
the
restriction of $g$ to the first curvilinear rectangle $R_{2j-1}'$
represents the
syllable $\mathfrak{s}_{2j-1}$  with right $tr$  boundary values and
the restriction of $g$ to the
second curvilinear rectangle $R_{2j}'$  represents the syllable with
left $tr$
boundary values. In this way each syllable except, maybe, the last one
is represented
with mixed boundary values by restricting $g$ to a member of the
collection of the obtained  pairwise
disjoint curvilinear rectangles.

If the left boundary values of $w$ are $tr$ then $g\mid R_1$ represents $\mathfrak{s}_1$ with mixed boundary values and we consider instead the rectangles $R_{2j,2j+1},\,j=1,\ldots$. We obtained the collection $R'_j$ for both cases of the left boundary values.

Repeating the procedure from right to
left gives the rectangles $R_j''$.
This finishes the proof of Statement 3. \hfill $\Box$

\medskip

We will  now give the proof of the lower bound in Theorem \ref{thm1}
and
in Theorem $1'$.

Figure 5 shows a mapping from a rectangle $R$ to $\mathbb{C} \setminus i\mathbb{Z}$ which represents the lift under $f_1 \circ f_2$ of an element $w \in \pi_1$ with $pb$ boundary values. The rectangle $R$ is covered by curvilinear rectangles. The right part of the figure shows the image of $R$ in $\mathbb{C} \setminus i\mathbb{Z}$ and a slalom curve that lifts a representative of $w$. For the elementary pieces of the slalom curve we indicate the element of $\pi_1$ a lift of which the piece represents. For the curvilinear rectangles contained in $R$ we indicate the element of the relative fundamental group with respective boundary values, a lift of which is represented by the restriction of the mapping. Notice that for some choices of the boundary values the curvilinear rectangles may intersect.

\begin{figure}[H]
\begin{center}
\includegraphics[width=13cm]{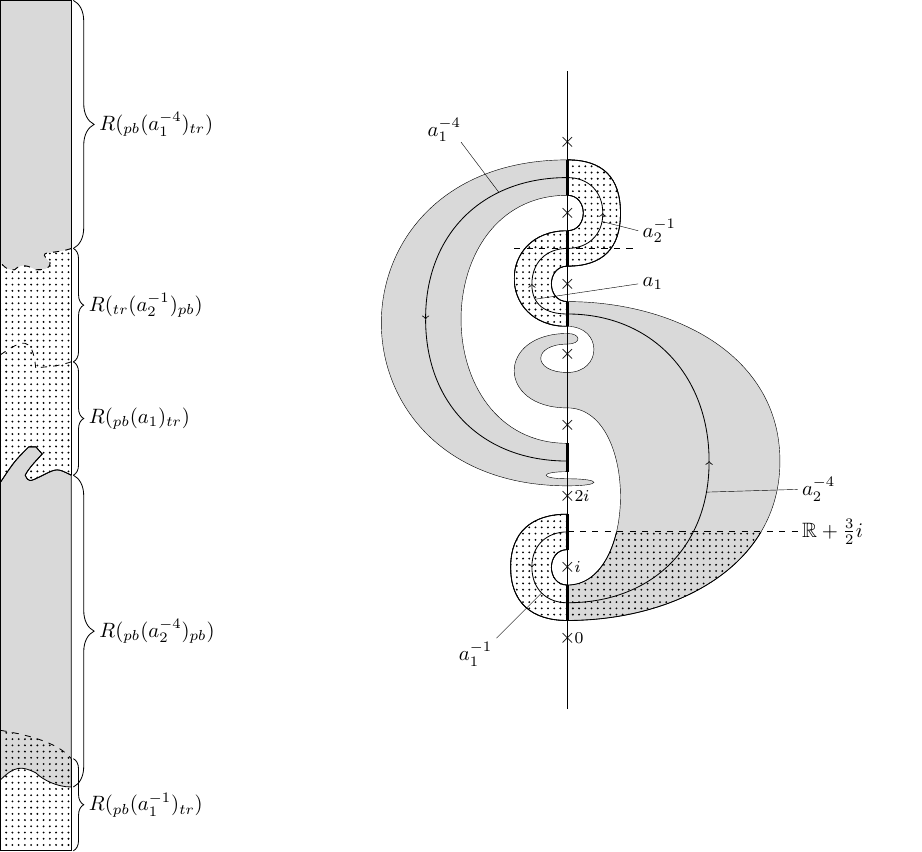}
\end{center}
\end{figure}
\centerline {Figure 5}

\medskip

\medskip

\noindent {\bf Proof of the lower bound in Theorem \ref{thm1} and in Theorem $1'$.} Identify $w$ with $\mathfrak{C}_*(b) \in \pi_1$. Let $g$
be a locally conformal mapping from a rectangle $R$ with sides parallel
to the axes into $\mathbb{C} \setminus \{-1,1\}$ which  represents the
word $w$ in the relative fundamental group of the twice punctured
complex plane with $tr$, $pb$, or mixed horizontal boundary values. We
assume first that the word has at least two syllables.
Let $\Gamma$ be the family of rectifiable curves in $R$ which join the
two
horizontal sides.

Suppose first that the word has at least one syllable of form (1). Let
$I$ be the set of the natural numbers $j$ for which the $j$-th
term of $w$ is a syllable of form (1) and let $R_j, \, j \in I,$ be the curvilinear
rectangles of Statement (1) of Lemma \ref{lemm15}. For $j\in I $ we
denote by
$\Gamma_j$ the
family of rectifiable arcs in $R_j$ which join the two horizontal
curvilinear
sides of $R_j$.
In Ahlfors' notation we have $\sum _{j \in I} \Gamma_j < \Gamma$.
By Ahlfors' Therorems A and B we obtain
\begin{equation}
 \sum _{j \in I} \lambda(\Gamma_j) \leq \lambda(\Gamma)\,.\nonumber
\end{equation}

The extremal length $\lambda(\Gamma_j)$  is estimated from
below as follows. Let $\mathcal{R}_j$ be a rectangle with
sides parallel to the axis which is conformally equivalent to
$R_j$ (i.e. there is  a conformal mapping $\psi_j: R_j \to
\mathcal{R}_j$ which takes vertical sides to vertical sides and
curvilinear horizontal sides to horizontal sides). Let
$\Gamma^*_j$ be the push forward of $\Gamma_j$ under the conformal
mapping. The extremal length is a conformal invariant
(see \cite{A1}), hence, $\lambda(\Gamma^*_j)=\lambda(\Gamma_j)$, and
$\lambda(\Gamma^*_j)$ is the extremal length of the rectangle
$\mathcal{R}_j$. The rectangle $\mathcal{R}_j$ admits a
holomorphic mapping $(g \mid R_j) \circ (\psi_j)^{-1}$ to
$\mathbb{C} \setminus \{-1,1\}$ which represents the syllable
corresponding to $R_j$ with $pb$   boundary values if the syllable is not at the end of the word, and, possibly, with mixed boundary values if the syllable is at the end of the word.

Hence, by Proposition  \ref{prop4b}  the extremal length
of  $\mathcal{R}_j$ is not smaller
than $\frac{1}{\pi} \log(4d_j-1)$ or
$\frac{2}{\pi} \log(2d_j-1) > \frac{1}{\pi} \log(4d_j-1)$
for $d_j\geq2$.
Here $d_j$ is the degree of the respective syllable.

Hence,
\begin{equation}\label{eq8}
\frac{1}{2}\lambda(\Gamma)\geq \frac{1}{2}\sum _{j \in I}
\lambda(\Gamma_j) \geq
\sum_{j \in I} \frac{1}{2\pi} \log(4d_j-1)\,.
\end{equation}
Suppose again that $w$ contains at least one syllable of form (1).
Denote by $J$ the set of all natural numbers $j$ for which
$\mathfrak{s}_j$ is
not of form (1). For $j \in J$ we denote
by $ R_j'$ the curvilinear rectangle of Statement 2 of Lemma \ref{lemm15} corresponding to $\mathfrak{s}_j$
(or to a syllable with one more letter than $\mathfrak{s}_j$).
Let $ \Gamma_j'$ be the family of rectifiable curves in
$R_j'$ which join the two curvilinear horizontal sides.

By Ahlfors' Theorems
$$
 \sum _{j \in J} \lambda( \Gamma_j') \leq \lambda(\Gamma)\,.
 $$
By the same arguments as before and
by Proposition  \ref{prop4b}  the extremal length
of  $R_j'$ is not smaller than
$\frac{1}{\pi} \log(4d_j-1)$
since $g\mid R_j'$ has mixed boundary values. (If $R_j'$ represents a syllable with one letter more than $\mathfrak{s}_j$ then the lower bound is $\frac{1}{\pi} \log(4(d_j+1)-1)$.)
As before $d_j\geq 1$ is the degree of
the syllable $\mathfrak{s}_j$.
Hence,
\begin{equation}\label{eq9}
\frac{1}{2}\lambda(\Gamma)\geq \frac{1}{2}\sum _{j \in J}
\lambda(\Gamma_j') \geq
\sum_{j \in J} \frac{1}{2\pi} \log(4d_j-1)\,.
\end{equation}
For the case when  $w$ contains syllables of form (1) and syllables not of form (1) we add the two inequalities
\eqref{eq8} and \eqref{eq9}. We obtain
\begin{equation}\label{eq10}
\lambda(\Gamma) \geq  \sum_{j \in J} \frac{1}{2\pi} \log(4d_j-1)
+ \sum_{j \in I} \frac{1}{2\pi} \log(4d_j-1) = \sum _{\mathfrak{s}_j}
\frac{1}{2\pi} \log(4d_j-1) \,.
\end{equation}
The last sum is extended over all syllables  $\mathfrak{s}_j$ of the word $w$, and $d_j$ is
the degree
of $\mathfrak{s}_j$.

Suppose that the word does not contain syllables of form (1) and the
syllables
are labeled from left to right by $j=1,\ldots,N$. Let $R_j'$ be the
rectangles from Statement 3 of Lemma \ref{lemm15} and let
$\Gamma'_j$ be the
family of curves in $R_j'$ joining the pair of horizontal sides of
$R_j'$.

Then by the same arguments as before

\begin{equation}\label{eq11}
\lambda(\Gamma) \geq \sum_{j=1}^{N-1} \lambda(\Gamma_j') \geq \sum
_{j=1}^{N-1}
 \frac{1}{\pi} \log(4d_j-1) \,.
\end{equation}

On the other hand, with the curvilinear rectangles $R_j''$ from
Statement 3 of Lemma \ref{lemm15} and the family  $\Gamma_j''$ of curves in $R_j''$
joining the pair of horizontal sides of $R_j''$ we obtain the
inequality
\begin{equation}\label{eq12}
\lambda(\Gamma) \geq \sum_{j=2}^{N} \lambda(\Gamma_j'') \geq \sum
_{j=2}^{N}
 \frac{1}{\pi} \log(4d_j-1) \,.
\end{equation}
It follows that for any word $w \in \pi_1$ with at least two syllables for the respective family $\Gamma$ the
inequality

\begin{equation}\label{eq13}
\lambda(\Gamma) \geq   \sum _{\mathfrak{s}_j}
 \frac{1}{2\pi} \log(4d_j-1)\,
\end{equation}
holds.

If the word consists of a single syllable, Proposition \ref{prop4b} implies \eqref{eq13} in the case of mixed boundary values as well as in the non-exceptional cases with both boundary values being $tr$ or $pb$.

Consider the exceptional cases. If $w=a_1^n$ with $n>0$ the mapping $\zeta \to -1+ e^{\zeta}, \; \zeta \in R,$ with $R=\{ \xi + i \eta,\, \xi \in (-\infty,0),\, \eta \in (0, 2\pi n)\}$, represents $w_{tr}$. Hence,  $\Lambda(w_{tr})=0$.

If $w=a_1 a_2\ldots$ and has degree $d \geq 2$ then the mapping $\zeta \to e^{\zeta}, \; \zeta \in R,$ with $R=\{ \xi + i \eta,\, \xi \in (1,\infty),\, \eta \in (\frac{\pi}{2},\frac{\pi}{2} + \pi d)\}$, represents $w_{pb}$. Hence, $\Lambda(w_{pb})=0$.

The other exceptional cases are similar.
The lower bound in Theorem \ref{thm1} and in Theorem $1'$ is proved. \hfill $\Box$

\medskip

\noindent {\bf Proof of the lower bound in Theorem \ref{thm2}.} Identify the class $\hat w$ with a conjugacy class of elements of the fundamental group $\pi_1$ of $\mathbb{C}\setminus \{-1,1\}$ with base point $0$. Let $\hat w$ be not among the exceptional cases of Theorem \ref{thm2}. Take any representative $w_1 \in \pi_1$ of $\hat w$ and write the infinite word $ \ldots, w_1 w_1 w_1, \ldots $ in reduced form.
If both letter $a_1$ and $a_2$ appear in the word $w_1$ and the terms of the infinite syllable reduced word do not have all equal power equal to either $1$ or $-1$ (i.e. $\hat w$ is not among the exceptional cases) then there is a decomposition of the infinite word into syllables of finite length.

Consider a locally conformal mapping $\hat g$ from an annulus $A$ into the twice punctured complex plane $\mathbb{C}\setminus \{-1,1\}$ that represents $\hat w$. The lift of $\hat g$ to the infinite strip $S$ that covers $A$ represents the infinite word.
Apply the argument of the proof of  Statement 1 of Lemma \ref{lemm15} to the infinite word and the infinite strip and project to the annulus $A$. We obtain a syllable reduced word $w\in \hat w$ and
a connected component $L_0^0$ of the level set $L_0= \{z \in A: \mbox{Re} \, g(z)=0\}$ that joins the two boundary circles of $A$ such that $A\setminus L_0^0$ is a curvilinear rectangle and the restriction of $\hat g$ to it represents the syllable reduced word $w$ with $pb$ boundary values. Lemma \ref{lemm15} applies to the restriction  $\hat g \mid (A \setminus L_0^0)$.
As in the proof of Theorem \ref{thm1} we obtain building blocks for estimating the lower bound of the extremal length. The estimates of the extremal length of the building blocks are put together in the same way as in the proof of Theorem \ref{thm1}.

Similarly as in the proof of the lower bound in Theorem 1 we see that
in the exceptional cases the extremal length of $\hat b$ equals zero. \hfill $\Box$

\bigskip

\noindent {\bf Proof of Lemma \ref{lemm1'}.} 
Let $\tau_3$ be the natural homomorphism from the braid group $\mathcal{B}_3$ to the symmetric group $\mathcal{S}_3$. When $\tau_3(b)$ is the identity then the braid $b$ is pure and the statement is clear.

Assume $\tau_3(b)$ is not the identity. We consider first the case when $\tau_3(b)$ is a transposition. If $\tau_3(b)=(13)$ then $b \Delta_3^{-1}$ is a pure braid, hence if $b$ is not a power of $\Delta_3$ it can written in the form \eqref{eq2'} with $k$ even and $\ell$ odd.
If $\tau_3(b)=(12)$ then $\sigma_1^{-1}b$ is a pure braid. Hence,  $b$
can be written in the form \eqref{eq2'} with $j=1$, $k$ odd and $\ell$ even. If $\tau_3(b)=(23)$ then $\sigma_2^{-1}b$ is a pure braid. Hence, $b$ can be written in the form \eqref{eq2'} with $j=2$, $k$ odd and $\ell$ even.

It remains to consider the case when $\tau_3(b)$ is a cycle. Suppose, $\tau_3(b)=(123)$. Since $(123)(13)=(12)$ the relation $\tau_3(b)(13)^{-1}=(12)$ holds. Hence, $\sigma_1^{-1}b \Delta_3^{-1}$ is a pure braid and $b$ has the form \eqref{eq2'} with $j=1$, and $k$ and $\ell$ odd. If $\tau_3(b)=(132)$ then, since $(132)(13)=(23)$ the braid  $\sigma_2^{-1}b \Delta_3^{-1}$ is a pure braid and $b$ can be written in the form \eqref{eq2'} with $j=2$, and $k$ and $\ell$ odd. \hfill $\Box$

\medskip

\noindent {\bf Proof of the lower bound in Theorem 3.}
Since $\Lambda(b_{tr})=\Lambda((b \Delta_3)_{tr})$ we may suppose that
$b=\sigma_j^{k} b_1$ for a pure braid $b_1$ which is a reduced word in $\sigma_1^2$ and $\sigma_2^2$. We may suppose that $k$ is an odd integer. The case when $k$ is even is contained in Theorem \ref{thm1}. Let $j=1$.
Consider a holomorphic mapping $g$ of a rectangle to $C_3(\mathbb{C)} \diagup \mathcal{S}_3$ that represents the braid and lift it to a mapping $\tilde g$ into $C_3(\mathbb{C)}$ so that the extension of the mapping to the closed rectangle takes the open lower side to  $\{(x_1,x_2,x_3) \in
\mathbb{R}^3: x_2<x_1<x_3\}$. Compose the mapping with the mapping $\mathfrak{C}$. We obtain a holomorphic mapping $\mathfrak{C}(\tilde g)$ of the rectangle to $\mathbb{C}\setminus \{-1,1\}$ whose continuous extension to the closure takes the open lower side of the rectangle to $(-\infty,-1)$ and the open upper side to $(-1,1)$.

Consider first the case when $|k|\geq 3$.
As in Lemma \ref{lemm14a} the restriction of the mapping $\mathfrak{C}(\tilde g)$ to the closure of an arbitrary maximal vertical segment in the rectangle is the union of the following two curves: an arc that is homotopic in $\mathbb{C}\setminus \{-1,1\}$ with endpoints in $\mathbb{R}\setminus \{-1,1\}$ to a half-circle in the lower half-plane if $k>0$, or in the upper half-plane if $k<0$, which joins a point in $(-\infty,-1)$ with a point in $(-1,1)$, and a curve representing the element $\mathfrak{C}_*(\vartheta(b)_{tr}) \in \pi_1^{tr}$. Hence, as in Lemma \ref{lemm15} the rectangle contains a curvilinear rectangle such that the restriction of the mapping $\mathfrak{C}(\tilde g)$ to it represents $\mathfrak{C}_*(\vartheta(b)_{tr})$. Therefore $\Lambda(b_{tr}) \geq \Lambda(\vartheta(b)_{tr})$. If $\vartheta(b)\neq \sigma_1^{2k'} $ for an integer $k'$, in other words, if $\vartheta(b)$ is not among the exceptional cases of
Theorem \ref{thm1} the lower bound holds.

Suppose $|k| =1$. If $b\neq \sigma_1^{\pm 1} $
then as in Lemma \ref{lemm14a} the restriction of the mapping $\mathfrak{C}(\tilde g)$ to the closure of each maximal vertical segment in the rectangle is the union of an arc that is homotopic in $\mathbb{C}\setminus \{-1,1\}$ with endpoints in $(-\infty,-1) \cup i \mathbb{R}$ to a quarter-circle in the upper or lower half-plane which joins a point in $(-\infty,-1)$ with a point in $i\mathbb{R}$, and a curve that represents $\mathfrak{C}_*(_{pb}\vartheta(b)_{tr})$.  As in Lemma \ref{lemm15} the rectangle contains a curvilinear rectangle such that the restriction of the mapping $\mathfrak{C}(\tilde g)$ to it represents $\mathfrak{C}_*(_{pb}\vartheta(b)_{tr})$. The lower bound
for this case follows from Theorem $1'$.

The proof
for the case $b=\sigma_2^k b_1$ (with $b_1$ a reduced word in $\sigma_1^2$ and $\sigma_2^2$)  is similar and is left to the reader. In this case the lift to $C_3(\mathbb{C})$ of a representing mapping for $b_{tr}$ is chosen with initial point in $\{x_1<x_3<x_2\}$.

The lower bound of Theorem \ref{thm3} is proved. \hfill $\Box$

 \medskip

\noindent {\bf Remark.} {\it For each conjugacy class $\hat w$ of elements of $\pi_1$ with at least two syllables and any syllable reduced representative $w \in  \pi_1$ of $\hat w$ the inequality
\begin{align}\label{eq14}
\Lambda(w_{pb}) \leq \Lambda(\hat w)
\end{align}
holds, but the inequality may be strict.}

\medskip

Indeed, let $\hat g: A \to \mathbb{C} \setminus \{-1,1\}$ be a locally conformal mapping representing $\hat w$.

As in the proof of Theorem \ref{thm2} there is a syllable reduced word $w\in \hat w$ and a connected component $L_0^0$ of the level set $L_0= \{z \in A: \mbox{Re} \, g(z)=0\}$ that joins the two boundary circles of $A$ such that $R \stackrel{def}{=}A\setminus L_0^0$ is a curvilinear rectangle and the restriction of $\hat g$ to $R$ represents the syllable reduced word $w$ with $pb$ boundary values.

Denote by $\Gamma_A$ the set
of loops in $A$ representing the positively oriented generator of the
fundamental group of
$A$ and by $\Gamma_R$ the set of arcs in $R$  whose extensions to the
boundary of $R$ join the lower side coming from $L_0^0$ with the respective upper side. Then the relation
$\Gamma_R<\Gamma_A$ holds, but not vice versa. (The elements of
$\Gamma_R$ may not close up to circles!)
We obtain \eqref{eq14}, but there are examples when the inequality is strict.

There is a statement analogous to the remark that concerns the extremal length of $w_{tr}$ (instead of $w_{pb}$). Take any conjugacy class $\hat w$ of elements of $\pi_1$ except conjugates of
$(a_1\,a_2)^n$ for an integer $n\neq 0$. Then each loop representing $\hat
w$ intersects the interval $(-1,1)$.
Let $\hat g: A \to \mathbb{C} \setminus \{-1,1\}$ be a locally conformal mapping on an annulus $A$ representing $\hat w$. As in the proof of Theorem \ref{thm2} there is a syllable reduced word $w\in \hat w$ and a connected component $L_0^0$ of the level set $L_0= \{z \in A: g(z) \in (-1,1)\}$ that joins the two boundary circles of $A$ such that $R \stackrel{def}{=}A\setminus L_0^0$ is a curvilinear rectangle and the restriction of $\hat g$ to $R$ represents the syllable reduced word $w$ with $tr$ boundary values.
Then
$$
\Lambda(\hat w) \geq \Lambda(w_{tr})
$$
and the inequality may be strict.

\medskip

\section{The upper bound for the extremal length. Plan of proof.}

Consider first the case of the extremal length with $pb$ boundary values.
Take a word  $w$ in the relative fundamental group $^{pb}\pi_1^{pb} $ and consider its lift under $f_1 \circ f_2$ to a homotopy class of slalom curves. We will represent this slalom class by a holomorphic mapping of a rectangle into $\mathbb{C}\setminus i\mathbb{Z}$. The extremal length of the rectangle provides an upper bound of the extremal length of $w$ with $pb$ boundary values.

For each syllable $\mathfrak{s}_j$ of $w$ (labeled from left to right) we represent the lift under $f_1\circ f_2$ of $_{pb}(\mathfrak{s}_j)_{pb}$ by a holomorphic mapping $\overset{\circ}{g_j}$ of a rectangle $\overset{\circ}{R_j}$ into $\mathbb{C} \setminus i\mathbb{Z}$. The rectangle $\overset{\circ}{R_j}$ is chosen as always with sides parallel to the axes.
We choose two points $\overset{\circ}{p_j}^-$
and $\overset{\circ}{p_j}^+$
in the boundary of $\overset{\circ}{R_j}$ so that the values of $\overset{\circ}{g_j}$ at these points are imaginary half-integers (but not imaginary integers). The rectangle $\overset{\circ}{R_j}$ will be normalized as follows. If a standard representative of the syllable, or a standard representative of the first term of the syllable, respectively, is contained in the right half-plane, then the derivative $\overset{\circ}{g_j}'$ at $\overset{\circ}{p_j}^-$ maps $i$ to the positive unit vector. It maps $i$ to the negative unit vector in the remaining case. Similarly, the derivative $\overset{\circ}{g_j}'$ at $\overset{\circ}{p_j}^+$ maps $i$ to the negative unit vector if a standard representative of the syllable, or a standard representative of the last term of the syllable, respectively, is contained in the right half-plane, and maps $i$ to the positive unit vector in the remaining case.

Each rectangle $\overset{\circ}{R_j}$ will be shifted by a complex number $c_j$ so
that $\overset{\circ}{p_j}^+ + c_j = \overset{\circ}{ p_{j+1}}^- + c_{j+1}$.
Put ${p_j}^{\pm}=\overset{\circ}{p_j}^{\pm} + c_j$ and $R_j=\overset{\circ}{R_j}+c_j$. Consider the
"shifted" functions  $ g_j$, $ g_j(z+c_j)= \overset{\circ}{g_j}(z)+i\, m_j$. The integers $m_j$ are chosen so that $g_j(p_j^+)= g_{j+1}(p_{j+1}^-)$. Since for each $j$ a standard representative of the (last term of the) syllable $\mathfrak{s}_j$ and a standard representative of the (first term of the) syllable $\mathfrak{s}_{j+1}$ are in different half-planes the derivatives of $g_j$ and $ g_{j+1}$ at $p_j^+ = p_{j+1}^-$ coincide. The shifted functions $g_j$ will be patched together in a quasiconformal way using an estimate for the second derivatives of the $g_j$. The obtained function is denoted by $g$. The union of the closures of the rectangles $\bar R_j $  will contain a curvilinear rectangle of the form $R_{J,\Phi, \frac{1}{18}}$.
Here for a real $C^1$-function $\Phi$ on a given open interval $J$ and a positive number $\textsf{b}$ we consider the curvilinear rectangle which is defined as  $R_{J,\Phi, \textsf{b}}\stackrel{def}{=} \{x+iy \in \mathbb{C}: y \in J,\, x \in (\Phi(y), \Phi(y)+\textsf{b})\}$. The construction produces a quasiconformal mapping of $R_{J, \Phi, \frac{1}{18}}$ into $\mathbb{C} \setminus i\mathbb{Z}$ that represents the lift of $_{pb}(w)_{pb}$ under $f_1 \circ f_2$. An upper bound  of the vertical length of each $R_j$, an estimate of the derivative of $\Phi$ and of the quasiconformal dilatation of the function $g$, and the following lemma will provide an upper bound of the extremal length of $R_{J,\Phi, \frac{1}{18}}$, and, hence of $\Lambda_{p}(w)$.

\begin{lemm}\label{lemm216} Let $\Phi$ be a real $C^1$-function on an open interval $J$ and let $\textsf{b}$ be a positive number.
Denote by $\Gamma_{\Phi}$ the set of curves in the rectangle $R_{J,\Phi, \textsf{b}}= \{x+iy \in \mathbb{C}: y \in J,\, x \in (\Phi(y),
\Phi(y)+\textsf{b})\}$ which join the two horizontal sides.  Suppose the absolute value $|\Phi'|$ of the derivative of $\Phi$ is bounded by the constant $C$. Then
$$
\lambda(\Gamma_{\Phi}) \leq (1 + C^2) \lambda(\Gamma_0),
$$
where $ \Gamma_0$ is the family corresponding to the function $\Phi_0$ which is identically equal to zero.
\end{lemm}

\noindent {\bf Proof}. The proof is similar to Example 1 in chapter 1 of \cite{A1}. For any measurable function $\varrho$ on $\mathbb{C}$ and any  $x \in (0,\textsf{b})$ we have
$$
\int_J \varrho(x + \Phi(y) +i y) \sqrt{1 +\Phi'(y)^2}dy \geq L_{\Gamma_{\Phi}}(\varrho).
$$
Integrate over the interval $(0,\textsf{b})$ and apply Fubini's Theorem and H\"older's inequality. Using the bound for $|\Phi'|$ we obtain
$$
(\int \int_{R_{\Phi, \textsf{b}}} dm_2 \int \int_{R_{\Phi, \textsf{b}}} \varrho^2 \cdot(1+ C^2) dm_2)^{\frac{1}{2}}\geq \textsf{b} L_{\Gamma_{\Phi}}(\varrho).
$$
Denote by $|J|$ the length of the interval $J$. We obtain
$$
\textsf{b} |J| (1+ C^2) A(\varrho) \geq \textsf{b}^2 L_{\Gamma_{\Phi}}(\varrho)^2.
$$
Hence,
$$
\frac{L_{\Gamma_{\Phi}}(\varrho)^2}{A(\varrho)} \leq \frac{|J|}{\textsf{b}}\cdot(1+C^2) = \lambda(\Gamma_0) (1+C^2).
$$
Taking the supremum over all measurable functions $\varrho$ with finite non-vanishing integral we obtain
$$
\lambda(\Gamma_{\Phi}) \leq (1+C^2)  \lambda(\Gamma_0).
$$
The lemma is proved. \hfill $\Box$

\medskip
The key for the construction of the rectangles $R_j$ and the mappings $g_j$ is the elliptic integral $\mathcal{F}_M$ (see equality \eqref{eq3} in Section 3) for integers and half-integers $M$. Recall that we take the branch of the square root which is positive on the positive real half-axis.  The mapping  $\mathcal{F}_M$  takes the left half-plane to a rectangle $R(M)$ contained in the left half-plane with sides parallel to  the axes. The mapping extends to a homeomorphism of the union of the closed half plane with $\infty$ which takes the interval $[-i(M+1),-iM]$ to the lower side 
and $[iM,i(M+1)]$ to the upper side of the rectangle, and it takes the interval $[-iM,iM]$ into the imaginary axis. Further, $0$ is mapped to $0$ which is the midpoint of one vertical side of the rectangle, and $\infty$ is mapped to the midpoint of the other vertical side. Moreover,
$$
\mathcal{F}_M (-z) = - \mathcal{F}_M (z)\; \mbox{for} \, z \in (-iM,iM).
$$
By the reflection principle, $\mathcal{F}_M$ extends holomorphically across each of the intervals $(-i(M+1),-iM)$, $(-iM,iM)$, $(iM,i(M+1))$, and across the complement of $[-i(M+1),i(M+1)]$ in $\mathbb{R} \cup \infty$, to a mapping from the right half-plane to the respective reflected rectangle. We will often denote the extension of $\mathcal{F}_M$ also by the same letter $\mathcal{F}_M$.

The integral, given by \eqref{eq3} defines also a conformal mapping from the right half-plane onto a rectangle contained in the right half-plane. Again, we use the branch of the square root which is positive on the positive real half-axis. The rectangle is the reflection in the imaginary axis of the rectangle obtained in the previous case. We denote the mapping given by \eqref{eq3} on the right half-plane by $\mathcal{F}_M^r$.
Notice that $\mathcal{F}_M^r$ is obtained by analytic continuation of  $\mathcal{F}_M$ across $(-iM,iM)$. 

The following lemma provides an estimate of the horizontal side length of the rectangle $R(M)$.

\begin{lemm}\label{lemm103} For $M \geq \frac{1}{2}$ the following
inequalities hold
\begin{equation}\label{eq106}
\frac{\pi}{2(M+1)}\leq |\mathcal{F}_M(i (M+1))-\mathcal{F}_M(i
M)|\leq\frac{\pi}{2 M}.
\end{equation}
\end{lemm}

\noindent {\bf Proof}. For $X  \in (0,1)$ we have
\begin{align}\label{eq108}\nonumber
\mathcal{F}_M(i (M+X))-\mathcal{F}_M(i M)& = \int_{iM}^{i(M+X)}
\frac{d\zeta}{\sqrt{(iM-\zeta)(iM+\zeta)(i(M+1)+\zeta)(i(M+1)-\zeta)}}\\
\nonumber
& =i\int_M^{M+X} \frac{dw}{\sqrt{(M-w)(M+1-w)(M+w)(M+1+w)}}\\
& =i \int_0^X \frac{dx}{\sqrt{x(x-1)c(M,x)}}.
\end{align}
Here $c(M, x) \overset{\mbox{def}} = (2M + x)(2M + 1 + x)$.
Integration is along the interval $(0,X) \subset (0,1)$ and we have
chosen the branch of the square root which is positive on the positive
part of the real axis. The term $c(M, x)$  is
positive and its square root is contained in the interval $(2M, 2M +
2)$.
With the choice of the branch of the square root we have
$\sqrt{x(x-1)}=	-i\sqrt{x(1-x)},\, x \in (0,1)$.

For $X=1$ we obtain
$$
\mathcal{F}_M(i (M+1))-\mathcal{F}_M(i M)= - \int_0^1
\frac{dx}{\sqrt{x(1-x)c(M,x)}}
$$
Since
$$
\int_0^1
\frac{dx}{\sqrt{x(1-x)}}= 2 \int_0^{\frac{1}{2}}
\frac{dx}{\sqrt{x(1-x)}}= 2\int _0^{\sqrt{\frac{1}{2}}}\frac{2udu}{\sqrt{u^2(1-u^2)}}=
$$
$$
4\int_0^{\sqrt{\frac{1}{2}}}\frac{du}{\sqrt{1-u^2}}=
4\arcsin\sqrt{\frac{1}{2}}=\pi,
$$ 
we obtain \eqref{eq106}.    \hfill $\Box$
\medskip

The following lemma gives an estimate of the vertical side length $\frac{2}{i}\mathcal{F}_M(iM)$ of the rectangle $R(M)$ for the case $M \geq 2 $.

\begin{lemm}\label{lemm104}
For $M \geq 2$ the inequalities
\begin{equation}\label{eq104a}
\frac{1}{2} \log(2M+1) < \frac{M+1}{i}
 \mathcal{F}_M(iM)= K(\frac{M}{M+1}) < \frac{1}{2}\big( \frac{12}{5 } + \log(2M+1)\big)
\end{equation}
hold.
\end{lemm}

{\bf Proof}. 
Use the following estimate  for the denominator in \eqref{eq5a}
$$
(1-x)^2 < \sqrt{(1-x^2)(1- (\frac{M}{M+1})^2 x^2)}
< 1- (\frac{M}{M+1})^2 x^2
$$
with $x \in (-1,1)$.
The lower bound in \eqref{eq104a} is obtained as follows.
\begin{align}\label{eq105}
\frac{M+1}{i}\mathcal{F}_M(i M) &> \int_0^1
\frac{dx}{1- (\frac{M}{M+1})^2 x^2}= \frac{1}{2}\int _0^1
(\frac{dx}{1- \frac{M}{M+1}x} + \frac{dx}{1+ \frac{M}{M+1}x}) \nonumber\\
& =  \frac{1}{2\frac{M}{M+1}} (\log(1+\frac{M}{M+1}) - \log(1-\frac{M}{M+1})) =
\frac{M+1}{2M}\log(2M+1) \nonumber\\
& > \frac{1}{2} \log(2M+1).
\end{align}
For the estimate from above we write
$$
\frac{M+1}{i}\mathcal{F}_M(i M)= \int_0^{1-\frac{1}{M+1}} +
\int_{1-\frac{1}{M+1}}^1 \overset{\mbox{def}}=\mathcal{I}_1 +
\mathcal{I}_2.
$$
We obtain for the term $\mathcal{I}_1$
$$
\mathcal{I}_1 < \int_0^{1-\frac{1}{M+1}} \frac{dx}{1-x^2}= \frac{1}{2} \int_0^{1-\frac{1}{M+1}}(\frac{dx}{1-x} + \frac{dx}{1+x})=
$$
$$
\frac{1}{2} (-\log(1-(1-\frac{1}{M+1})) + \log(1+(1-\frac{1}{M+1})))=
$$
$$
\frac{1}{2}\log(2M+1)
$$
For the estimate of $\mathcal{I}_2$  we use the fact that for $1-\frac{1}{M+1} \leq x \leq 1$ the
denominator $\sqrt{(1-x^2)(1-(\frac{M}{M+1}x)^2)}$ is not smaller than
$$
\sqrt{1-x} \sqrt{2- \frac{1}{M+1}}\sqrt{1-(\frac{M}{M+1})^2}.
$$
Hence,
\begin{align}\label{eq105'}
\mathcal{I}_2 & \leq
\frac{\sqrt{(M+1)^3}}{\sqrt{2M+1} \sqrt{(M+1)^2-M^2}}\int_{1-\frac{1}{M+1}}^1 \frac{dx}{\sqrt{1-x}}\nonumber \\ = & \frac{\sqrt{(M+1)^3}}{2M+1} \cdot  2\cdot \sqrt{\frac{1}{M+1}}=
2 \cdot \frac{M+1}{2M+1} \leq 2 \cdot  \frac{3}{5}
\end{align}
for $M\geq 2$. 
For the last inequality we used that for $M \geq 2 $ the inequality $ \frac{M+1}{2M+1} \leq \frac{3}{5}$ holds.

The lemma is proved.  \hfill $\Box$
\medskip

The inverse mappings of the versions of the elliptic integral represent
classes of curves as follows.
The integral \eqref{eq3} which defines $\mathcal{F}_M(z)$  with $z \in \mathbb{C}_{\ell}$
is the inverse of a conformal mapping of a rectangle $R(M)$ representing the
class $\gamma^*_{\ell,M}$. The class of curves $\gamma^*_{\ell,M,-}$ obtained from $\gamma^*_{\ell,M}$ by inverting the orientation is represented by
the inverse of the mapping $z \to
\mathcal{F}^-_M(z) \stackrel{def}{=}
-\mathcal{F}_M(z) - |\mathcal{F}_M(i(M+1))-\mathcal{F}_M(iM)|,\, z \in \mathbb{C}_{\ell}\,,$ which takes the left half-plane to the same rectangle $R(M)$.  (Here $|\mathcal{F}_M(i(M+1)-\mathcal{F}_M(iM)|$ is the horizontal side length of the rectangle $R(M) = \mathcal{F}_M(\mathbb{C}_{\ell})$.)

Consider the function $\mathcal{F}_M^r$ represented by the integral  \eqref{eq3} with $z \in \mathbb{C}_{r}$. The inverse $(\mathcal{F}_M^r)^{-1}$ of $\mathcal{F}_M^r$ is a conformal mapping onto the right half-plane from a rectangle $R(M)^r$ in the right half-plane with one vertical side on the imaginary axis. The mapping  $(\mathcal{F}_M^r)^{-1}$ represents the class $\gamma^*_{r,M,-}$ which is obtained from   $\gamma^*_{r,M}$ by inverting orientation. The class $\gamma^*_{r,M}$ is represented by the inverse of the mapping
$z \to
\mathcal{F}^{r,-}_M(z) \stackrel{def}{=}
-\mathcal{F}^r_M(z) + |\mathcal{F}_M(i(M+1)-\mathcal{F}_M(iM)|,\, z \in \mathbb{C}_{r}\,,$
which is defined on the same rectangle $R(M)^r$.

\section{Holomorphic maps representing lifts of syllables with ${pb}$ boundary values}

\medskip

\noindent {\bf 1. Syllables of form (1) with $pb$ boundary values of degree at least $5$.}

\noindent We consider first syllables of the form $a^n$ with $pb$ boundary values where $a$ is one of the standard generators or its inverse and $n\geq 5$. Recall that $d=n$ is the degree of the syllable in this case.
Assume $a=a_1$. All other cases of $a$ are treated in the same way. According to Proposition \ref{prop4b}    
the class of representatives of the syllable with $pb$ boundary values lifts to 
the class $\gamma_{M,\ell}^*+iM$ with
$M= \frac{d-1}{2}$. Notice that with the choice of $n$ the parameter $M$ is at least $2$.
Let $\mathcal{F}_M(z), z \in \mathbb{C}_{\ell}$, be the mapping \eqref{eq3} defined in Section 3
and let $\mathcal{F}_M^{-1}$ be its inverse.
The derivative  $\mathcal{F}_M'(\pm i(M+\frac{1}{2}))$  is equal to

\begin{align}\label{eq301}
\mathcal{F}_M'(\pm i(M+\frac{1}{2})) & = \frac{1}{\sqrt{((i(M+\frac{1}{2}))^2-(iM)^2)((i(M+\frac{1}{2}))^2-(i(M+1))^2)}}\nonumber\\
& =\frac{1}{\sqrt{((M+\frac{1}{2})^2-M^2)((M+\frac{1}{2})^2-(M+1)^2)}}\nonumber\\
& = \frac{\pm i}{\sqrt{(M+\frac{1}{4})(M+\frac{3}{4})}}\;.
\end{align}

Put
\begin{align}\label{eq301a}
r_{(1),M}  & \stackrel{def}{=}  \sqrt{(M+\frac{1}{4})(M+\frac{3}{4})}\,,\nonumber\\
R_{(1),M} & \stackrel{def}{=} r_{(1),M}\, R(M)\,,\nonumber\\
f_{(1),M} &  \stackrel{def}{=} r_{(1),M} \, \mathcal{F}_M\;\;\;\;\;, \mbox{and}\nonumber\\
g_{(1),M} & \stackrel{def}{=}(f_{(1),M})^{-1}.
\end{align}
Let $\xi^{\pm}_M$ be the point $f_{(1),M}(\pm i(M+\frac{1}{2}))$ on the upper (lower, respectively) horizontal side of the rectangle. Then
$(g_{(1),M})'(\xi^{\pm}_M)=\frac{1}{(f_{(1),M})'(\pm i(M+\frac{1}{2}))}=\mp i\,.$ Thus $g_{(1),M}'(\xi^+_M)$ maps $i$ to the unit vector in positive direction and  $g_{(1),M}'(\xi^-)$ maps $i$ to the unit vector in negative direction. Moreover, $g_{(1),M}$ maps $R_{(1),M}$ onto $\mathbb{C}_{\ell}$ so that the upper horizontal side is mapped to $(iM,i(M+1))$, the lower  horizontal side is mapped to $(-i(M+1),iM)\,,$ and the points $\xi^{\pm}_M$ are mapped to $\pm i (M+\frac{1}{2})$.

The following lemma holds.

\begin{lemm}\label{lemm217} For $\zeta \in \mathbb{C}$, $|\zeta \pm i(M+\frac{1}{2})|< \frac{1.03 \sqrt{2}}{18}$, the inequality
$$
|\frac{1}{(f_{(1),M})'(\zeta)}|< 1.03
$$
holds.
\end{lemm}

\noindent \textbf{Proof.} Put $\zeta=i(M+\frac{1}{2})+z$ with $|z|< \frac{1.03 \sqrt{2}}{18} $. For the function $\mathcal{F}_M$ (see \eqref{eq3}) we have

\begin{align}\label{eq221}
&(\frac{1}{\mathcal{F}'( i(M+\frac{1}{2}) +z   )})^2
= ((M+\frac{1}{2}-iz)^2-M^2)((M+\frac{1}{2}-iz)^2-(M+1)^2)\nonumber\\
=& (M+\frac{1}{2} -iz-M)(M+\frac{1}{2} -iz+M)(M+\frac{1}{2} -iz-(M+1))(M+\frac{1}{2} -iz+(M+1))\nonumber \\
=&(\frac{1}{2}-iz)(-\frac{1}{2}-iz)\big((2M+\frac{1}{2})-iz\big)\big((2M+\frac{3}{2})
-iz\big).
\end{align}

We obtained
\begin{equation}\label{eq222}
\frac{1}{\mathcal{F_M}'(i(M+\frac{1}{2}) +z)}= -i \sqrt{(\frac{1}{4}+z^2)\big(2M+\frac{1}{2}-iz\big)\big(2M+\frac{3}{2}
-iz\big)},
\end{equation}
hence, for $|z|< \frac{1.03 \sqrt{2}}{18}$
\begin{align}\label{eq223}
&|\frac{1}{(f_{(1),M})'(i(M+\frac{1}{2})+z)}| \nonumber\\
&\leq \frac{\sqrt{\frac{1}{4}+2 \frac{1.03^2}{18^2}}\sqrt{2M+\frac{1}{2} +\frac{1.03 \cdot \sqrt{2}}{18}}\sqrt{ 2M+\frac{3}{2}+\frac{1.03 \cdot \sqrt{2}}{18}}}{\sqrt{(M+\frac{1}{4})(M+\frac{3}{4})}}.
\end{align}

The expression on the right hand side of \eqref{eq223} is decreasing in $M$, hence for $M\geq 2$ it does not exceed its value at $2$, i.e.
\begin{align}\label{eq224}
|\frac{1}{(f_{(1),M})'(i(M+\frac{1}{2})+z)}| & \leq 4 \frac{\sqrt{\frac{1}{4} +2 \frac{1.03^2}{18^2}}\sqrt{\frac{9}{2}+ \frac{1.03 \sqrt{2}}{18}}\sqrt{\frac{11}{2}+\frac{1.03 \sqrt{2}}{18}}}{\sqrt{99}}\nonumber\\
&<1.0296<1.03
\,.
\end{align}
The estimate for $i(M+\frac{1}{2})$ replaced by $-i(M+\frac{1}{2})$ is the same. \hfill $\Box$

\begin{cor}\label{cor203} The inverse $g_{(1),M}=(f_{(1),M})^{-1}$ of $f_{(1),M}$ maps the sets $E_{(1),M}^{\pm}\stackrel{def}{=}\{|\xi \mp \xi_M^{\pm}|< \frac{\sqrt{2}}{18}\}$ into the sets $\{|\zeta \mp i(M+\frac{1}{2})|< \frac{1.03 \cdot \sqrt{2}}{18}\}$.
The following estimate holds for the derivative of the inverse
\begin{equation}\label{eq224a}
|(g_{(1),M})'(\xi)|< 1.03
\end{equation}
for $|\xi- \xi_M^{\pm}|< \frac{\sqrt{2}}{18}.$
\end{cor}

\noindent{\textbf{Proof.}} 
If for $\zeta= g_{(1),M}(\xi)$ the inequality  $|\zeta \mp i(M+\frac{1}{2})|< \frac{1.03 \sqrt{2}}{18} $    holds we have $|(g_{(1),M})'(\xi)|=    |\frac{1}{(f_{(1),M})'((f_{(1),M})^{-1}(\xi))}|< 1.03$.
Suppose \eqref{eq224a} is not true for some $\xi, \,|\xi- \xi_M^{\pm}|< \frac{\sqrt{2}}{18}$. Since the inequality \eqref{eq224a} holds
at $\xi^{\pm}_M$ there is a point $\xi'_{\pm}\,, |\xi'_{\pm}- \xi_M^{\pm}|< \frac{\sqrt{2}}{18},$ of closest distance to $\xi^{\pm}_M$ among points for which
\eqref{eq224a} is violated. Then $|(g_{(1),M})'(\xi'_{\pm})|= 1.03$ and \eqref{eq224a} is satisfied for all $\xi$ with $|\xi- \xi_M^{\pm}|<|\xi'_{\pm}- \xi_M^{\pm}|$.
Write
 $g_{(1),M}(\xi'_{\pm}) - g_{(1),M}(\xi^{\pm}_M )= \int_{{\xi}^{\pm}_M}^{\xi'_{\pm}} (g_{(1),M})'(\xi)d\xi$ where integration in the last integral is along the straight line segment $[\xi^{\pm}_M,\xi'_{\pm}]$ joining the two points.
 The estimate of the length of the segment and the fact that \eqref{eq224a} is satisfied on the open segment of integration yields $|g_{(1),M}(\xi'_{\pm}) - g_{(1),M}(\xi^{\pm}_M )|< 1.03 \cdot \frac{\sqrt{2}}{18}$ which contradicts
 the assumption. The corollary is proved.   \hfill $\Box$

\medskip

The corollary implies that $g_{(1),M}$ maps the sets $E_{(1),M}^{\pm}{=}\{\xi \in \mathbb{C}: |\xi - \xi_M^{\pm}|< \frac{\sqrt{2}}{18}\}$ into $\mathbb{C} \setminus i \mathbb{Z}$.
Indeed, the image $g_{(1),M}$ of a point $\xi \in E_{(1),M}^{\pm}$ has distance less than $\frac{1.03 \sqrt{2}}{18} <\frac{1}{2}$ from $\pm i(M+\frac{1}{2})$.

\medskip
To estimate the second derivative of $g_{(1),M}$ we need the following simple
lemma.

\begin{lemm}\label{lemm218} For any conformal mapping $f$ on a small disc $\Delta$ and its inverse $f^{-1}$ the following formula holds
$$
(f^{-1})''(f(z))= \frac{1}{2}\big(\frac{1}{f'(z)^2}\big)'\,.
$$
\end{lemm}

\noindent {\bf Proof}. Since $(f^{-1})'(f(z))= \frac{1}{f'(z)}$, we obtain for the second derivative
$$(f^{-1})''(f(z)) \cdot f'(z) = (\frac{1}{f'(z)})'.
$$
Hence
$$
(f^{-1})''(f(z))= \frac{1}{f'(z)} \cdot (\frac{1}{f'(z)})' = \frac{1}{2}\big(\frac{1}{f'(z)^2}\big)'\,.
$$
\hfill $\Box$

\medskip
\begin{lemm}\label{lemm219} For $\xi \in E_{(1),M}^{\pm}$
the following inequality holds for the second derivative of $g_{(1),M}$
\begin{equation}\label{eq225}
|(g_{(1),M})''(\xi)| < 2.75 .
\end{equation}
\end{lemm}

\medskip

\noindent {\bf Proof.}
Under the condition of the lemma $\xi = f_{(1),M}(\pm i(M+\frac{1}{2}) +z)$ with $|z|< \frac{1.03 \, \sqrt{2}}{18}$. By lemma \ref{lemm218} (see also \eqref{eq221})
\begin{align}
&(g_{(1),M})''(\xi)= ((f_{(1),M})^{-1})''(f_{(1),M}(\pm i(M+\frac{1}{2}) +z))) \nonumber\\
&=\frac{1}{2}\frac{1}{(r_{(1),M})^2} \cdot\big(((M+\frac{1}{2}\mp iz)^2-M^2)
((M+\frac{1}{2}\mp iz)^2-(M+1)^2)\big)'\nonumber\\
&= \frac{1}{2}\frac{1}{(r_{(1),M})^2} \cdot \frac{1}{((\mathcal{F}_M)'(\pm i(M+\frac{1}{2})+z))^2}
\cdot \big(\frac{\mp i}{\frac{1}{2}\mp iz}+\frac{\pm i}{\frac{1}{2}\mp iz}-\frac{\pm i}
{2M+\frac{1}{2}\mp iz}-\frac{\pm i}{2M+\frac{3}{2}\mp iz}\big)\,. \nonumber
\end{align}
Hence, \eqref{eq224a} implies for $M\geq 2$

\begin{align}
|(g_{(1),M})''(\xi)|\leq \frac{1}{2} 1.03^2 \big(\frac{2}{\frac{1}{2}-\frac{1.03 \sqrt{2}}{18}} + \frac{1}{\frac{9}{2}-\frac{1.03 \sqrt{2}}{18}} + \frac{1}{\frac{11}{2}-\frac{1.03 \sqrt{2}}{18}}\big) <2.75.
\end{align}
The lemma is proved. \hfill $\Box$

\medskip

The relations \eqref{eq104a} and \eqref{eq301a} imply the following upper bound for the
vertical side length $\mbox{vsl}(R_{(1),M})$  of the normalized rectangle $R_{(1),M}=r_{(1),M} R(M)$ with $M=\frac{d-1}{2}\geq 2$
\begin{align}\label{eq302b}
\mbox{vsl}(R_{(1),M})& \leq
({\frac{12}{5}+ \log(2M+1)})\frac{\sqrt{(M+\frac{1}{4})(M+\frac{3}{4})}}{{M+1}}  < \frac{12}{5}-\log(4-\frac{1}{5}) + \log((4-\frac{1}{5})d) \nonumber \\
& < 0.362 \cdot \log 19 + \log(4d-1) < 1.362 \cdot \log(4d-1).
\end{align}
Here $d=2M+1\geq 5$, $(4-\frac{1}{5})d<4d-1$ for  $d \geq 5$. We used the inequality  ${\frac{12}{5}-\log(4-\frac{1}{5})} < 0.362 \cdot \log 19$.

\bigskip
\bigskip

\noindent {\bf 2. Syllables of form (2) with $pb$ boundary values of degree at least $5$}.

\noindent Consider first syllables of degree at least $5$. Without loss of generality we consider the syllables $(a_2^{-1} a_1^{-1})^k a_2^{-1}$ or $(a_2^{-1} a_1^{-1})^{k+1}$  for $k\geq 2$ with $pb$ boundary values. The degree $d$ is equal to $2|k|+1$ or $2|k|+2$ respectively, and is at least $5$. The other syllables of form (2) of degree at least $5$ are treated in the same way.   The lift under $f_1 \circ f_2$ of a suitable representative is contained in the punctured unit disc $\mathbb{D}\setminus \{ 0\}$, has endpoints on the imaginary axis and makes $d$ half-turns around $0$. This lift is a slalom curve which is the composition of $d \geq 5$ trivial elementary slalom curves.

We construct now a holomorphic mapping on a curvilinear rectangle which represents the homotopy class of this slalom curve. Take $M= \frac{d+1}{2}$. Restrict $\mathcal{F}_M$ to the domain $\{z \in \mathbb{C}_{\ell}: |\mbox{Im}\, z| < M-\frac{1}{2}\}$. Recall that $\mathcal{F}_M$ maps $\mathbb{C}_{\ell}$ onto a rectangle $R(M)$ whose right vertical side is contained in the imaginary axis and has midpoint $0$. The extension of
$\mathcal{F}_M$ takes the point $\infty$ to the midpoint of the left vertical side of $R(M)$. Each of the two lines  $\{z \in \mathbb{C}_{\ell}:\mbox{Im}\, z = \pm (M-\frac{1}{2})\}$ is mapped to a curve in $R(M)$ with one endpoint being a point on the right vertical side. 
The other endpoint is the image of $\infty$. Put $X_M\stackrel{def}{=}\mathcal{F}_M(\{z \in \mathbb{C}_{\ell}: |\mbox{Im}\, z| < M-\frac{1}{2}\})$. (See also Figure 7 below.)

Let $r_{(2),M}$ be a positive constant that will be specified later. Put

\begin{align}\label{eq302c}
R_{(2),M} & \stackrel{def}{=} r_{(2),M} R(M) ,\; \tilde R_{(2),M} \stackrel{def}{=} r_{(2),M} X_M \cap \{- \frac{1}{18}<\mbox {Re}\, \xi < 0\} , \nonumber\\
f_{(2),M} &  \stackrel{def}{=} r_{(2),M}\cdot \mathcal{F}_M \;\;\;\;\; \mbox{and}\nonumber\\
g_{(2),M}(\xi) & \stackrel{def}{=} \exp(\pi((f_{(2),M})^{-1}(\xi)+i (M-1))).
\end{align}
Note that with these definitions  $r_{(2),M} X_M=f_{(2),M}(\{z \in \mathbb{C}_{\ell}: |\mbox{Im}\, z| < M-\frac{1}{2}\})$.

(If e.g. the syllable equals 
$(a_1 a_2)^k a_1$ or $(a_1 a_2)^{k+1}, \, k\geq 2$, we take a rectangle in the right half-plane defined similarly as $ R_{(2),M}$ using $f_{(2),M}^r   \stackrel{def}{=} r_{(2),M}\cdot\mathcal{F}_M^r$, and define 
$g_{(2),M}^r(\xi)  \stackrel{def}{=} \exp(\pi(({f_{(2),M}}^r)^{-1}(\xi)+i (M-1)))$.)

\begin{lemm}\label{lemm18} The set $\tilde R_{(2),M}$ is a curvilinear rectangle with horizontal sides $f_{(2),M}(\{z \in \mathbb{C}_{\ell}: \mbox{Im}\, z = \pm  (M-\frac{1}{2})\}) \cap \{- \frac{1}{18}<\mbox {Re} \,\xi < 0\} $ being graphs over the a segment of the real axis of functions with absolute value of the derivative not exceeding $0.05$.
\end{lemm}
We postpone the proof of the lemma.

Restrict the mapping $(f_{(2),M})^{-1}$ to the curvilinear rectangle $\tilde R_{(2),M}$.

The mapping $\zeta \to \exp(\pi (\zeta + i(M-1)))$ takes the point $-i(M-\frac{1}{2})$ to $-i$. For each $x \in (-\infty,0)$ the curve $y \to \exp(\pi (x + i y + i(M-1)))\,, y \in (-(M-\frac{1}{2 }), (M-\frac{1}{2 }))\,,$ runs along the circle of radius $e^x<1$ and center $0$. It has initial point in $(-i,0)$ and makes $d=2M-1$ half-turns around zero.
Hence, the composition $\xi \to g_{(2),M}(\xi)=   \exp(\pi (f_{(2),M})^{-1}(\xi) + i(M-1)),\, \xi \in \tilde R_{(2),M},\,$ represents a lift of the syllable.
The image $g_{(2),M}(\tilde R_{(2),M})$ is contained in the punctured unit disc.

Put $\eta^{\pm} \stackrel{def}{=}f_{(2),M}(\pm i(M-\frac{1}{2}))$.
Note that
\begin{align}\label{eq302d}
(g_{(2),M})'(f_{(2),M}(\zeta))& = \pi \exp(\pi ((f_{(2),M})^{-1}(f_{(2),M}(\zeta)) +i (M-1))) \cdot ((f_{(2),M})^{-1})'(f_{(2),M}(\zeta))
\nonumber\\
& = \pi \exp(\pi (\zeta +i (M-1))) \cdot \frac{1}{(f_{(2),M})'(\zeta)}\;.
\end{align}
The derivative of $g_{(2),M}$ at $\eta^{\pm}$ equals
\begin{align}
(g_{(2),M})'(\eta^{\pm})& = \pi \exp(\pm \pi i((M-\frac{1}{2})+(M-1)))\cdot \frac{1}{r_{(2),M} (\mathcal{F}_M)'(\pm i (M-\frac{1}{2}))}\;.\nonumber
\end{align}
Here
\begin{align}
\frac{1}{ (\mathcal{F}_M)'(\pm i (M-\frac{1}{2}))} & =  \sqrt{((M-\frac{1}{2})^2-M^2)((M-\frac{1}{2})^2-(M+1)^2)} \nonumber\\
& = \sqrt{3} \sqrt{(M-\frac{1}{4})(M+\frac{1}{4})  }. \nonumber
\end{align}
Put
\begin{equation}\label{eq302e}
r_{(2),M}= \pi \sqrt{3} \sqrt{(M-\frac{1}{4})(M+\frac{1}{4})  }.
\end{equation}
Then, see \eqref{eq302c}, $ (g_{(2),M})'(\eta^-)= - i$ and $(g_{(2),M})'(\eta^+)$ equals $i$ if $d$ is odd and equals $-i$ if $d$ is even.

\begin{lemm}\label{lemm218a} For $\zeta \in \mathbb{C}$, $|\zeta\mp i(M-\frac{1}{2})|< \frac{0.4 \sqrt{2}}{18} $, the inequality
$$
\frac{1}{|(f_{(2),M})'(\zeta)|}< 0.3343
$$
holds.
\end{lemm}

\medskip

\noindent \textbf{Proof.}
Put here and below $\zeta= \pm i(M-\frac{1}{2}) +z$. Then
\begin{align}\label{eq25}
&\frac{1}{\mathcal{F}_M'(\zeta)^2}= \nonumber \\
&((M-\frac{1}{2}\mp iz)^2-M^2)
((M-\frac{1}{2}\mp iz)^2-(M+1)^2)= \nonumber\\
&(-\frac{1}{2}\mp iz)(-\frac{3}{2}\mp iz)(2M-\frac{1}{2}\mp iz)
(2M+\frac{1}{2}\mp iz).
\end{align}
Hence for $|z|<\frac{0.4 \sqrt{2}}{18}$
\begin{equation}\label{eq26}
\frac{1}{|(f_{(2),M})'(\zeta)|}<\frac{2\sqrt{\frac{1}{2}+
\frac{0.4 \sqrt{2}}{18}}
\sqrt{\frac{3}{2}+\frac{0.4 \sqrt{2}}{18}}
\sqrt{M+\frac{1}{4}+\frac{0.2 \sqrt{2}}{18}}
\sqrt{M-\frac{1}{4}+\frac{0.2 \sqrt{2}}{18}}}{\pi \sqrt{3}\sqrt{(M+\frac{1}{4})(M-\frac{1}{4})}}.
\end{equation}
Since the right hand side is decreasing in $M$,  the expression on the right for $M\geq 2$ does not exceed the expression for $M=2$, i.e. for $M\geq 2$
$$
\frac{1}{|(f_{(2),M})'(\zeta)|}<\frac{8\sqrt{\frac{1}{2}+\frac{0.4 \sqrt{2}}{18}}
\sqrt{\frac{3}{2}+\frac{0.4 \sqrt{2}}{18}}
\sqrt{\frac{9}{4}+\frac{0.2 \sqrt{2}}{18}}
\sqrt{\frac{7}{4}+\frac{0.2 \sqrt{2}}{18}}}{\pi \sqrt{3}\sqrt{63}}< 0.3343.
$$
\hfill $\Box$

The proof of the following corollary is the same as the proof of corollary \ref{cor203}.
\begin{cor}\label{cor204} The inverse $(f_{(2),M})^{-1}$ of $f_{(2),M}$ maps the sets $E_{(2),M}^{\pm}\stackrel{def}{=}\{|\xi-\eta_M^{\pm}|< \frac{\sqrt{2}}{18}\}$ into the sets $\{|\zeta\mp i(M-\frac{1}{2})|< \frac{0.3343 \cdot\sqrt{2}}{18}\}$. The derivative of the inverse satisfies the following inequality
\begin{equation}\label{eq26a}
|((f_{(2),M})^{-1})'(\xi)|<0.3343
\end{equation}
for $|\xi- \eta_M^{\pm}|<  \frac{\sqrt{2}}{18}.$
\end{cor}

The corollary says that the images $f_{(2),M}(\{|\zeta \mp i(M-\frac{1}{2})|< \frac{0.3343 \cdot\sqrt{2}}{18}\})$ cover the sets $\{|\xi-\eta_M^{\pm}|< \frac{\sqrt{2}}{18}\}$. Since $0.3343\cdot \frac{\sqrt{2}}{18}< 0.0263<0.03$,
Lemma \ref{lemm18} is an immediate consequence of the following lemma.

\begin{lemm}\label{lemm19} For $\zeta$ in a $\frac{3}{100}$-neighbourhood of $\pm i (M-\frac{1}{2})$ the (principal branch) of the argument of $f_{(2),M}'(\zeta)$ satisfies the inequality
$$
|\arg(f_{M,(2)}'(\zeta))|<\arctan(0.05).
$$
\end{lemm}

\medskip

\noindent \textbf{Proof of Lemma \ref{lemm19}.} We have
$$
\arg ( f_{M,(2)}'(\zeta))= \frac{1}{2}\arg (f_{M,(2)}'(\zeta))^2 =-\frac{1}{2}\arg ( f_{M,(2)}'(\zeta))^{-2} .
$$
Also,  $(f_{M,(2)}'(\zeta))^{-2}$ is the product of a real number and four complex numbers, hence its argument is the sum of the arguments of the factors provided all these arguments are small.
For $|z|<\frac{3}{100}$ and $M\geq 2$ we have
$$
|\arg(-\frac{1}{2}\mp iz)| \leq \arctan(\frac{3}{100} \frac{1}{\frac{1}{2} -\frac{3}{100}}) \leq 0.06375,
$$
$$|\arg(-\frac{3}{2}\mp iz)| \leq \arctan(\frac{3}{100} \frac{1}{\frac{3}{2} -\frac{3}{100}})\leq 0.02041,
$$
$$
|\arg(2M-\frac{1}{2} \mp iz)|\leq \arctan (\frac{3}{100} \frac{1}{  2M-\frac{1}{2} -\frac{3}{100}}) \leq \arctan(\frac{3}{100} \frac{1}{4             -\frac{1}{2} -\frac{3}{100}}) \leq 0.00865,
$$
and
$$
|\arg(2M+\frac{1}{2} \mp iz)|\leq \arctan(\frac{3}{100} \frac{1}{4             +\frac{1}{2} -\frac{3}{100}}) \leq 0.006712.
$$
The sum of the four numbers does not exceed $0.0996$. Hence,
$$|\arg (f_M'(\zeta)|< \frac{0.0996}{2} = 0.0498 < \arctan(0.05).$$
\hfill $\Box$
\begin{lemm}\label{lemm220} For $|\xi- \eta^{\pm}_M|< \frac{\sqrt{2}}{18}$ the following inequality holds for the second derivative of $g_{(2),M}$
\begin{equation}\label{eq27}
|(g_{(2),M})''(\xi)| < 1.863.
\end{equation}
\end{lemm}

\medskip
\noindent {\bf Proof}. Put $\xi= f_{(2),M} (\zeta)$. By formulas \eqref{eq302d} and \eqref{eq25} and by lemma \ref{lemm218}
\begin{align}\label{eq28}
&(g_{(2),M})''( f_{(2),M}(\zeta))=\nonumber\\
&  \pi^2 \exp(\pi((f_{(2),M})^{-1}(f_{(2),M}(\zeta))+i (M-1)))\cdot \frac{1}{((f_{(2),M})'(\zeta))^2}+  \nonumber\\
&  \pi \exp(\pi((f_{(2),M})^{-1}(f_{(2),M}(\zeta))+i (M-1)))\cdot ((f_{(2),M})^{-1})''(f_{(2),M}(\zeta))=\nonumber\\
&\pi^2 \exp(\pi ( \zeta +i (M-1))) (\frac{1}{(f_{(2),M})'(\zeta)})^2 + \pi \exp(\pi (\zeta +i (M-1)))\frac{1}{2} ((\frac{1}{(f_{(2),M})'(\zeta)})^2)' = \nonumber\\
&\pi^2 \exp(\pi( \zeta +i (M-1))) (\frac{1}{(f_{(2),M})'(\zeta)})^2 \cdot\nonumber\\
& \big(1+\frac{1}{2\pi} (\frac{i}{\frac{1}{2}\pm iz}+\frac{i}{\frac{3}{2}\pm iz}-\frac{i}
{2M-\frac{1}{2}\mp iz}-\frac{i}{2M+\frac{1}{2}\mp iz})\big)
\end{align}
If $|\xi- \eta^{\pm}_M|< \frac{\sqrt{2}}{18}$ then by Corollary \ref{cor204}
$$
|(f_{(2),M})^{-1}(\xi)-(f_{(2),M})^{-1}(\eta^{\pm}_M)|=
|\zeta\mp i(M-\frac{1}{2})|<0.3343 \frac{\sqrt{2}}{18}<0.4 \frac{\sqrt{2}}{18}.
$$
Hence, by \eqref{eq28} and Corollary \ref{cor204} for $M\geq 2$
\begin{align}
|(g_{(2),M})''(\xi)|  < & \pi^2 \exp(\frac{ 0.4 \sqrt{2}}{18} \pi)\cdot (0.3343)^2 \cdot\nonumber \\
&(1+\frac{1}{2\pi}(\frac{1}{\frac{1}{2}-\frac{0.4 \sqrt{2}}{18}}+
\frac{1}{\frac{3}{2}-0.4 \frac{\sqrt{2}}{18}}+\frac{1}{\frac{7}{2}-\frac{0.4 \sqrt{2}}{18}}
+\frac{1}{\frac{9}{2}-\frac{0.4 \sqrt{2}}{18}}) ) < 1.863.
\end{align}
\hfill $\Box$

\medskip

The lemma implies that $g_{(2),M}$ maps
the sets $E_{(2),M}^{\pm}{=}\{\xi \in \mathbb{C}: |\xi - \eta_M^{\pm}|< \frac{\sqrt{2}}{18}\}$ into $\mathbb{C} \setminus i \mathbb{Z}$.
Indeed, by the lemma
for $\xi \in E_{(2),M}^{\pm}$
$$
|(g_{(2),M})'(\xi)|\leq 1 + 1.863 \cdot \frac{\sqrt{2}}{18}< 1.15 .
$$
Hence, the image $g_{(2),M}(\xi)$ of a point $\xi \in E_{(2),M}^{\pm}$ has distance less than $\frac{1.15 \sqrt{2}}{18} <\frac{1}{2}$ from $\pm i(M-\frac{1}{2})=g_{(2),M}(\eta_M^{\pm})$.

The inequality \eqref{eq104a} and the equation \eqref{eq302e} imply the following upper bound for the
vertical side length $\mbox{vsl}(R_{(2),M})$  of the normalized rectangle $R_{(2),M}=r_{(2),M} R(M)$ with $M=\frac{d+1}{2}$ and $d\geq 5$

\begin{equation}\label{eq302f}
\mbox{vsl}(R_{(2),M}) \leq
(\frac{12}{5}+ \log(2M+1))\frac{\pi \sqrt{3}\sqrt{(M-\frac{1}{4})(M+\frac{1}{4})}}{{M+1}} < 1.504 \cdot\pi \sqrt{3} \log(4d-1).
\end{equation}
The last inequality follows
from the inequalities
$$
{\frac{12}{5} -\log 2.5} < 0.504 \log 19 < 0.504 \log(4d-1)
$$
for $d \geq 5$ and
$$
\log (\frac{5}{2}(d+2)) \leq    \log (4d-1)
$$
for $d\geq 5$.

\medskip

\noindent {\bf 3. Short syllables.}

\noindent  {\bf a) Syllables of form (1) of degree at most $4$ and of form (3)}.
Without loss of generality we assume that the syllable has the form $a_1^n$ where $n$ is a natural number at most equal to $4$. The other cases are treated symmetrically.
Put $M=\frac{d-1}{2}\leq \frac{3}{2}$ (with $d=n$). Consider the rectangle $R_{(1),M}$ in the plane with vertices $\pm i \frac{\pi}{2} (M+\frac{1}{2})$ on the imaginary line and vertices $(M+\frac{1}{2}) \log \frac{M}{M+1}\pm i \frac{\pi}{2} (M+\frac{1}{2})$ in the open left half-plane. Take $\xi_M^{\pm} =  -(M+\frac{1}{2}) \log \frac{M+\frac{1}{2}}{M}\pm i \frac{\pi}{2} (M+\frac{1}{2})$. Put
\begin{equation}\label{eq302'}
g_{(1),M}(\xi)= - M \exp(\frac{-\xi}{M+\frac{1}{2}}).
\end{equation}
The mapping $g_{(1),M}$ takes the rectangle $R_{(1),M}$ conformally onto a half-annulus in the left half-plane.
It maps the vertices $\pm \frac{\pi i}{2} (M+\frac{1}{2})$ of $R_{(1),M}$ contained the imaginary axis to $\pm i M$, and the vertices of $R_{(1),M}$ in the open left half-plane to $\pm i(M+1)$.
Hence, the mapping $g_{(1),M}:R_{(1),M} \to \mathbb{C} \setminus i \mathbb{Z}$ represents a lift of $a_1^n$.

The points $\xi_M^{\pm}$ are mapped to $\pm i(M+\frac{1}{2})$.
The derivative of $g_{(1),M}$ at the points  $\xi_M^{\pm}$ equals $\mp i$.
The second derivative of $g_{(1),M}$ can be estimated by
\begin{equation}\label{eq30}
|(g_{(1),M})''(\xi)|\leq  \frac{M+1}{(M+\frac{1}{2})^2} \leq \frac{3}{2}
\end{equation}
for $\xi \in R_{(1),M}$ and $M\geq\frac{1}{2}$.

By the same reasoning as above $g_{(1),M}$ maps the sets $E_{(1),M}^{\pm}\stackrel{def}{=}\{\xi \in \mathbb{C}: |\xi - \xi_M^{\pm}|< \frac{\sqrt{2}}{18}\}$ into $\mathbb{C} \setminus i\mathbb{Z}$.

For the considered values $M, \frac{1}{2}\leq M\leq \frac{3}{2}$, the horizontal side length of the rectangle $R_{(1),M}$ is at least $(M+\frac{1}{2})\log(1+\frac{1}{M})\geq \log \frac{5}{3}\geq 0.5> \frac{1}{18}$.
The vertical side length of the rectangle $R_{(1),M}$ equals
\begin{equation}\label{eq29'}
\mbox{vsl}(R_{(1),M})=\pi(M+\frac{1}{2})=\frac{\pi}{2} d.
\end{equation}
Since
\begin{align}\label{eq30'''}
  \frac{\pi}{2}   < 1.43 \cdot  \log 3, \;\; \frac{2\pi}{2}  < 1.615 \cdot  \log 7,\;\; \frac{3\pi}{2}  < 1.97 \cdot  \log 11, \;\;  \frac{4\pi}{2}< 2.33  \cdot \log 15,
\end{align}
we obtain the estimate
\begin{equation}\label{eq30'}
\mbox{vsl}(R_{(1),M}) < \frac{5}{2} \cdot \log(4d-1)
\end{equation}
for $d=1,2,3,4$ and $M=\frac{d-1}{2}$.

\medskip

\noindent {\bf  b) Short syllables of form (2) of degree $2$, $3$, and $4$.} Consider the syllables $a_2^{-1} a_1^{-1}$, $a_2^{-1} a_1^{-1} a_2^{-1}$, and $a_2^{-1} a_1^{-1} a_2^{-1} a_1^{-1}$. Other syllables of this form are treated similarly. Put $M=\frac{d}{2}$. Consider the rectangle $\tilde R_{(2),M}$ with vertices $0$ and $i d \frac{\pi}{2}$ on the imaginary axis and vertices $\log \frac{1}{2}$ and   $\log \frac{1}{2} +  i d \frac{\pi}{2}$ in the open left half-plane. $\tilde R_{(2),M}$ is contained in the rectangle $ R_{(2),M}$ with vertices $-i\frac{1}{18}$,
$i ( d \frac{\pi}{2} +\frac{1}{18}) $, $\log \frac{1}{2} - i\frac{1}{18}$, and   $\log \frac{1}{2} +  i (d \frac{\pi}{2} + \frac{1}{18})$. The horizontal side length of $ R_{(2),M}$ equals $\log 2 >\frac{1}{18}$.

Let 
$\eta^-_M=0$ and $\eta^+_M=i d\frac{\pi}{2}$. Put
\begin{equation}\label{eq30a'}
g_{(2),M}(\xi) = - \frac{i}{2} \exp(2\xi) , \;\xi \in R_{(2),M}.
\end{equation}
$g_{(2),M}$ maps the rectangle $R_{(2),M}$ into the annulus $\{\frac{1}{4}<|z|<\frac{1}{2}\}$.
Denote the extension of $g_{(2),M}$ to the closed rectangle by the same letter. We obtain   $g_{(2),M}(\eta^-_M)=-\frac{i}{2}$ and $ g_{(2),M}'(\eta^-_M)=-i$. Further, $g_{(2),M}(\eta^+_M)= -\frac{i}{2}$  if $d=2$ and $d=4$, and $+\frac{i}{2}$ if $d=3$, and       $(g_{(2),M})'(\eta^+_M)= -i$ if $d=2$ and $d=4$, and $+i$ if $d=3$.
The restriction of $g_{(2),M}$ to each maximal vertical segment in $\tilde R_{(2),M}$ makes $d=2M$ half-turns around $0$. Hence the mapping $g_{(2),M} : \tilde R_{(2),M} \to \mathbb{C} \setminus i \mathbb{Z}$  represents a lift of the syllable.

For the second derivative we have the estimate
\begin{equation}\label{eq30a}
|(g_{(2),M})''(\xi)|\leq 2
\end{equation}
on $R_{(2),M}$.

By the same reasoning as above $g_{(2),M}$ maps the sets $E_{(2),M}^{\pm}=\{\xi \in \mathbb{C}: |\xi - \xi^{\pm}|< \frac{\sqrt{2}}{18}\}$ into $\mathbb{C} \setminus i \mathbb{Z}$.

The vertical side length of the rectangle $R_{(2),M}$ equals
\begin{equation}\label{eq30*}
\frac{\pi}{2} d + \frac{1}{9}.
\end{equation}
Since $\pi + \frac{1}{9} < 1.672 \cdot \log 7$,
$\frac{3}{2} \pi + \frac{1}{9} < 2.0116  \cdot \log 11$, and
$2 \pi +\frac{1}{9} < 2.362  \cdot \log 15$
we obtain the estimate
\begin{equation}\label{eq30''}
\mbox{vsl}(R_{(2),M})< \frac{5}{2} \log(4d-1)
\end{equation}
for $d=2,3,4$ and $M=\frac{d}{2}$.

\medskip

\section{Quasiconformal gluing}
We give the {\bf proof of the upper bound in Theorem \ref{thm1} and Theorem $1'$.} Identify $w$ with $\mathfrak{C}_*(b) \in \pi_1$. Assume first that $w$ has at least two syllables and we are interested in $pb$ boundary values.
We associated to each syllable $\mathfrak{s}_j,\, j=1,\ldots,N,\,$ of the word $w$ a rectangle $\overset{\circ} {R_j}$ of the form $R_{(1),M}$ or $R_{(2),M}$ for a number $M=M_j$. If the rectangle is of the form
$R_{(1),M}$ we consider the holomorphic mapping  $g_{(1),M} + i M$ defined on it (see \eqref{eq301a} and \eqref{eq302'})and put $\overset{\circ}p_j^{\pm}=\xi_M^{\pm}$, otherwise we consider the mapping $g_{(2),M} + iM$ on $\overset{\circ}{ R_j}$ (see \eqref{eq302c}and \eqref{eq30a'}) and put $\overset{\circ}p_j^{\pm}=\eta_M^{\pm}$.
Denote the chosen  mapping on $ \overset{\circ}{R_j}$
by $\overset{\circ} {g_j}$. The mapping $\overset{\circ}g_j$ takes the rectangle $\overset{\circ}{R_j}$ to $\mathbb{C}\setminus i\mathbb{Z}$. In the first case the mapping represents a lift of $_{pb}(\mathfrak{s}_j)_{pb}$ under $f_1 \circ f_2$, in the second case the restriction of the mapping to the curvilinear rectangle $\tilde{R}_{(2),M}$
contained in $\overset{\circ}{R_j}$ has this property.

Choose complex numbers 
$c_j$ so that $\overset{\circ} p_j^+ + c_j= \overset{\circ} p_{j+1}^- + c_{j+1}$ for all except the last number $j$.
Put $p_j^{\pm} = \overset{\circ} p_j^{\pm} + c_j $.
Since for all $j$ the real part of $\overset{\circ} p_j^+$ coincides with the real part of $\overset{\circ} p_j^-$ we may choose the numbers $c_j$ so that the $p_j^{\pm}$ lie on the imaginary axis for all $j$. Denote by $R_j$ the translated rectangles
$\overset{\circ}{ R_j} + c_j$.
Consider the shifted mappings $g_j: R_j \to \mathbb{C} \setminus i\mathbb{Z}$, $g_j(z+c_j)= \overset{\circ}{g}_j(z) + i\, m_j$, where the integers $m_j$ are chosen so that $g_j(p_j^+)=g_{j+1}(p_{j+1}^-)$. The mappings $g_j$
are called normalized representatives of the lifts of the syllables $ \mathfrak{s}_j$. See Figure 6.

We describe now the choice of a closed curvilinear rectangle $\bar R$ contained in the union of the closed rectangles $ \overline {R_j}$, see Figure 6.
For each rectangle $R_j$ of the form $R_{(1),M}$ we consider the rectangle $R_j^-$ contained in
$R_j$ of the same vertical side length as $R_j$ and of horizontal side length $\frac{1}{18}$, for which $p_j^{\pm}$ are the midpoints of the horizontal sides.

For each rectangle $R_j$ of the form $R_{(2),M}$ we consider the rectangle $R_j^-$ contained in
$R_j$ of the same vertical side length as $R_j$ and of horizontal side length $\frac{1}{18}$,
one vertical side of which is contained in the imaginary axis. Notice that the points $p_j^{\pm}$ lie on this side.

\medskip
\begin{figure}[H]
\begin{center}
\includegraphics[width=9cm]{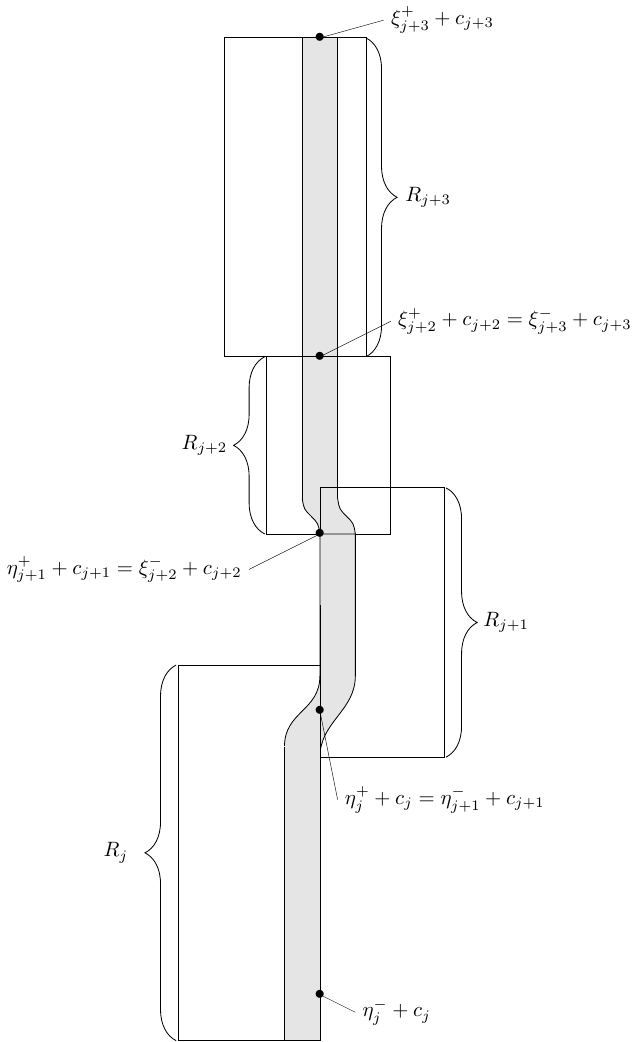}
\end{center}
\end{figure}
\centerline {Figure 6}

\bigskip

For $1<j<N$ we denote by  $R^0_j$ the subset of this rectangle that is contained in the horizontal strip $\{\mbox{Im}\,p_j^- +\frac{1}{18}< \mbox{Im}\, \xi <\mbox{Im}\, p_{j}^+ - \frac{1}{18}\}$.
Put $R_1^0= R_1^- \cap \{ \mbox{Im} \,\xi <\mbox{Im}\, p_{j}^+ - \frac{1}{18}\}$  and $R^0_N = R_N^- \cap \{\mbox{Im}\,p_j^- +\frac{1}{18}< \mbox{Im}\, \xi \}$.
The rectangles $R^0_j$ are contained in the $\frac{1}{18}$-neighbourhood of the imaginary axis. Let $q_j^{-}$ be the left endpoint of the lower side of $R^0_j$, and let $q_j^{+}$ be the left endpoint of the upper side of  $R^0_j$, respectively. For $j<N$ the inequality $|\mbox{Re} \, (q_j^+-q_{j+1}^-)|\leq \frac{1}{18}$ holds. Hence, for each $\varepsilon >0$  in each horizontal strip  $\{\mbox{Im}\,p_j^+ -\frac{1}{18}< \mbox{Im} \,\xi <\mbox{Im}\, p_j^+ + \frac{1}{18}\}$
around the $p_j^+=p_{j+1}^-$, $j=1,\ldots,N-1,$  there is a closed curvilinear rectangle whose horizontal sides are the upper side of $R^0_j$ and the lower side of $R^0_{j+1}$ which has the following two properties. Each horizontal line contained in $\{|\mbox{Im}\,(\xi-p_j^+)|\leq \frac{1}{18}\}$ intersects the curvilinear rectangle along a segment of length $\frac{1}{18}$. The vertical curvilinear sides are graphs of smooth functions over the segment $[i(\mbox{Im}\,p_j^+-\frac{1}{18}),i(\mbox{Im}\,p_j^++ \frac{1}{18})]$ of the imaginary axis whose derivatives vanish near the endpoints and have absolute value not exceeding $\frac{1}{2}+\varepsilon$. The number $\varepsilon >0$ will be chosen later as close to $0$ as needed. The closed curvilinear rectangle $\bar R$ is the union of the $\overline{R^0_j}$ and all obtained closed curvilinear rectangles. See Figure 6.

We will now patch together the functions $g_j$,  $  g_ j(\xi + c_j) = \overset{\circ} {g_j}(\xi) + i\, m_j$, on $R_j $, using their holomorphic extensions to a $\frac{\sqrt{2}}{18}$-neighbourhood of $p_j^{\pm}$. Recall that the integers $m_j$ are chosen in such a way that for  all $j$ (except the last one) the values of $g_j $ and of $g_{j+1}$ at $p_j^+=p_{j+1}^-$ coincide and are equal to an imaginary half-integer which is not an imaginary integer. Moreover, by the normalization of the ${R}_j$ the values of their derivatives at the point $p_j^+$
coincide. Both functions,  $ g_j $ and $g_{j+1}$ map the $\frac{\sqrt{2}}{18}$-neighbourhood of the point into a $\frac{1}{2}$-neighbourhood of an imaginary half-integer which is not an imaginary integer.

Consider the $C^1$-function $\chi_0$ on the interval $[0,1]$, $\chi_0(t)=6 \int_0^t \tau (1-\tau) d\tau$. Then $\chi_0(0)=0$, $\chi_0(1)=1$, and $0\leq \chi_0'(t)\leq 6t(1-t)\leq \frac{3}{2}$. 
Put $\chi(t)=\chi_0(9 t)$.

Define a function $g$ on $R$ as follows. 
Each point $\xi$ in $R$ for which $|\mbox{Im}\,(\xi -p_j^+)|> \frac{1}{18}$ for all $j<N$
belongs to a single rectangle $\overline{R_k }$ (depending on $\xi$) and we put $g(\xi)\stackrel{def}{=}g_k(\xi)$ for such a point.

Fix a number $j<N$ 
and consider the set $Q_j\stackrel{def}{=}\{\xi:|\mbox{Im}\,(\xi -p_j^+)|\leq \frac{1}{18}\}$. Put  
$\chi_j(\xi)=\chi(\mbox{Im}\,(\xi -p_j^+)+\frac{1}{18})\, $ for $\xi \in  Q_j$.
Let $g\stackrel{def}{=} (1-\chi_j)  g_j + \chi_j  g_{j+1}$ on $ \bar R \cap Q_j$.
For $\{\xi \in Q_j: \mbox{Im}\,\xi = p_j^+ -\frac{1}{18}\}$ the equalities   
$\chi_j(\xi )= \chi_0(0)=0$ and $\chi'_j(\xi)=\chi_0'(0)=0$ hold. Hence, the function $g$ is $C^1$ smooth near such points $\xi$. Further, for $\{\xi \in Q_j:\mbox{Im}\,\xi = p_j^+ +\frac{1}{18}\}$ the equalities   $\chi(\mbox{Im}\,(\xi -p_j^+)+\frac{1}{18})= \chi_0(1)=1$ and $\chi_j'(\xi)=\chi_0'(1)=0$ hold, hence, the function $g$ is $C^1$ smooth near such $\xi$. Since both functions $ g_j$ and $ g_{j+1}$ map the $\frac{\sqrt{2}}{18}$-neighbourhood of $p_j^+$ into a $\frac{1}{2}$-neighbourhood of an imaginary half-integer which is not an imaginary integer, the convex combination $g$ of these two functions has the same property. Make the same definition for all but the last number $j$. We obtain a smooth mapping $g$ from the curvilinear rectangle $\bar R$ to $\mathbb{C} \setminus i\mathbb{Z}$ which represents a lift of $_{pb}(w)_{pb}$ under $f_1 \circ f_2$.

\begin{lemm}\label{lemm221} The mapping $g$ is a quasiconformal mapping from $R$ onto its image. The Beltramy differential $\mu_g$ of $g$ has absolute value $|\mu_g|< 0.1712$.
\end{lemm}

\medskip
\noindent {\bf Proof}. Put $\xi=u+i v$. If $\xi \in R$, $|\mbox{Im}\, (\xi -p_j^+)|\leq \frac{1}{18}$ for some $j, \, 1 \leq j <N,\,$  then the Beltrami differential at $\xi$ equals
$$
\mu_g(\xi)= \frac{\frac{\partial}{\partial \overline \xi}g(\xi)}{\frac{\partial}{\partial \xi}g(\xi)} = \frac{\frac{i}{2} (\frac{\partial}{\partial v}{\chi_j}\cdot (g_{j+1} -g_{j}))(\xi)}{(\frac{-i}{2} \frac{\partial}{\partial v}{\chi_j}\cdot (g_{j+1} -g_{j})+ (1-\chi_j)\cdot g_j' + \chi_j\cdot g_{j+1}')(\xi)}.
$$
On the rest of the rectangle $R$ the function is analytic.

By the Lemmas \ref{lemm219} and \ref{lemm220} and inequalities  \eqref{eq30}
and  \eqref{eq30a}
the estimate
$$
\max\{|g_j''(\xi)|,|g_{j+1}''(\xi)|\}< 2.75
$$
holds for the considered points $\xi$. Since $g_{j+1}-g_{j}$ vanishes together with its derivative at $p_j^+$, the estimate
$$
|(g_{j+1}-g_{j})(\xi)|\leq 2 \cdot \frac{1}{2} \cdot (\frac{\sqrt{2}}{18})^2 \cdot 2.75
$$
holds for $\xi$ in the $\frac{\sqrt{2}}{18}$-neighbourhood of $p_j^+$. Further
$$
|\frac{\partial}{\partial v}{\chi_j}|\leq\frac{3}{2} \cdot 9,
$$
on $Q_{j}$, and
$$
\max\{|g_j'-g_j'(p_j^+)|,|g_{j+1}'-g_j'(p_j^+)|\}<2.75 \cdot \frac{\sqrt{2}}{18}
$$
on the  $\frac{\sqrt{2}}{18}$-neighbourhood of $p_j^+$ by the estimate for the second derivative of the $g_k$, since $g_j'(p_j^+)=g_{j+1}'(p_j^+)$.
Hence,
$$
k = \sup_R |\mu_g(\xi)| < \frac{\frac{1}{24} \cdot 2.75}{1-\frac{1}{24} \cdot 2.75-\frac{\sqrt{2}}{18}\cdot 2.75} < 0.1712  .
$$
We used that $|g_j'(p_j^+)|=1$.
\hfill $\Box$

\medskip

The quasiconformal dilatation $K=\frac{1+k}{1-k}$ does not exceed $1.414$.

Let $\omega$ be the normalized solution of the Beltrami equation
$$
\frac{\partial}{\partial \overline z} \omega = \tilde \mu_g \frac{\partial}{\partial z} \omega
$$
on the complex plane. Here $\tilde \mu_g$ equals $\mu_g$ on $\bar R$ and equals $0$ outside $ \bar R$. $\omega$ is a H\"older continuous self-homeomorphism of the complex plane. The mapping $g \circ  {\omega}^{-1}$ is holomorphic on $\omega(R)$ (see \cite{A1}, Chapter I C).
The image $\omega(R)$  can be considered as a curvilinear rectangle. The curvilinear sides are the images of the sides of $R$. By \cite{A1} (chapter I, Theorem 3) the extremal length of $\omega(R)$  does not exceed $K \cdot \lambda(R)$. In other words, there is a conformal mapping $\psi$ of a true rectangle $\mathcal{R}$ of extremal length not exceeding  $K \cdot \lambda(R)$ onto $\omega(R)$, which takes the sides of $\mathcal{R}$ to the respective curvilinear sides of $\omega(R)$ . The mapping $g\circ {\omega}^{-1} \circ \psi: \mathcal{R} \to \mathbb{C} \setminus i \mathbb{Z}$ is a holomorphic mapping from the rectangle $\mathcal{R}$ of extremal length not exceeding  $K \cdot \lambda(R)$ to $\mathbb{C}\setminus i\mathbb{Z}$ that represents a lift of $_{pb}(w)_{pb}$.

Estimate the extremal length $\lambda(R)$ of $R$. By Lemma \ref{lemm216}
$$
\lambda(R) \leq (1+C^2)\; \sum \frac{\mbox{vsl}( R_j)}{\textsf{b}}.
$$
Here $\textsf{b}=\frac{1}{18}$ and the rectangle $\bar R$ was chosen so that $C$ does not exceed $\frac{1}{2} + \varepsilon$.

We obtain
\begin{align}\label{eq31a}
\lambda(\mathcal{R}) &\leq  1.414 \cdot (1+ (\frac{1}{2} + \varepsilon)^2)    \cdot 18
\cdot \sum \mbox{vsl}(R_j).
\end{align}
Further, by \eqref{eq302b},  \eqref{eq302f}, \eqref{eq30'}, and \eqref{eq30''}
\begin{align}\label{eq31b}
\mbox{vsl}(R_j) \leq 1.504 \cdot \pi \sqrt{3} \log(4d_j-1),
\end{align}
where $d_j$ is the degree of the $j$-th syllable. Hence,
\begin{align}\label{eq31}
\lambda(\mathcal{R}) &\leq  1.414 \cdot (1+ (\frac{1}{2} + \varepsilon)^2)    \cdot 18 \cdot 1.504 \cdot \sqrt{3} \pi \cdot \sum \log(4d_j-1)
\end{align}
Since $1.414 \cdot (1+ \frac{1}{4})    \cdot 18 \cdot 1.504 \cdot \sqrt{3} \pi < 260.4 $,
the number $\varepsilon$ may be chosen so that the right hand side of \eqref{eq31} does not exceed $300 \cdot \sum \log(4d_j-1)$, i.e.
\begin{align}\label{eq31c}
\lambda(\mathcal{R}) \leq  300 \cdot \sum \log(4d_j-1).
\end{align}

The upper bound is proved for $pb$ boundary values.

\medskip

Suppose now, for instance, that the left boundary values are $tr$. We assume again that $w$ has at least two syllables.
Suppose the first syllable $\mathfrak{s}_1$ is of form (1) or (3), say, without loss of generality, $\mathfrak{s}_1=a_1^n$ for $n\geq 1$. A lift of $_{tr}(\mathfrak{s}_1)_{pb}$  is a half-slalom class with $M=d-\frac{1}{2}$ (see Proposition \ref{prop4b}).

Consider the rectangle $R_{(1),2M}$ and the mapping $g_{(1),2M}$ on this rectangle. The restriction of $g_{(1),2M}$ to the ''upper half '' $R_{(1),2M}^{\#}$ of $R_{(1),2M}$ represents a lift of $_{tr}(\mathfrak{s}_1)_{pb}$.
Shift the rectangle $R_{(1),2M}^{\#}$ and the restriction of $g_{(1),2M}$ to it
in the needed way and denote the obtained rectangle by $R_1$ and the function by $g_1$.

If $n=d \geq 3$ then 
by the first inequality in \eqref{eq302b} the vertical side length of
$R_{(1),2M}$ does not exceed $\frac{1}{2}(\frac{12}{5}+ \log (4M+1))=
\frac{1}{2}(\frac{12}{5}+ \log (4d-1))$.
For $d\geq 3$ we have $\log(4d-1) \geq \log 11 > 2.3$. Since $1.1 \cdot 2.3 >2.5 >\frac{12}{5}$  we obtain
\begin{equation}\label{eq31d}
\mbox{vsl}({R_j})\leq \frac{1}{2} ( 2.1 \cdot \log(4d-1))  = 1.05 \cdot \log(4d-1).
\end{equation}

If $n=d\leq 2$ then by \eqref{eq29'} the vertical side length of $R_{(1),M}$ does not exceed $\pi (M+\frac{1}{2})=\pi d$. Hence by \eqref{eq30'''}
$$\mbox{vsl}({R_j}) \leq  \frac{1}{2}\pi d <2 \log(4d-1)$$
for $d=1$ and $2$.

Suppose the first syllable $_{tr}(\mathfrak{s}_1)_{pb}$ is of form (2). 
We may assume that it has the form $(a_2^{-1} a_1^{-1})^k$
or $(a_2^{-1} a_1^{-1})^k a_2^{-1}$ for $k\geq 1$.

Consider first the case when $d\geq 5$. With $d=2M-1$ the curve $y \to \frac{1}{2}\exp(\pi(iy + i(M-1))), y \in [-(M-\frac{1}{2}), (M-\frac{1}{2}) ],$ represents a lift of the syllable $_{tr}(\mathfrak{s}_1)_{pb}$ to $\mathbb{C} \setminus i\mathbb{\mathbb{Z}}$.
We will extend the curve to a mapping of a curvilinear rectangle contained in $\tilde {R}_{(2),M}$ (see \eqref{eq302c}) which represents $_{tr}(\mathfrak{s}_1)_{pb}$.

For this purpose we notice first that for any curve $\alpha$ with interior in $\{z \in \mathbb{C}_{\ell}: -(M-\frac{1}{2}) < \mbox{Im} \, z < ( M-\frac{1}{2})\}$, with initial point in 
$[\frac{1}{\pi}\log \frac{1}{2},0]-i(M-\frac{1}{2})$ and with endpoint in $(-\infty,0) + i(M-\frac{1}{2})$,
the curve $\exp (\pi \alpha +i(M-1))$ intersects $-\frac{i}{2} + \mathbb{R}$. (See Figure 7.)
For $x \in [\frac{1}{\pi}\log \frac{1}{2},0]$ we let $y(x)$ be the first parameter of intersection of the curve $y \to \exp(\pi(x+iy + i(M-1))), y \in [-(M-\frac{1}{2}),(M-\frac{1}{2}) ],$ with $-\frac{i}{2}+\mathbb{R}$. We obtain an arc $\beta$, $x \to \beta(x)=(x,y(x)), \, x \in  [\frac{1}{\pi}\log\frac{1}{2},0],$ with left endpoint equal to $\frac{1}{\pi}\log\frac{1}{2}-i(M-\frac{1}{2})$ and right endpoint on $(-i(M-\frac{1}{2}),i(M-\frac{1}{2}))$. The interior of the arc is contained in $\mathbb{C}_{\ell} \cap \{|\mbox{Im}\, \zeta| < M-\frac{1}{2}\}$.
The image of $\beta$ under $f_{(2),M}$  is an arc with interior in
$f_{(2),M}(\{z \in \mathbb{C}_{\ell}: -(M-\frac{1}{2}) < \mbox{Im} \, z < ( M-\frac{1}{2})\})$,
with initial point equal to $f_{(2),M}(\frac{1}{\pi}\log\frac{1}{2}-i(M-\frac{1}{2}))$ and right endpoint on $(-\eta_M^-,\eta_M^+)$. (See Figure 7.)

The initial point of $f_{(2),M} \circ \beta$ is not contained in the closure $\overline{ \tilde R_{(2),M}}$ of $\tilde R_{(2),M}$. Indeed, the point lies on the curve $f_{(2),M}((-\infty,0) -i(M-\frac{1}{2}))$. If it was contained in the closure
$\overline{ \tilde R_{(2),M}}$ of $\tilde R_{(2),M}$ then by Lemma \ref{lemm19} it would have distance not exceeding $\frac{1}{18 \cos(\arctan(0.05))}$ from $\eta_M^-$ and by Corollary  \ref{cor204} the distance of $\frac{1}{\pi}\log\frac{1}{2}-i(M-\frac{1}{2})$
from $i(M-\frac{1}{2})$ would not exceed $ \frac{1}{18 \cos(\arctan(0.05))} \cdot 0.3343  $.
But the latter distance equals $\frac{1}{\pi}\log 2 > \frac{{1}}{18 \cos(\arctan(0.05))} \cdot 0.3343$.

As a consequence, the arc $f_{(2),M} \circ \beta$ intersects $\tilde R_{(2),M}$ along a union of relatively closed arcs. There is exactly one connected component $\tilde{\beta}$ of the intersection that has an endpoint on  $(-\eta_M^-,\eta_M^+)$. The other endpoint of the component $\tilde{\beta}$ lies on the left side of $\tilde R_{(2),M}$. The arc $\tilde{\beta}$  divides  $\tilde R_{(2),M}$ into two connected components. Denote by  $\tilde R_{(2),M}^0$ the component of  $\tilde R_{(2),M} \setminus \tilde{\beta}$ whose closure contains $\eta_M^+$.
The restriction of $g_{(2),M}$ to $\tilde R_{(2),M}^0$ represents the syllable $_{tr}(\mathfrak{s}_j)_{pb}$. The vertical side length of $\tilde R_{(2),M}^0$ does not exceed the vertical side length of $ R_{(2),M}$.

\medskip

\begin{figure}[H]
\begin{center}
\includegraphics[width=15cm]{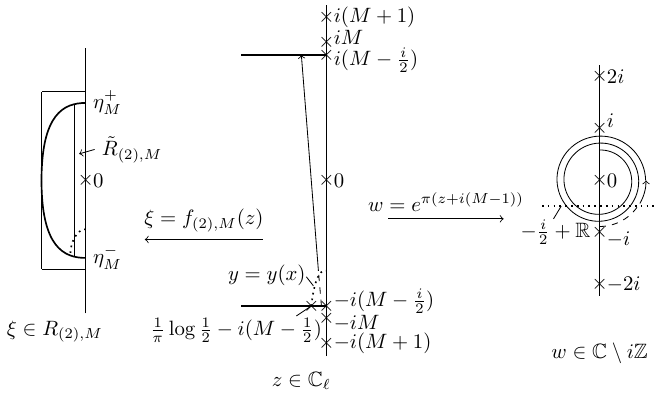}
\end{center}
\end{figure}

\centerline {Figure 7}

\bigskip

The case when the degree of the syllable is either $2$, $3$, or $4$ can be treated
similarly but is simpler because the mapping $g_{(2),M}$ is simpler. Also in this case
we obtain a rectangle $\tilde R_{(2),M}^0 \subset \tilde R_{(2),M}$ such that the restriction $g_{(2),M}|\tilde R_{(2),M}^0$ represents the syllable $_{tr}(\mathfrak{s}_j)_{pb}$.
The vertical side length of $\tilde R_{(2),M}^0$ does not exceed the vertical side length of $ R_{(2),M}$.

The case when $w$ has totally real right boundary values or both boundary values are totally real, is treated in the same way. The quasiconformal gluing is done as in the case of $pb$ boundary values.

If $w$ consists of a single syllable the upper bound in the case of mixed boundary values or in the non-exceptional cases of $tr$ or $pb$ boundary values follows directly from Proposition \ref{prop4b}. The exceptional cases were treated in the proof of the lower bound in Theorem \ref{thm1}.

Theorem 1 and Theorem $1'$ are proved.  \hfill $\Box$

\medskip
\noindent {\bf Proof of the upper bound in Theorem 2}. 
Take any syllable reduced word $w$ representing the conjugacy class $\mathfrak{C}_*(\hat b)$. Consider all syllables $\mathfrak{s}_j , \, j=1,\ldots, N,$ of $w$, labeled from left to right, with $pb$ boundary values. Consider for each $j$ the rectangle $\overset{\circ}{R_j}$ and the holomorphic mapping $\overset{\circ}{g_j}$ on it.
Extend the finite sequence $\mathfrak{s}_j$ of syllables to an infinite sequence of syllables $\mathfrak{s}_j , \, j \in \mathbb{Z},$
with $\mathfrak{s}_j = \mathfrak{s}_{j+N}$ for all $j$. Put also $\overset{\circ}{R_j}= \overset{\;\circ}{R}_{j+N}$ and $\overset{\circ}{g_j}= \overset{\;\circ}{g}_{j+N}$ for integers $j$ and $N$.  In the same way as in the case of finite sequences of syllables we consider for each $j$ the normalized representative $g_j: R_j \to \mathbb{C} \setminus i\mathbb{Z}$  of the lift under $f_1\circ f_2$ of the syllable $ _{pb}(\mathfrak{s}_j)_{pb}$.
The normalization is done in such a way that the value and the first derivative of $g_j$ and $g_{j+1}$ coincide at $p^+_j=p_{j+1}^-$. Recall that all $p_j^{\pm}$ are on the imaginary axis. Hence, $R_{N+j}=R_j + i\textsf{a}$ for all $j$
and a positive number $\textsf{a}$,
and $g_{j+N}(\xi+ i  \textsf{a} )   =g_j(\xi) + i\, m$ for all $j \in \mathbb{Z}$ and an integer $m$.

Do quasiconformal gluing for all $j\in\mathbb{Z} $ by the same procedure as described before. We may perform quasiconformal gluing in such a way that for the obtained quasiconformal mapping $g$ on the infinite curvilinear strip $R_{\infty}$ of width $\frac{1}{18}$ the following relation
$$
g(\xi+i \textsf{a})=  g(\xi)+ m\, i
$$
holds.  The quasiconformal dilatation $K$ of the mapping $g$ on the infinite strip has the same estimate $K\leq 1.414$ as the quasiconformal dilatation of the mapping constructed in the proof of Theorem \ref{thm1}.
Since $f_1 \circ f_2$ has period $i$ the composition $G=f_1 \circ f_2 \circ g$ has period $i\, \textsf{a}$.

Let the annulus $A$ be the quotient of  $R_{\infty}$ by the equivalence relation $\xi \sim \xi + i \,\textsf{a}$. The function $G$ descends to a function
$\hat G: A \to \mathbb{C } \setminus \{-1,1\}$  which has quasiconformal dilatation 
$K$ and represents $\hat b$. The annulus $A$ is conformally equivalent to the annulus $A' \subset \mathbb{C}$, which is the image of  $R_{\infty}$ under the mapping $\xi \to \exp( \frac{2\pi \xi}{ \textsf{a}}) $. Consider the infinite strip
of width $\frac{1}{18}$, that is bounded by two vertical lines.
Let  $A_0$ be the quotient of this strip by the equivalence relation $\xi \sim \xi + i\, \textsf{a} $. The inequality for the extremal length $\lambda(A') \leq (1+C^2)\lambda(A_0)$ is obtained by the same arguments that are used for the proof of
Lemma \ref{lemm216}.

Putting together the estimates for the quasiconformal dilatation of $\hat G$, for the constant $C$, and for the estimate of the extremal length of $A_0$ through the vertical side length's of the $R_j$,  we obtain the same upper bound $300 \cdot \sum \log(4d_j-1)$ as in Theorem 1. Theorem 2 is proved. \hfill $\Box$

\medskip

\noindent {\bf Proof of the upper bound of Theorem 3}.
We prove the upper bound for $3$-braids that are not pure and not among the exceptional cases. Write such a braid in the form \eqref{eq2'}. We may assume that the braid equals $\sigma_j^k b_1$ with $j=1$ or $j=2$, an odd integral number $k$ and a pure braid $b_1$ as in Lemma \ref{lemm1'} which is not the identity. Indeed, multiplying by a power of $\Delta_3$ does not change the extremal length. Suppose $j=1$. Then each curve $\gamma$ in $C_3(\mathbb{C})\diagup \mathcal{S}_3$ representing $b_{tr}$ can be decomposed into two curves, $\gamma_0$ representing $_{tr}(\sigma_1^k)_{pb}$, and $\gamma_1$ representing $_{pb}(b_1)_{tr}$. This can be seen by lifting $\gamma$ to a curve $\tilde \gamma$ in $C_3(\mathbb{C})$ with initial point
on $\{x_2<x_1<x_3\}$, and applying the analog of Lemma \ref{lemm14a} to $\mathfrak{C}(\tilde \gamma)$. Note that the terminating point of $\tilde \gamma$
lies on $C_3(\mathbb{C})^0=\{x_1<x_2<x_3\}$. Denote the lift of $\gamma_0$ by $\tilde{\gamma}_0$.

The arc $\mathfrak{C}(\tilde \gamma_0)$  is an arc in $\mathbb{C} \setminus \{-1,1\}$
with initial point in $(-\infty,-1)$ and terminating point in the imaginary axis. Recall that for each integer $k'$ the mapping $f_1 \circ f_2$ takes the interval $(-\infty,0)+i k'$
onto $(-\infty,-1)$, and for each  $j' \in\mathbb{Z}$ the mapping $f_1 \circ f_2$ takes the interval $(i j', i( j' +1))$ onto the imaginary axis. (See also Figure 1.)
More precisely,
put $k=2\ell+1$ for an integer $\ell$. Assume $\ell$ is non-negative.
(The case of negative integers $\ell$ is treated similarly.)
Take for $\Gamma_0$ the lift of $\mathfrak{C}(\tilde \gamma_0)$ to $\mathbb{C} \setminus i\mathbb{Z}$ with initial point in $(-\infty,0)$. Then the terminating point of the lift is contained in $(i \ell, i( \ell +1))$.(See also Figure 1.)

Associate to $\mathfrak{C}(_{tr}(\sigma_1^k)_{pb})$ a rectangle $\overset{\circ}{R_0}$ and a holomorphic function $\overset{\circ}{g_0}: \overset{\circ}{R_0}\to \mathbb{C} \setminus i\mathbb{Z}$ that represents the   class  of curves that are homotopic in $\mathbb{C} \setminus i\mathbb{Z}$ to $\Gamma_0$ with initial point in $(-\infty,0)$ and terminating point in  $(i \ell, i( \ell +1))$.
For $\ell \geq 1$  the class of $\Gamma_0$ can be represented by the restriction  of the mapping $g_{(1),M}$ to the upper half of the normalized  rectangle  $R_{(1),M}$ with $M=\ell=\frac{k-1}{2} $. If $\ell \geq 3$ the vertical side length
of this "half-rectangle" does not exceed
\begin{align}\label{eq33}
\frac{1}{2}(\frac{12}{5}+ \log(2 \ell +1))< \frac{1}{2}(\log(32 \ell +16)) < \log(4\ell -1).
\end{align}
(We used the first inequality in \eqref{eq302b} and the relation $\frac{12}{5}<2.7     <\log 16$.)

If $\ell=1$ or $2$ the
vertical side length of this "half-rectangle" does not exceed
\begin{align}\label{eq34a}
\frac{\pi}{2} (M+\frac{1}{2})= \frac{\pi}{2} (\ell+\frac{1}{2}) < 3 \log(4\ell -1).
\end{align}
(We used \eqref{eq29'} and the fact that $\frac{\pi}{2} \frac{3}{2} < 3 \log 3$ and $\frac{\pi}{2} \frac{5}{2} < 3 \log 7$, see also \eqref{eq30'''}.)

Suppose $k=1$ (i.e. $\ell=0$). In this case the class of $\Gamma_0$ can be represented by a conformal mapping of a rectangle to a quarter of an annulus. Indeed, consider the normalized rectangle with vertices $\frac{\pi i}{4}$ and $0$ on the real axis and vertices
$\frac{\pi i}{4}-\log 2$ and $-\log 2$ in the open left half-plane. The mapping $\xi\to
-\frac{1}{4} e^{-2\xi}$ takes $\frac{\pi i}{4}$ to $\frac{i}{4}$, $0$ to $-\frac{1}{4}$,  $\frac{\pi i}{4}-\log 2$ to $i$ and $-\log 2$ to $-1$.
The image of the rectangle under this mapping is the upper left quarter of the annulus $\{\frac{1}{4}<|z|<1\}$.

The mapping represents the class of $\Gamma_0$. Put $\xi_0^+= -\frac{1}{2}\log 2 +   \frac{\pi i}{4}$. Then the derivative at this point equals $-i$.
We obtained a normalized rectangle $\overset{\circ}{R}_0$ and a mapping $\overset{\circ}{g}_0$. The normalized rectangle
has vertical side length $\frac{\pi}{4}$.

Let us prove the theorem first in the case when $k=1$ or $k=-1$. Then the class
$_{tr}b_{tr}$  is the product of  $_{tr}(\sigma_1^{\pm 1})_{pb}$ and  $_{pb} (b_1)_{tr}$. The class  $_{tr}(\sigma_1^{\pm 1})_{pb}$ can be represented  by a holomorphic mapping on a normalized rectangle of vertical side length not exceeding $\frac{\pi}{4}$. The class $_{pb} (b_1)_{tr}$  is a pure braid which can be decomposed into parts that can be treated as in the proof of the upper bound of Theorem $1'$. 
Quasiconformal gluing as in Theorem $1'$ gives
the estimate

\begin{align}\label{eq35}
\Lambda(b_{tr}) \leq
1.414 \cdot (1+ (\frac{1}{2} + \varepsilon)^2)    \cdot 18 \cdot (1.504 \cdot \sqrt{3} \pi \cdot \sum_{\mathfrak{s}_j}  \log(4d(\mathfrak{s}_j)-1) + \frac{1}{4} \pi),
\end{align}
where the ${\mathfrak{s}_j}$ run over the syllables of the image of $ \vartheta(b)$ in the braid group modulo its center. Since $\mathfrak{C}_*(\vartheta(b))$ is not the identity the expression $\mathcal{L}(\mathfrak{C}_*(\vartheta(b))) = \sum_{\mathfrak{s}_j}  \log(4d(\mathfrak{s}_j)-1)$ is not smaller than $\log 3$.
Further, the inequality
\begin{align}\label{eq36}
\frac{\pi}{4}\leq 0.715 \cdot \log 3 \leq  0.715 \cdot \mathcal{L}(\mathfrak{C}_*(\vartheta(b))),
\end{align}
holds.
Since $1.414 \cdot \frac{5}{4}    \cdot 18 \cdot ( 1.504 \cdot \sqrt{3} \pi + 0.715)<  283.12 < 300$, the number $\epsilon$ can be chosen small so that
\begin{align}\label{eq37}
\Lambda(b_{tr}) & \leq 1.414 \cdot (1+ (\frac{1}{2} + \varepsilon)^2)   ( \cdot 18 \cdot 1.504 \cdot \sqrt{3} \pi + 0.715) \cdot \mathcal{L}(\mathfrak{C}_*(\vartheta(b))) \nonumber \\
& < 300 \cdot \mathcal{L}(\mathfrak{C}_*(\vartheta(b))).
\end{align}


Consider the case when $|k|\geq 3$ is an odd integer and  $\mathfrak{C}_*(\sigma_1^{q(k)})$ is a syllable of $\mathfrak{C}_*(\vartheta(b))$.  By our estimates  \eqref{eq33} and \eqref{eq34a} for the vertical side length of the normalized rectangle $\overset{\circ}{R_0}$ associated to $\mathfrak{C}_*(_{tr}(\sigma_1^{q(k)})_{pb})$,  quasiconformal gluing gives us the same estimate for $\Lambda(b_{tr})$  as Theorem \ref{thm1} gives for $\Lambda(_{tr}(\vartheta(b))_{tr})$. Hence we obtain again
\begin{equation}\label{eq37'}
\Lambda(b_{tr}) < 300 \cdot \mathcal{L}(\mathfrak{C}_*(\vartheta(b))).
\end{equation}

It remains to consider the case when the term $\mathfrak{C}_*(\sigma_1^{q(k)})$ is not the identity and not a syllable of $\mathfrak{C}_*(\vartheta(b))$. In this case the exponent $k$ equals $\pm 3$.
The vertical length of the normalized rectangle corresponding to $\mathfrak{C}_*(_{tr}(\sigma_1^{\pm 3})_{pb})$ does not exceed
$3 \log 3$ (see \eqref{eq34a}). Since the term $\mathfrak{C}_*(\sigma_1^{q(k)})= \mathfrak{C}_*(\sigma_1^{\pm 2}) $ is not a syllable of $\mathfrak{C}_*(\vartheta(b))$, the first  syllable $\mathfrak{s}_1$ of $\mathfrak{C}_*(\vartheta(b))$
has the form $\mathfrak{C}_*(\sigma_1^2 \sigma_2^2 \dots)$ or $\mathfrak{C}_*(\sigma_1^{-2} \sigma_2^{-2} \dots)$
and has degree at least two. Hence $\mathcal{L}(\mathfrak{C}_*(\vartheta(b)))\geq \log 7$.

Suppose first the syllable $\mathfrak{s}_1$ equals $\mathfrak{s}_1= \mathfrak{C}_*(\sigma_1^2 \sigma_2^2) $ or $\mathfrak{C}_*(\sigma_1^{-2} \sigma_2^{-2})$ .
The normalized rectangle corresponding to $_{pb}(\sigma_2^{\pm 2})_{pb}$ or $_{pb}(\sigma_2^{\pm 2})_{tr}$  has vertical side length not exceeding $3 \log 3$ (see \eqref{eq30'}). Hence, as in Theorem 1 and Theorem $1'$ quasiconformal gluing  of the mappings representing $\mathfrak{C}_*(_{tr}(\sigma_1^{\pm 3})_{pb})$, $\mathfrak{C}_*(_{pb}(\sigma_2^{\pm 2})_{\#})$, and  syllables $\mathfrak{s}_j,\, j \geq 2,$  of the word $\mathfrak{C}_*(\vartheta(b))$ (if there are such)  gives the estimate
\begin{align}\label{eq38}
\Lambda(b_{tr}) &\leq 1.414 \cdot (1+ (\frac{1}{2} + \varepsilon)^2)    \cdot 18 \cdot ( 6 \cdot \log 3 + 1.504 \cdot \sqrt{3} \pi \cdot \sum_{j\geq 2} \log(4d(\mathfrak{s}_j)-1)),
\end{align}
where  ${\mathfrak{s}_j}$ runs over the syllables of $\mathfrak{C}_*(b_1) $.
Since
\begin{align}\label{eq39}
 6\log 3 \leq  0.414  \cdot 1.504 \cdot \sqrt{3}\cdot \pi\cdot  \log 7 < 1.504 \cdot \sqrt{3} \pi \cdot \log(4 d(\mathfrak{s}_1) -1),
\end{align}
we obtain the upper bound \eqref{eq37'} of Theorem \ref{thm3} for this case.

For the case when the degree of $\mathfrak{s}_1$ equals $3$ the same arguments give
the estimate
\begin{align}\label{eq40}
\Lambda(b_{tr}) &\leq 1.414 \cdot (1+ (\frac{1}{2} + \varepsilon)^2)    \cdot 18 \cdot ( 3 \log 3 + \pi + \frac{1}{9}      + 1.504 \cdot \sqrt{3} \pi \cdot \sum_{j\geq 2} \log(4d(\mathfrak{s}_j)-1)).
\end{align}
In case the degree of $\mathfrak{s}_1$ equals $4$ the estimate is
\begin{align}\label{eq41}
\Lambda(b_{tr}) &\leq 1.414 \cdot (1+ (\frac{1}{2} + \varepsilon)^2)    \cdot 18 \cdot ( 3 \log 3 + \frac{3}{2}\pi + \frac{1}{9}      + 1.504 \cdot \sqrt{3} \pi \cdot \sum_{j\geq 2} \log(4d(\mathfrak{s}_j)-1)).
\end{align}
(See \eqref{eq30*} for $d=2$ and $d=3$ and the argument in the proof of Theorem \ref{thm1} for other than $pb$ boundary values.)
Since
\begin{align}\label{eq42}
 3\log 3 + \pi + \frac{1}{9} \leq  0.334  \cdot 1.504 \cdot \sqrt{3}\cdot \pi\cdot  \log 11 < 1.504 \cdot \sqrt{3} \pi \cdot \log(4 d(\mathfrak{s}_1) -1),
\end{align}
in the first case, and
\begin{align}\label{eq43}
 3\log 3 + \frac{3}{2}\pi + \frac{1}{9} \leq  0.3664  \cdot 1.504 \cdot \sqrt{3}\cdot \pi\cdot  \log 15 < 1.504 \cdot \sqrt{3} \pi \cdot \log(4 d(\mathfrak{s}_1) -1),
\end{align}
in the second case, we obtain the upper bound in Theorem \ref{thm3} also in these cases.

In the case when the degree of $\mathfrak{s}_1$ is at least $d_1=5$, we use that the vertical length of the normalized rectangle corresponding to the expression $\mathfrak{C}_*(_{pb}(\sigma_2^{\pm 2} \sigma_1^{\pm 2} \ldots)_{\#})$ (of degree $d_1-1$) does not exceed the vertical length of the rectangle corresponding to $\mathfrak{s}_1$. We obtain
\begin{align}\label{eq38}
\Lambda(b_{tr}) &\leq 1.414 \cdot (1+ (\frac{1}{2} + \varepsilon)^2)    \cdot 18 \cdot (  3  \log 3 + 1.504 \cdot \sqrt{3} \pi \cdot \sum_{j\geq 1} \log(4d(\mathfrak{s}_j-1)).
\end{align}
Since $\log(4d(\mathfrak{s}_1)-1) \geq \log 19$ and
\begin{equation}
1.414 \cdot \frac{5}{4} \cdot 18 \cdot (1.504 \sqrt{3}  \pi + \frac{3 \log 3} {\log 19})\leq 296,
\end{equation}
we obtain the upper bound in Theorem \ref{thm3} also in this case.

Theorem 3 is proved. \hfill $\Box$


\begin{thebibliography}{99}





\bibitem{A1} {\sc Ahlfors}, Lars, Lecture on Quasiconformal
    Mappings,
    Van Nostrand, Princeton (1966).





\bibitem{Ar} {\sc Arnol'd}, Vladimir, On some Topological
    Invariants
    of Algebraic Functions,
{\it Trudy Moskov. Mat. Obsc.} {\bf 21} (1970), 27--46. Engl.
Transl.: {\it Trans. Moscow Math. Soc.} {\bf 21} (1970),
30--52.






\bibitem {BH} M. {Bestvina}, and M. {Handel}, \textit{Train Tracks for
    Surface
Homeomorphisms.}  Topology {\bf 34}, 1 (1995), 109--140.
\bibitem{Go} {\sc Goluzin}, Gennadii Michailovich, Geometric
    theory of
    functions of a complex variable, {\it Translations of
    Mathematical
    Monographs}, Vol. 26, American Mathematical Society,
    Providence,
    R.I. 1969.


\bibitem{Jo1} {\sc J\"oricke}, Burglind, Braids, Conformal
    Module and
Entropy (145 p.),  arXiv:1412.7000, to appear Lecture Notes in Math.
\bibitem{Jo2} {\sc J\"oricke}, Burglind, Braids, Conformal
    Module and
Entropy, {\em C.R.Acad.Sci.Paris}, Ser.I, {\bf 351} (2013)
289-293.
\bibitem{Jo3} {J\"oricke}, Burglind, Fundamental groups, slalom curves and extremal length, Operator Theory: Advances and Applications, {261}, Birkh\"auser Basel, 2018, 307-315.

\bibitem{Le} {\sc Lehner}, Joseph, Discontinuous Groups and Automorpic Functions, Mathematical Survays VIII, Am. Math. Soc., Providence, Rhode Island, 1964.

\bibitem {Ma} {\sc Markushevich}, Alexei Ivanovich,
    Theory of functions of a complex variable. Vol. III,
    Chelsea
    Publishing Co., 1977.

\end{thebibliography}
\end{document}